\documentclass[11pt, reqno]{amsart}

\usepackage[usenames,dvipsnames]{color}

\usepackage[OT2, T1]{fontenc}
\usepackage{url}
\usepackage{amsmath}
\usepackage{array}
\usepackage{graphicx}
\usepackage{amssymb}
\usepackage{amsthm}
\definecolor{link}{RGB}{11,0,128}
\usepackage[colorlinks=true,citecolor=NavyBlue,linkcolor=OliveGreen,urlcolor=link]{hyperref}
\usepackage{hyperref}
\usepackage{paralist}
\usepackage{colonequals}			
\usepackage{color}
\usepackage{mathabx}			
\usepackage{amscd}
\usepackage[all,cmtip]{xy}			
\usepackage{sseq}				
\usepackage{verbatim}
\usepackage{parskip}
\usepackage{microtype}
\usepackage[in]{fullpage}
\usepackage{mathrsfs} 			
\usepackage{cleveref}
\usepackage{enumitem}
\usepackage{moreenum}			
\usepackage[alphabetic,lite]{amsrefs} 	
\usepackage{stmaryrd}			
\usepackage[foot]{amsaddr}
\usepackage{subfiles}
\usepackage{ragged2e}
\usepackage{stackengine}

\DeclareSymbolFont{cyrletters}{OT2}{wncyr}{m}{n}
\DeclareMathSymbol{\Sha}{\mathalpha}{cyrletters}{"58}	

\newcommand{\ce}{\colonequals}

\newcommand{\llb}{\llbracket}			
\newcommand{\llp}{(\!(}			
	
\newcommand{\rrb}{\rrbracket}			
\newcommand{\rrp}{)\!)}			

\newcommand{\Zar}{\mathrm{Zar}}		

\newcommand{\csub[1]}{\centering\subsection{\text{#1} \nopunct} \hfill \justify}

\providecommand{\up}[1]{{\upshape(}#1{\upshape)}}
\providecommand{\uref}[1]{{\upshape\ref{#1}}}
\providecommand{\uS}{{\upshape\S}}

\renewcommand{\b}{\textbf}
\providecommand{\ucolon}{{\upshape:} }
\providecommand{\uscolon}{{\upshape;} }

\newcommand{\brems}{\begin{rems} \hfill \begin{enumerate}[label=\b{\thenumberingbase.},ref=\thenumberingbase]}

\newcommand{\erems}{\end{enumerate} \end{rems}}
\newcommand{\begs}{\begin{egs} \hfill \begin{enumerate}[label=\b{\thenumberingbase.},ref=\thenumberingbase]}

\newcommand{\eegs}{\end{enumerate} \end{egs}}
\newcommand{\m}{\item}
\newcommand{\bsm}{\begin{smallmatrix}}
\newcommand{\esm}{\end{smallmatrix}}
\newcommand{\blem}{\begin{lemma}}
\newcommand{\elem}{\end{lemma}}

\newcommand{\bconj}{\begin{conj}}
\newcommand{\econj}{\end{conj}}
\newcommand{\bprob}{\begin{Problem}}
\newcommand{\eprob}{\end{Problem}}

\newcommand{\bq}{\begin{Q}}
\newcommand{\eq}{\end{Q}}
\newcommand{\benum}{\begin{enumerate}[label={{\upshape(\alph*)}}]}
\newcommand{\benuma}{\begin{enumerate}[label={{\upshape(\arabic*)}}]}
\newcommand{\benumb}{\begin{enumerate}[label={{\upshape\b{\arabic*.}}}]}
\newcommand{\benumr}{\begin{enumerate}[label={{\upshape(\roman*)}}]}
\newcommand{\eenum}{\end{enumerate}}
\newcommand{\bitem}{\begin{itemize}}
\newcommand{\eitem}{\end{itemize}}
\newcommand{\bc}{}
\newcommand{\bd}{\begin{defn}}
\newcommand{\ed}{\end{defn}}

\newcommand{\beg}{\begin{eg}}
\newcommand{\eeg}{\end{eg}}

\newcommand{\bcl}{\begin{claim}}
\newcommand{\ecl}{\end{claim}}

\newcommand{\x}{\text}

\newcommand{\q}{\quad}

\newcommand{\qq}{\quad\quad}

\newcommand{\qqqq}{\quad\quad\quad\quad}

\newcommand{\tst}{\textstyle}
\newcommand{\ba}{\begin{aligned}}
\newcommand{\ea}{\end{aligned}}
\newcommand{\be}{\begin{equation}}
\newcommand{\ee}{\end{equation}}
\newcommand{\bpf}{\begin{proof}}
\newcommand{\epf}{\end{proof}}
\newcommand{\bthm}{\begin{thm}}
\newcommand{\ethm}{\end{thm}}
\newcommand{\bprop}{\begin{prop}}
\newcommand{\eprop}{\end{prop}}
\newcommand{\bcor}{\begin{cor}}
\newcommand{\ecor}{\end{cor}}
\newcommand{\brem}{\begin{rem}}
\newcommand{\erem}{\end{rem}}

\usepackage{aliascnt}
\newaliascnt{numberingbase}{subsubsection}
\numberwithin{equation}{numberingbase}

\newtheoremstyle{thms}{0.5em}{0.5em}{\itshape}{}{\bfseries}{.}{ }{}
\theoremstyle{thms}
\newtheorem{conj}[numberingbase]{Conjecture}
\newtheorem{corollary}[numberingbase]{Corollary}
\newtheorem{cor}[numberingbase]{Corollary}

\newtheorem{lemma}[numberingbase]{Lemma}

\newtheorem{prop}[numberingbase]{Proposition}

\newtheorem{proposition}[numberingbase]{Proposition}
\newtheorem{Q}[numberingbase]{Question}
\newtheorem{thm}[numberingbase]{Theorem}
\newtheorem{theorem}[numberingbase]{Theorem}

\newtheoremstyle{claims}{0.5em}{0.5em}{}{}{\itshape}{.}{ }{}
\theoremstyle{claims}
\newtheorem{claim}[equation]{Claim}

\newtheoremstyle{defs}{0.5em}{0.5em}{}{}{\bfseries}{.}{ }{}
\theoremstyle{defs}
\newtheorem{defn}[numberingbase]{Definition}
\newtheorem{definition}[numberingbase]{Definition}

\newtheorem{eg}[numberingbase]{Example}
\newtheorem*{egs}{Examples}
\newtheorem{rem}[numberingbase]{Remark}

\newtheorem{remark}[numberingbase]{Remark}
\newtheorem*{rems}{Remarks}

\Crefname{claim}{Claim}{Claims}
\Crefname{bclaim}{Claim}{Claims}
\Crefname{sublemma}{Lemma}{Lemmas}
\Crefname{conj}{Conjecture}{Conjectures}
\Crefname{cor}{Corollary}{Corollaries}
\Crefname{defn}{Definition}{Definitions}
\Crefname{eg}{Example}{Examples}
\Crefname{prop}{Proposition}{Propositions} 
\Crefname{Q}{Question}{Questions}
\Crefname{rem}{Remark}{Remarks}
\Crefname{thm}{Theorem}{Theorems}
\Crefname{Theorem}{Theorem}{Theorems}
\Crefname{variant}{Variant}{Variants}
\Crefname{caution}{Caution}{Cautions}

\theoremstyle{thms}
\newtheorem{thm-tweak}[subsection]{Theorem}
\Crefname{thm-tweak}{Theorem}{Theorems}
\newtheorem{lemma-tweak}[subsection]{Lemma}
\Crefname{lemma-tweak}{Lemma}{Lemmas}
\newtheorem{cor-tweak}[subsection]{Corollary}
\Crefname{cor-tweak}{Corollary}{Corollaries}
\newtheorem{prop-tweak}[subsection]{Proposition}
\Crefname{prop-tweak}{Proposition}{Propositions} 
\newtheorem{conj-tweak}[subsection]{Conjecture}
\Crefname{conj-tweak}{Conjecture}{Conjectures} 
\newtheorem{q-tweak}[subsection]{Question}
\Crefname{q-tweak}{Question}{Questions} 

\theoremstyle{defs}
\newtheorem{defn-tweak}[subsection]{Definition}
\Crefname{defn-tweak}{Definition}{Definitions}
\newtheorem{eg-tweak}[subsection]{Example}
\Crefname{eg-tweak}{Example}{Examples}
\newtheorem*{rems-tweak}{Remarks}
\newtheorem{rem-tweak}[subsection]{Remark}
\Crefname{rem-tweak}{Remark}{Remarks}

\newtheoremstyle{subsection-tweak}
   {2pt}
   {3pt}%
   {}
   {}%
   {\bfseries}
   {}%
   {.5em}
   {\thmnumber{\@{#1}{}\@{#2}.}%
    \thmnote{~{\bfseries#3.}}}    
    
\theoremstyle{subsection-tweak}
\newtheorem{pp}[numberingbase]{}
\newcommand{\bpp}{\begin{pp}}
\newcommand{\epp}{\end{pp}}
\newcommand{\bppt}{\begin{pp-tweak}}
\newcommand{\eppt}{\end{pp-tweak}}

\theoremstyle{subsection-tweak}
\newtheorem{pp-tweak}[subsection]{}

\makeatletter
\def\@tocline#1#2#3#4#5#6#7{
    \begingroup 
    \@ifempty{#4}{}{}

    \parindent\z@ \leftskip#3\relax \advance\leftskip\@tempdima\relax
    #5\hskip-\@tempdima
      \ifcase #1
       \or\or \hskip 2em \or \hskip 1em \else \hskip 3em \fi%
      #6\nobreak\relax
    \dotfill\hbox to\@pnumwidth{\@tocpagenum{#7}}\par
    \nobreak
    \endgroup
 }
 \def\l@section{\@tocline{1}{0pt}{1pc}{}{}}

\renewcommand{\tocsection}[3]{%
  \indentlabel{\@ifnotempty{#2}{\makebox[1.3em][l]{%
    \ignorespaces#1 \bfseries{#2}.\hfill}}}\bfseries{#3}
    \vspace{-5pt}}

\renewcommand{\tocsubsection}[3]{%
  \indentlabel{\@ifnotempty{#2}{\hspace*{-0.5em}\makebox[2.1em][l]{%
    \ignorespaces#1#2.\hfill}}}#3
    \vspace{-5pt}}
   
\makeatother 

\makeatletter 
\newcommand\appendix@section[1]{%
  \refstepcounter{section}%
  \orig@section*{Appendix \@Alph\c@section. #1}%
}
\let\orig@section\section
\g@addto@macro\appendix{\let\section\appendix@section}
\makeatother

\makeatletter
\@namedef{subjclassname@2020}{%
  \textup{2020} Mathematics Subject Classification}
\makeatother

\author{Alexis Bouthier}
\author{K\k{e}stutis \v{C}esnavi\v{c}ius}
\author{Federico Scavia}
\address{Institut de Math\'{e}matiques de Jussieu-PRG, Universit\'{e} Pierre et Marie Curie, 4 place Jussieu, 75005 Paris, France}
\address{CNRS, Universit\'{e} Paris-Saclay,   Laboratoire de math\'{e}matiques d'Orsay, F-91405, Orsay, France}
\address{CNRS, Institut Galil\'ee, Universit\'e Sorbonne Paris Nord, 99 avenue Jean-Baptiste
Cl\'ement, 93430, Villetaneuse, France}
\email{alexis.bouthier@imj-prg.fr, kestutis@math.u-psud.fr, scavia@math.univ-paris13.fr}

\date{\today}

\begin{document}


\title{Generically trivial torsors under constant groups}

\maketitle

\begin{abstract}
We resolve the Grothendieck--Serre question over an arbitrary base field~$k$: for a smooth $k$-group scheme~$G$ and a smooth $k$-variety~$X$, we show that every generically trivial $G$-torsor over~$X$ trivializes Zariski semilocally on~$X$. This was known when $G$ is reductive or when $k$ is perfect, and to settle it in general we uncover a wealth of new arithmetic phenomena over imperfect~$k$. 

We build our arguments on new purity theorems for torsors under pseudo-complete, pseudo-proper, and pseudo-finite $k$-groups, for instance, respectively, under wound unipotent $k$-groups, under pseudo-abelian varieties, and under the kernels $\mathrm{Ker}(i_G)$ of comparison maps $i_G$ that relate pseudo-reductive groups to restrictions of scalars of reductive groups. We then deduce an Auslander--Buchsbaum type extension theorem for torsors under quasi-reductive $k$-groups; for instance, we show that torsors over $\mathbb{A}^2_k \setminus \{(0,0)\}$ under wound unipotent $k$-groups extend to torsors over $\mathbb{A}^2_k$, in striking contrast to the case of split unipotent groups. For a quasi-reductive $k$-group $G$, this extension theorem allows us to quickly classify $G$-torsors over $\mathbb{P}^1_k$ by an argument that already simplifies the reductive case and to establish Birkhoff, Cartan, and Iwasawa decompositions for $G(k\llp t \rrp )$. 

We combine these new results with deep inputs from recent work on the structure of pseudo-reductive and quasi-reductive $k$-groups to show an unramifiedness statement for the Whitehead group (the unstable $K_1$-group) of a quasi-reductive $k$-group, and then use it to argue that, for a smooth $k$-group $G$ and a semilocal $k$-algebra $A$, every $G$-torsor over $\mathbb{P}^1_A$ trivial at $\{t = \infty\}$ is also trivial at $\{t = 0\}$, which is known to imply the Grothendieck--Serre conclusion via geometric arguments.  To achieve all this, we develop and heavily use the structure theory of $k$-group schemes locally of finite type.
 \end{abstract}

\vspace{-45pt}

\hypersetup{
    linktoc=page,     
}
\renewcommand*\contentsname{}
\q\\
\tableofcontents

\section{The Grothendieck--Serre question}

\csub[The main result] \label{sec:main-result}

In 1958, Grothendieck~\cite{Gro58}*{pages 26--27, remarques 3${}^\circ$} and Serre~\cite{Ser58b}*{page 31, remarque} predicted that, for an algebraically closed field~$k$, a finite type, smooth $k$-group scheme~$G$, and a smooth $k$-scheme~$X$, every generically trivial $G$-torsor over~$X$ trivializes \emph{Zariski} locally on~$X$. In~\cite{CTO92}, Colliot-Th\'el\`ene and Ojanguren proved that this is the case, and they also established several special cases of the analogous prediction over any base field~$k$.

We fully resolve the Grothendieck--Serre question over an arbitrary base field~$k$, with the main novelty being in the case of an imperfect~$k$ with a nonreductive~$G$. Over a general (imperfect)~$k$, group schemes locally of finite type have a rich structure, in which wound unipotent groups, pseudo-reductive and quasi-reductive groups, and pseudo-abelian varieties play major roles, see~\S\ref{sec:group-structure} for a review of this and for these terms. Consequently, the Grothendieck--Serre question over an imperfect~$k$ is vastly more intricate. On the other hand, allowing general~$G$ is important: for instance, recent ``inverse Galois'' type results of~\cite{Flo23}, \cite{BS24} say that essentially any~$G$ occurs as the automorphism group of a projective $k$-variety~$Y$, and then $G$-torsors amount to forms of~$Y$.

We work directly with geometrically regular, semilocal $k$-algebras: in the special case of local rings of smooth varieties~$X$, this recovers the setting above. Allowing semilocal rings is an additional complication because \Cref{eg:semilocal} rules out direct reductions to the local case.

\begin{theorem}[\Cref{thm:GS-main}]\label{thm:main-GS}
Let $k$ be a field, let $R$ be a geometrically regular, semilocal $k$-algebra with $K := \mathrm{Frac}(R)$, and let $G$ be a smooth $k$-group scheme \up{more generally, a locally of finite type $k$-group scheme such that every $\overline{k}$-torus of $G_{\overline{k}}$ lies in $(G^{\mathrm{gred}})_{\overline{k}}$, where $G^{\mathrm{gred}} \leq G$ is the largest smooth $k$-subgroup}. Every generically trivial $G$-torsor over~$R$ is trivial, that is,
\begin{equation}\label{eqn:GS-inj}
\mathrm{Ker}\big(H^1(R, G) \to H^1(K, G)\big) = \{*\}.
\end{equation}
\end{theorem}

The parenthetical condition on~$G$ was introduced by Gabber in~\cite{Gab12} 
for a different purpose, and it holds if~$k$ is perfect or if~$G$ is locally of finite type and either is a normal $k$-subgroup of a smooth $k$-group or is nilpotent, see Remark~\ref{rem:condition-star}. These cases cover many groups that appear ``in nature,'' for instance, Picard group schemes, which are always commutative, even though in general the parenthetical condition seems delicate to check.

\Cref{thm:main-GS} is optimal in that its parenthetical condition on~$G$ cannot be dropped and its conclusion~\eqref{eqn:GS-inj} cannot be strengthened to injectivity; relatedly, the parenthetical condition is lost by inner forms, as examples of Florence--Gille~\cite{FG21}*{Examples 7.2, Remark 7.3} show. In~\Cref{eg:condition-star}, we give a further example of this with~$G$ being the automorphism group of a pseudo-semisimple $k$-group and illustrate how to use \Cref{thm:main-GS} to show that this automorphism group is not a normal subgroup of any smooth $k$-group.

\bpp[Known cases]
\Cref{thm:main-GS} has so far been established in the following cases.
\begin{itemize}
    \item When $k$ is perfect by Colliot-Th\'el\`ene--Ojanguren \cite{CTO92} if $k$ is infinite and by Gabber (unpublished) if $k$ is finite. These works make additional assumptions on the $k$-group $G$ locally of finite type that are, however, not difficult to remove because their $k$ is perfect.
    
    \item When $G$ is reductive by Raghunathan \cite{Rag94}, \cite{Rag94e} if $k$ is infinite and by Gabber (unpublished) if $k$ is finite. The reductive case has subsequently been taken much further, culminating in the equal characteristic setting in the works of Fedorov--Panin \cite{FP15} if $k$ is infinite and of Panin \cite{Pan20a} if $k$ is finite, where they showed that \eqref{eqn:GS-inj} also holds if the reductive $G$ is defined merely over $R$ and need not descend to $k$. In \Cref{eg:fiberwise-unipotent,eg:congruence}, we give new examples showing the failure of such a generalization beyond reductive $G$.
    
    \item When $G$ is affine and $R \cong k\llb t\rrb$ by Florence--Gille in \cite{FG21}*{Theorem 6.3}.
\end{itemize}
\epp

\Cref{thm:main-GS} is genuinely new over imperfect $k$ and beyond reductive $G$, although Gabber also considered this direction in unpublished work using methods different from ours (private communication).
\footnote{As far as we are aware, Gabber's approach was based on establishing the Birkhoff decomposition \Cref{thm:BC-ann}~\ref{m:BCA-a} for quasi-reductive $k$-groups $G$ (compare with the approach of \cite{Rag94} in the reductive case) by building a suitable generalization of a refined Tits system on $G(k\llp t\rrp )$ in the sense of Kac--Peterson \cite{KP85}.} 
In fact, our \Cref{thm:main-GS} is the culmination of a succession of intermediate results of independent interest about purity and torsors under smooth (and often quasi-reductive) $k$-group schemes overviewed in \S\S\ref{sec:purity-ann}--\ref{sec:Whitehead-ann} below.

We recall that a $k$-group scheme $G$ is \emph{quasi-reductive} if it is connected, smooth, affine, and has no nontrivial split unipotent normal $k$-subgroups. If $k$ is perfect, then this is nothing else than being reductive, but it is much more general otherwise: as special cases, quasi-reductive groups include both smooth, connected, wound unipotent groups, which are abundant over imperfect $k$ due to the additivity of $p$-polynomials (example: $\{x = x^p + ty^p\} \le \mathbb{G}_{a,\, \mathbb{F}_p(t)}^2$), and pseudo-reductive groups (example: $\mathrm{Res}_{\mathbb{F}_p(t^{1/p})/\mathbb{F}_p(t)}(\mathrm{GL}_{n,\, \mathbb{F}_p(t^{1/p})})$). Already in the case of pseudo-reductive groups, the structure theory is vast and complex, although we have the enormous benefit of having the recent books \cite{CGP15}, \cite{CP16}, as well as the survey \cite{CP17}, where this theory developed by Tits and Conrad--Gabber--Prasad is presented. The classification of pseudo-reductive groups is not the key to the Grothendieck--Serre problem for their torsors---much like in the reductive case, where the usefulness of the classification in terms of root data only helps to understand the overall landscape---still, we apply and refine a significant number of results from these works while arguing \Cref{thm:main-GS}.

We stress that our attention to quasi-reductive groups or, for that matter, to wound unipotent groups, pseudo-reductive groups, or pseudo-abelian varieties, all of which play important roles in the proof of \Cref{thm:main-GS}, is not dictated by cravings for baroque generalizations but rather by sober realities of the situation. Indeed, all of these groups are subquotients of the fundamental filtration describing the structure of a general $k$-group scheme locally of finite type (see \S\ref{pp:fundamental-filtration}), so we must handle them to obtain \Cref{thm:main-GS}. The reason why quasi-reductive groups and not, for instance, pseudo-abelian varieties form the core case of the Grothendieck--Serre problem is that pseudo-abelian varieties $G$ over $k$ satisfy more: for them, not only is every generically trivial $G$-torsor $E$ over $R$ trivial, but also trivializations extend uniquely, that is, $E(R) \cong E(K)$ (see \Cref{thm:pseudo-extend}~\ref{m:PE-ii}), as was observed in the abelian variety case already by Serre himself \cite{Ser58b}*{p.~22, preuve du lemme~4}. This allows for stronger d\'evissage in exact sequences, so we may ``peel off'' the pseudo-abelian variety part of the fundamental filtration when proving \Cref{thm:main-GS}. Likewise, since $\mathbb{G}_a$ has no nontrivial torsors over affine schemes, we may also ``peel off'' the split unipotent part. What is left is a quasi-reductive group, for which further direct reductions of the Grothendieck--Serre problem appear delicate. This matches experience with the reductive case, in which reducing to semisimple or simply connected groups is both desirable and remarkably complex: such a reduction was the main goal of \cite{Pan20b}.

\csub[Purity and extension theorems for torsors] \label{sec:purity-ann}

We build \Cref{thm:main-GS} on the purity \Cref{thm:purity-ann} for torsors under groups that satisfy the following generalizations of finiteness or properness.

\bd
A finite type, separated scheme $X$ over a field $k$ is
\benumr
    \item (\Cref{def:pseudo-finite}) \emph{pseudo-finite} if $X(k^s)$ is finite;
    \item (\Cref{def:pseudo-proper}) \emph{pseudo-proper} (resp.,~\emph{pseudo-complete}) if it satisfies the valuative criterion of properness with respect to those discrete valuation rings over $k$ that are geometrically regular (resp.,~whose residue field is separable over $k$).
\eenum
\ed

As an example, a pseudo-proper $X$ is required to satisfy the valuative criterion of properness with respect to local rings of smooth curves over $k$, but not necessarily with respect to local rings of smooth curves over purely inseparable field extensions of $k$. Restrictions of scalars of proper schemes along such extensions are always pseudo-proper (granted that they are schemes and not merely algebraic spaces), but in most cases they are not proper. If $k$ is perfect, then an $X$ is pseudo-finite (resp.,~pseudo-proper or pseudo-complete) if and only if it is finite (resp.,~proper), but in general only the `if' holds and we have strict implications
\[
\text{pseudo-finite} \Rightarrow \text{pseudo-proper} \Rightarrow \text{pseudo-complete}.
\]
We show in \Cref{prop:PP-groups} that a finite type $k$-group scheme $G$ is pseudo-finite precisely when its largest connected, smooth $k$-subgroup is trivial, and that $G$ is pseudo-proper (resp.,~pseudo-complete) precisely when its largest connected, smooth, affine $k$-subgroup is strongly wound unipotent (resp.,~is wound unipotent) in the sense that it is unipotent and has no nontrivial unirational $k$-subgroups. For instance, pseudo-abelian varieties in the sense of Totaro \cite{Tot13} are pseudo-proper, and they are proper precisely when they are abelian varieties.

Pseudo-properness is a robust geometric refinement of the older notion of pseudo-completeness: for instance, conditionally on resolving singularities, a $k$-smooth, integral $X$ is pseudo-proper if and only if for every proper, regular compactification $X \subset \overline{X}$ over $k$ with $\overline{X} \setminus X$ a divisor, the $k$-smooth locus of $\overline{X}$ is precisely $X$. As for pseudo-completeness, Borel--Tits \cite{BT78}*{Proposition~1}, Tits \cite{Tit13}*{cours 1992--1993, section~2.5}, and Conrad--Gabber--Prasad \cite{CGP15}*{Proposition C.1.6} showed that $G/P$ is pseudo-complete for every pseudo-parabolic $k$-subgroup $P$ of a connected, smooth, affine $k$-group $G$. In \Cref{thm:pseudo-flag}, we show (via a new argument) that such $G/P$ are even pseudo-proper.

\begin{theorem}\label{thm:purity-ann}
Let $k$ be a field, let $S$ be a geometrically regular $k$-scheme, let $Z \subset S$ be a closed subset of codimension $\ge 2$, and let $G$ be a finite type $k$-group scheme. If either
\benumr
    \item\label{m:PA-i} \up{\Cref{thm:purity-pseudo-finite}}. $G$ is pseudo-finite and commutative\uscolon or
    \item \up{\Cref{thm:purity-pproper}}. $G$ is pseudo-proper and smooth\uscolon or
    \item\label{m:PA-iii} \up{\Cref{thm:purity-pcomplete}}. $G$ is pseudo-complete and smooth, and every $z \in Z$ of codimension $2$ in $S$ lies in a geometrically regular $k$-subscheme $S_z \subset S$ of codimension $> 0$ \up{when $k_z/k$ is separable, we may take $S_z = z$}\uscolon
\eenum
then pullback induces an equivalence of categories
\[
\{G\text{-torsors over }S\} \xrightarrow{\sim} \{G\text{-torsors over }S \setminus Z\}.
\]
\end{theorem}

For instance, \Cref{thm:purity-ann}~\ref{m:PA-iii} says that for a smooth, wound unipotent $k$-group $G$, every $G$-torsor over $\mathbb{A}^2_k \setminus \{ (0, 0)\}$ extends uniquely to a $G$-torsor over $\mathbb{A}^2_k$. This came as a surprise because wound groups tend to have many nontrivial torsors (see \cite{Ros25}*{Theorem 1.6}), whereas $\mathbb{G}_a$ has many nontrivial torsors over $\mathbb{A}^2_k \setminus \{ (0, 0)\}$ none of which extends (which is how one sees, via \v{C}ech cohomology, that $\mathbb{A}^2_k \setminus \{ (0, 0)\}$ is not affine).

Thus, for the extendability of torsors, wound unipotent groups are much closer to reductive groups than to split unipotent groups. For instance, reductive group torsors over $\mathbb{A}^2_k \setminus \{ (0, 0)\}$ extend to those over $\mathbb{A}^2_k$, and likewise over all regular schemes of dimension $2$, thanks to the Auslander--Buchsbaum formula, which gives the key case of vector bundles (see \cite{CTS79}*{Corollary 6.13}, or \cite{torsors-regular}*{Section~1.3.9} for a review). We use the purity \Cref{thm:purity-ann} to generalize this extendability result to quasi-reductive groups as follows.

\begin{theorem}[\Cref{thm:AB-extn}]\label{thm:AB-ann}
Let $k$ be a field, let $S$ be a geometrically regular $k$-scheme of dimension $2$, let $z \in S$ be a point of height $2$, and let $G$ be a quasi-reductive $k$-group. Suppose that either $G$ is pseudo-reductive or $z$ lies on a geometrically regular $k$-subscheme $S_z \subset S$ of codimension $> 0$. Pullback induces an equivalence of categories
\[
\{G\text{-torsors over }S\} \xrightarrow{\sim} \{G\text{-torsors over }S \setminus z\}.
\]
\end{theorem}

\Cref{thm:AB-ann} directly reduces to vector bundles only in the reductive case: by the Matsushima theorem \cite{Alp14}*{Theorems 9.4.1 and 9.7.5}, a connected, smooth subgroup $G \le \mathrm{GL}_n$ is reductive if and only if the homogeneous space $\mathrm{GL}_n/G$ is affine, and this affineness is critical because it implies that reductions of $\mathrm{GL}_n$-torsors to $G$-torsors over $S \setminus z$ extend uniquely to those over $S$. Moreover, \Cref{thm:AB-ann} is specific to dimension $2$, as nontrivial vector bundles exist already over the punctured spectrum of $\mathbb{C}\rrb s, t, u\rrb$: indeed, the kernel of the map $\mathbb{C}\llb s, t, u\rrb^{\oplus 3} \xrightarrow{(s,\, t,\, u)} \mathbb{C}\llb s, t, u\rrb$ is such by the Auslander--Buchsbaum formula applied to the cokernel. We argue \Cref{thm:purity-ann,thm:AB-ann} by extensive d\'evissage both in $S$ and in $G$ that uses the ``classical'' cases of these theorems (finite groups, abelian varieties, or reductive groups) to reduce to when $S$ is (roughly) $\mathrm{Spec}(k\llb s, t\rrb)$ and $G$ is wound unipotent, given by the vanishing of some $p$-polynomial $F$, and then by computing with this $F$ (woundness amounts to a nonvanishing property of the principal part of $F$, which is critical in the computation). The main inputs to the d\'evissage in $S$ are the Gabber--Quillen geometric presentation theorem and the Popescu theorem. The d\'evissage in $G$ is more intricate, for instance, it requires a fresh common perspective on the theories of pseudo-reductive groups and pseudo-abelian varieties. Our point of view on both is that it is fruitful to study them via the comparison map
\[
i_G \colon G \rightarrow \mathrm{Res}_{k'/k}(\overline{G})
\]
where $k'/k$ is the finite, purely inseparable field of definition of the geometric unipotent radical $\mathscr R_{\mathrm{u},\, \overline{k}}(G_{\overline{k}})$ and $\overline{G} := G_{k^\prime}/\mathscr{R}_{\mathrm u,\, k'}(G_{k'})$ is the associated reductive $k^\prime$-group (resp.,~an abelian variety over $k^\prime$). Except perhaps for the pseudo-abelian variety aspect, this is not new: the map $i_G$ also appeared in \cite{CGP15}, then to a much larger extent in \cite{CP16} and in \cite{CP17}. What is new is the affineness of the homogeneous space $\mathrm{Res}_{k'/k}(\overline{G})/i_G(G)$ that we show in \Cref{prop:Brauer-affine} by building on ideas from \cite{brauer-purity}*{Lemma 2.1}. As we already saw when discussing the Matsushima theorem above, affineness of homogeneous spaces is both delicate and key for handling torsors. As for the kernel $\mathrm{Ker}(i_G)$, it is unipotent, pseudo-finite, and, in situations to which it is easy to reduce to, also commutative (see \S\ref{pp:pav-iG} and \S\ref{pp:pred-iG}), so \Cref{thm:purity-ann}~\ref{m:PA-i} applies to it. This control of the kernel and the ``cokernel'' translates into the control of torsors when passing from $G$ to $i_G(G)$, then to $\mathrm{Res}_{k'/k}(\overline{G})$, and, finally, to the ``classical'' case of $\overline{G}$, so it enables our d\'evissage in $G$. This sequence of reductions from $G$ to $\overline{G}$ is also how we argue the aforementioned pseudo-properness of $G/P$ in \Cref{thm:pseudo-flag}. Another argument for the latter is to note that \Cref{thm:AB-ann} implies that the affine Grassmannian of a quasi-reductive $k$-group is ind-pseudo-proper and to then realize $(G/P)_{k^s}$ as a closed subscheme of some such affine Grassmannian. For the sake of focus, we do not include this alternative approach but we hope to return to it, especially since ind-pseudo-properness combined with the loop rotation action also gives different proofs for the Birkhoff and Cartan decompositions of \Cref{thm:BC-ann} below.

\csub[Classification of $G$-torsors over $\mathbb{P}^1_k$] \label{sec:P1k-ann}

For a field $k$, torsors over $\mathbb{P}^1_k$ under a reductive $k$-group $G$ form a well-studied subject, with key classification results of Grothendieck \cite{Gro57}, Harder \cite{Har68}, and Biswas--Nagaraj \cite{BN09}, and subsequent simplifications of Ansch\"utz \cite{Ans18} and Wedhorn \cite{Wed24}, among others. We extend this classification to when $G$ is merely quasi-reductive, and simultaneously quickly reprove the reductive case by combining \Cref{thm:AB-ann} with results of Wedhorn \cite{Wed24}, alternatively, with the theory of coherent completeness of algebraic stacks from \cite{AHR20}. Our approach is more robust already in the reductive case, for instance, we do not need our smooth $k$-group to be connected or even affine.

\bthm[\Cref{thm:torsors-P1}] \label{thm:P1k-ann} For a field $k$ and a smooth $k$-group scheme $G$ whose largest connected, smooth, affine $k$-subgroup $G^{\mathrm{sm},\, \mathrm{lin}} \le G$ is quasi-reductive,
\[
H^1(\mathbb P^1_k, G) \cong H^1(\mathbf B\mathbb G_{m}, G) \quad \text{and} \quad H^1_{\Zar}(\mathbb{P}^1_k,G) \cong \mathrm{Hom}_{k\text{-\upshape{gp}}}(\mathbb{G}_m, G)/G(k),
\]
moreover, a $G$-torsor $E$ over $\mathbb P^1_k$ is Zariski locally trivial if and only if it is trivial at a single $k$-point, in which case it reduces to the $\mathbb{G}_m$-torsor $\mathscr{O}(1)$ along some $k$-homomorphism $\lambda\colon \mathbb{G}_{m,\, k} \rightarrow G$.
\ethm

The key idea is to view $\mathbb P^1_k$ as $[(\mathbb{A}^2_k\setminus\{0\})/\mathbb{G}_m]$, that is, as the open complement of the stacky origin $\mathbf B \mathbb G_{m} \subset [\mathbb A^2_k/\mathbb G_m]$, and to then apply \Cref{thm:AB-ann} (in its finer form given in \Cref{thm:AB-extn}) to uniquely extend $G$-torsors over $\mathbb P^1_k$ to those over $[\mathbb A^2_k/\mathbb G_m]$. We then classify $G$-torsors over $[\mathbb A^2_k/\mathbb G_m]$ for any smooth $k$-group $G$ by reducing to results from \cite{Wed24} that are specific to the case when $G$ is also affine, see \Cref{lem:maps-full}.

In contrast, it seems difficult to directly adapt known arguments from the reductive case because they critically use the reductivity, for instance, the Tannakian approach of Ansch\"utz \cite{Ans18} rests on Haboush's theorem from \cite{Hab75}. Relatedly, the classification \Cref{thm:P1k-ann} fails for general connected, smooth, affine $k$-groups $G$, see \Cref{rem:classification-fail}.

\csub[The Birkhoff, Cartan, and Iwasawa decompositions] \label{sec:Birkhoff-ann}

\Cref{thm:P1k-ann} and its proof method allow us to quickly establish the Birkhoff and the Cartan decompositions of a quasi-reductive $k$-group, once more even after dropping the connectedness or the even affineness assumptions. Our argument for the Birkhoff decomposition is new and simpler already in the reductive case, it is inspired by the Alper--Heinloth--Halpern-Leistner approach to the reductive case of the Cartan decomposition from \cite{AHHL21}.

\bthm \label{thm:BC-ann} Let $k$ be a field, let $G$ be a $k$-group scheme $G$ locally of finite type, and let $G^{\mathrm{sm},\, \mathrm{lin}} \le G$ be the largest connected, smooth, affine $k$-subgroup \up{see \uS\uref{pp:fundamental-filtration}~\ref{m:Gsmlin}}.

\benum
    \item \label{m:BCA-a} \up{Birkhoff decomposition, Theorem \uref{thm:Birkhoff}}. If $G^{\mathrm{sm},\, \mathrm{lin}}$ is quasi-reductive, then
    \[
    \textstyle G(k\llp t\rrp) = \coprod_{\lambda \in \mathrm{Hom}_{k\text{-\upshape{gp}}}(\mathbb{G}_m,\, G)/G(k)} G(k[t^{-1}]) t^\lambda G(k\llb t\rrb);
    \]
    \item \up{Cartan decomposition, Theorem \uref{thm:Cartan}}. If $G^{\mathrm{sm},\, \mathrm{lin}}$ is quasi-reductive, then
    \[
    \textstyle G(k\llp t\rrp) = \coprod_{\lambda \in \mathrm{Hom}_{k\text{-\upshape{gp}}}(\mathbb{G}_m,\, G)/G(k)} G(k\llb t\rrb) t^\lambda G(k\llb t\rrb);
    \]
    \item \up{Iwasawa decomposition, Theorem \uref{thm:Iwasawa}}. For each pseudo-parabolic $k$-subgroup $P \le G^{\mathrm{sm},\, \mathrm{lin}}$, 
    \[
    \textstyle G(k \llp t\rrp) = P(k\llp t\rrp) G(k\llb t\rrb ).
    \]
\end{enumerate}
\ethm

The Iwasawa decomposition results by combining the pseudo-properness of (the connected components of) $G/P$ with a special case of \Cref{thm:GS-main}. In the Birkhoff and Cartan decompositions, for a maximal split $k$-torus $S \le G$, the indexing set $\mathrm{Hom}_{k\text{-gp}}(\mathbb{G}_m,\, G)/G(k)$ may be identified with $\mathrm{Hom}_{k\text{-gp}}(\mathbb{G}_m,\, S)/N_G(S)(k)$, see \Cref{lem:SGconj}. Thus, if $G$ contains no nontrivial split $k$-tori, then the Cartan decomposition gives $G(k\llp t\rrp ) = G(k\llb t\rrb)$, see \Cref{cor:anisotropic}. This consequence leads to simple new proofs of some key results of \cite{CGP15}*{Appendix C}: in \Cref{cor:Ga-Gm,cor:split-parabolic}, we show that the maximal split unipotent $k$-subgroups are precisely the unipotent radicals of the minimal pseudo-parabolic $k$-subgroups (a generalization of a result of Borel--Tits from \cite{BT71}), and that if a quasi-reductive $k$-group has $\mathbb{G}_{a,\, k}$ as a $k$-subgroup, then it also has $\mathbb G_{m,\, k}$ as a $k$-subgroup. 

To link the Birkhoff decomposition to \Cref{thm:P1k-ann} it suffices to note that, by patching, the set of double cosets  $G(k[t^{-1}]) \backslash G(k\llp t\rrp)/G(k\llb t\rrb )$ is identified with the set of isomorphism classes of those $G$-torsors over $\mathbb P^1_k$ that trivialize over both $\mathbb P^1_k \setminus \{ t = 0\}$ and $\{t = 0\}$, so with the set of isomorphism classes of Zariski locally trivial $G$-torsors over $\mathbb P^1_k$. Similarly, the set $G(k\llb t\rrb) \backslash G(k\llp t\rrp )/G(k\llb t\rrb)$ is identified with the set of isomorphism classes of those $G$-torsors over the glueing $\mathrm{Spec}(k\llb t\rrb) \bigcup_{\mathrm{Spec}(k\llp t\rrp )} \mathrm{Spec}(k\llb t\rrb )$ that trivialize over both copies of $\mathrm{Spec}(k\llb t\rrb)$. After noting that this glueing is the open complement of $\{s = s^\prime = 0\}$ in the quotient stack $[\mathrm{Spec}(k\llb t\rrb[s,s^{\prime}]/(ss^{\prime}-t))/\mathbb G_m]$, where $\mathbb{G}_m$ acts over $k\llb t\rrb$ by scaling $s$ (resp.,~$s^{\prime}$) via the character of weight $1$ (resp.,~$-1$), we apply the Auslander--Buchsbaum extension \Cref{thm:AB-ann} to extend them to $G$-torsors over the entire quotient stack. We then classify $G$-torsors over the latter by once more applying the results of Wedhorn \cite{Wed24} or by adapting the deformation-theoretic method of \cite{AHHL21} from the reductive case. 

In the case when $k$ is finite but the quasi-reductive group $G$ is defined merely over $k\llp t\rrp $, Solleveld established Cartan decompositions in \cite{Sol18}*{Theorem 5} using methods from Bruhat--Tits theory. Although the finiteness of $k$ is very restrictive, he manages to decompose $G(k\llp t\rrp )$ with respect to double cosets of more general subgroups than our $G(k\llb t\rrb)$. It would be very interesting to find a way to adapt our geometric approach to these more general subgroups and to obtain the general Cartan decompositions without restrictions on $k$.

\csub[Unramifiedness of the Whitehead group and torsors over $\mathbb P^1_A$] \label{sec:Whitehead-ann}

Geometric simplifications of our geometrically regular, semilocal $k$-algebra $R$, more precisely, the geometric approach to the Grothendieck--Serre conjecture developed in \cite{CTO92}, \cite{FP15}, \cite{Pan20a}, \cite{split-unramified}, and \cite{totally-isotropic}, reduce our goal \Cref{thm:main-GS} to the study of $G$-torsors over $\mathbb P^1_R$, more precisely, to arguing that a $G$-torsor over $\mathbb P^1_R$ that is trivial at $\{t = \infty\}$ is also trivial at $\{t = 0\}$. In the toy case when $R$ is a field, this sectionwise triviality follows from the classification \Cref{thm:P1k-ann}, so the problem becomes that of bootstrapping this statement from the residue fields of $R$. For reductive groups, this was carried out in \cite{totally-isotropic}*{Theorem 3.5} (see also \cite{PS25}) using the geometry of $\mathrm{Bun}_G$. In our setting, a relevant extension of \Cref{thm:P1k-ann} is the following theorem.

\begin{theorem}[\Cref{thm:swise}] \label{thm:swise-ann} 
For a field $k$, a smooth $k$-group $G$, a semilocal $k$-algebra $A$, and a $G$-torsor $E$ over $\mathbb{P}^1_A$, if $E|_{\{t = \infty\}}$ is trivial, then so is $E|_{\{ t = 0\}}$.
\end{theorem}

\Cref{thm:swise-ann} is the most technically demanding part of the proof of \Cref{thm:main-GS} and rests on some of the deepest aspects of the structure theory of pseudo-reductive and quasi-reductive $k$-groups that were recently developed in \cite{CGP15} and \cite{CP16}. For instance, after reducing to a $G$ that is quasi-semisimple, that is, quasi-reductive and perfect, we critically rely on the theory of the simply connected cover of a quasi-semisimple $k$-group supplied by \cite{CP16}*{Theorem 5.1.3}, as well as on inputs from \cite{CGP15}*{Appendix C} about the existence of Levi subgroups of quasi-reductive groups and about the subgroup $G(k)^+ \le G(k)$ generated by the ``elementary matrices'' (by the $k$-points of the unipotent radicals of the pseudo-parabolic $k$-subgroups). Indeed, as in the reductive case, a crucial step towards \Cref{thm:swise-ann} is the so-called unramifiedness of the Whitehead group
$$
W(k, G) \colonequals G(k)/G(k)^+
$$
that we argue in a sufficient for our purposes pseudo-split case in \Cref{prop:W-unr}.

The Whitehead group is an invariant of $K$-theoretic flavor, for instance, the stabilization (the direct limit over $n$) of the Whitehead groups  $W(k, \mathrm{GL}_n)$ is $K_1(k)$. The strategy for bootstrapping \Cref{thm:swise-ann} from the case of the residue fields of $A$ supplied by \Cref{thm:P1k-ann} is to modify $E$ along a well-chosen $A$-(finite \'etale) closed $Y \subset \mathbb G_{m,\, A}$ in order to force $E$ to be residually trivial over $A$ and to then use the rigidity of $G$-torsors over $\mathbb P^1$ that results from deformation theory. Thus, restricting for the sake of illustration to when $A$ is local with residue field $k$ and $Y$ is an $A$-point cut out by a $\tau \colonequals t - y$, the relevance of the Whitehead group $W(k\llp\tau\rrp, G)$ stems from the fact that $G(k\llp\tau\rrp)$ parametrizes patchings of torsors along the formal completion of $\{\tau = 0\}$ in $\mathbb P^1_k$, while $G(k\llp \tau\rrp)^+$ parametrizes ``elementary'' patchings, which are straightforward to lift to patchings along the formal completion of $\{\tau = 0\}$ in $\mathbb P^1_A$. The problem of controlling the difference between general patchings and the ``elementary'' ones becomes the problem of controlling the Whitehead group, and the unramifiedness of the latter, expressed concretely as
$$
G(k\llp\tau\rrp) = G(k\llp\tau\rrp)^+ G(k\llb\tau\rrb),
$$
becomes the key to the liftability of the relevant patchings to $\mathbb P^1_A$, so also to the bootstrap argument.

\bppt[Fixed base field]\label{pp:notation}
Throughout this article, we fix an arbitrary base field $k$. The main case to keep in mind is when $k$ is imperfect because that of a perfect $k$ is much simpler. 
\eppt

\bppt[Notation and conventions] \label{pp:conv} 
For a field $K$, we let $K^s$ (resp.,~$\overline{K}$) denote a choice of its separable closure (resp.,~denote the algebraic closure of $K^s$). We let $\alpha_p$ (resp.,~$\mu_p$) denote the kernel of the endomorphism of $\mathbb{G}_{a,\, K}$ (resp., $\mathbb{G}_{m,\, K}$) given by $t\mapsto t^p$, where $p$ is the characteristic exponent of $K$ (so $p = 1$ if $\mathrm{char}(K) = 0$, and else $p = \mathrm{char}(K)$). A $K$-algebra $R$ is \emph{geometrically regular} if $R \otimes_K K^\prime$ is a regular Noetherian ring for every finite field extension $K^\prime/K$; by Popescu's theorem \cite{SP}*{Theorem~\href{https://stacks.math.columbia.edu/tag/07GC}{07GC}} (a deep result!), this amounts to $R$ being Noetherian and a filtered direct limit of smooth $k$-algebras. We let $\mathrm{Frac}(-)$ denote the total ring of fractions. For a ring $A$, we write $A\{t\}$ for the Henselization of $A[t]$ along $tA[t]$.

We freely use various widely-known properties of restriction of scalars reviewed in \cite{BLR90}*{Section~7.6} and in \cite{CGP15}*{Section A.5}, in particular, its left exactness and commutation with quotients by faithful actions of smooth groups (which is in essence immediate from \cite{SP}*{Lemma~\href{https://stacks.math.columbia.edu/tag/04GH}{04GH}}, see also \cite{CGP15}*{Corollary~A.5.4~(3)}). We will also freely use the representability by algebraic spaces of restrictions of scalars of algebraic spaces along finite flat maps, see \cite{SP}*{Proposition~\href{https://stacks.math.columbia.edu/tag/05YF}{05YF}}. In general, analogous representability for schemes needs the quasi-compact opens to all be quasi-projective, but when dealing with $k$-group schemes locally of finite type this will not be an issue because the results that we review in \S\ref{pp:k-groups-lft} supply such quasi-projectivity.

We recall from \cite{SGA3II}*{expos\'e XVII, d\'efinition 1.1, propositions 1.2, th\'eor\`eme 3.5, lemme 3.9} that a $k$-group scheme $U$ is unipotent if $U_{\overline{k}}$ is a finite successive extension of closed $\overline{k}$-subgroups of $\mathbb{G}_{a,\, \overline{k}}$ (which may all be taken to be $\mathbb{G}_{a,\, \overline{k}}$ if $U$ is smooth and connected), equivalently, if $U$ is a closed $k$-subgroup of the group of upper unitriangular matrices of some $\mathrm{GL}_{n,\, k}$. We recall from \cite{SGA3II}*{expos\'e~XVII, propositions 2.1, 2.2} that unipotent groups are affine, of finite type, and stable under closed subgroups, quotients (which are therefore affine, see also \Cref{lem:unip-quotient}), extensions, base change, and, by \cite{CGP15}*{Proposition A.5.12}, \cite{SGA3Inew}*{expos\'e VII${}_{\x{\upshape{A}}}$, proposition~8.3}, and embedding a large Frobenius kernel into $\mathrm{GL}_{n,\, k}$, also under restrictions of scalars along field~extensions.

We do not assume algebraic spaces to be quasi-separated, that is, we use \cite{SP}*{Definition~\href{https://stacks.math.columbia.edu/tag/025Y}{025Y}}. For a group fppf sheaf $G$ over a scheme $S$, as in \cite{Ray70b}*{chapitre VI, d\'efinitions VI 1.1}, a \emph{homogeneous space} (resp.,~a \emph{torsor}) under $G$ is an fppf $S$-sheaf $E$ that fppf locally on $S$ has a section and is equipped with a right $G$-action such that the map $G \times_S E \xrightarrow{(g,\, e) \mapsto (eg,\, e)} E \times_S E$ is an fppf cover (resp.,~is an isomorphism), in particular, throughout we work with torsors for the fppf topology. We freely use well-known representability properties of torsors and of quotients reviewed in \cite{torsors-regular}*{Section 1.2.3}.
\eppt

\subsection*{Acknowledgements}  
We thank Ofer Gabber for many helpful interactions; as the reader will notice, this article owes a significant intellectual debt to his ideas. We thank Michel Brion, Brian Conrad, Christophe Cornut, Roman Fedorov, Philippe Gille, Mathieu Florence, João Lourenço, Siddharth Mathur, Zev Rosengarten, and Anis Zidani for helpful conversations and correspondence. We thank the Institute for Advanced Study for ideal conditions while working on parts of this project. This project has received funding from the European Research Council (ERC) under the European Union's Horizon 2020 research and innovation programme (grant agreement No.~851146). This project is based upon work supported by the National Science Foundation under Grant No.~DMS-1926686. This project was supported by the Simons Collaboration on Perfection in Algebra, Geometry, and Topology, award ID~MP-SCMPS-00001529-10.

 \section{The structure of $k$-group schemes locally of finite type}
 
 Our main result, \Cref{thm:main-GS}, deals with arbitrary group schemes $G$ locally of finite type over a field, so we begin by reviewing the structure theory of such $G$ in this chapter. This both prepares us for subsequent work by reviewing critical notions specific to imperfect fields (wound unipotent, quasi-reductive, pseudo-(abelian variety), pseudo-parabolic, etc.) and also shows how these notions arise naturally from an arbitrary $G$. More precisely, in \S\ref{sec:group-structure}, we review a fundamental filtration of $G$ by $k$-subgroups whose study is a fruitful way to approach an arbitrary $G$, and in subsequent \S\S\ref{sec:pass-to-smooth}--\ref{sec:pseudo-parabolic} we present techniques for attacking critical subquotients in this filtration. Most of this material is essentially a review, even if of facts that deserve to be known more widely, although the affineness of the ``cokernel'' of $i_{G,\, \overline{G}}$ established in \Cref{prop:Brauer-affine} is a new result that will be critically important in subsequent chapters.
 
 \csub[The fundamental filtration by $k$-subgroups and its subquotients] \label{sec:group-structure}
 
 We start by reviewing the structure theory of arbitrary $k$-group schemes locally of finite type.

\bpp[$k$-groups locally of finite type] \label{pp:k-groups-lft}
Throughout \S\ref{sec:group-structure}, we fix a $k$-group scheme $G$ locally of finite type. We recall from \cite{BLR90}*{Section 7.1, Lemma 2} that such a $G$ is automatically separated, and, from \cite{SGA3Inew}*{expos\'e VI${}_{\x{\upshape{A}}}$, proposition 2.5.2 (c)}, that a quasi-compact monomorphism of $k$-groups locally of finite type is necessarily a closed immersion. Certainly, this fails for non-quasi-compact monomorphisms, such as for $\underline{\mathbb{Z}}_{k} \xrightarrow{1 \mapsto 1} \mathbb{G}_{a,\, k}$ with $k$ of characteristic $0$, for which, relatedly, the $k$-group algebraic space $\mathbb{G}_{a,\, k}/\underline{\mathbb Z}_{k}$ locally of finite type is not quasi-separated (see also \S\ref{pp:conv}). 

On the other hand, we recall from \cite{Art69a}*{Lemma 4.2} that every quasi-separated $k$-group algebraic space is representable by a scheme. This means that we would not gain much by allowing $G$ to be an algebraic space and, more significantly, that $G/H$ is representable by a $k$-group scheme locally of finite type for any closed normal $k$-subgroup $H \lhd G$ (see also \S\ref{pp:conv}). Moreover, by the Chevalley theorem \cite{SGA3Inew}*{expos\'e VI${}_{\x{\upshape{B}}}$, th\'eor\`eme 11.17}, these quotient groups $G/H$ are all affine as soon as so is $G$. As far as their representability goes, however, much more is true: by \cite{Ray70b}*{chapitre~VI, corollaire~2.6} (with \cite{SGA3Inew}*{expos\'e~VI${}_{\x{\upshape{A}}}$, th\'eor\`eme 3.2 (ii), (iv) (a$'$); expos\'e V, th\'eor\`eme~4.1~(iv)}), every quasi-separated homogeneous space under $G$ over $k$ is representable by a $k$-scheme each of whose quasi-compact opens is quasi-projective. Knowing this quasi-projectivity is useful when dealing with the representability question of restrictions of scalars, see \S\ref{pp:conv}.
\epp

\bpp[The fundamental filtration]\label{pp:fundamental-filtration}
To study our general $G$, it is useful to keep in mind its following (closed) $k$-subgroups.

\[
\xymatrix@C=38pt{
{}\save[]+<3.3cm,-0.5cm>*\txt<9pc>{%
\fontsize{8pt}{8pt}\selectfont quotient is quasi-reductive \\ $\overbrace{ \quad \quad \quad \quad\ \ \quad \quad \quad \quad \quad}_{}$ } \restore& {}\save[]+<2cm,-6.25cm>*\txt<6pc>{%
\fontsize{8pt}{8pt}\selectfont quotient is pseudo-reductive}  \restore    &  \ar@{}[d]|-{\text{\scalebox{1.7}{\rotatebox{-90}{$\le$}}}} P_\lambda {}\save[]+<-1.4cm,-3.9cm>*\txt<6pc>{%
\rotatebox{90}{$\xrightarrow{\hspace{3.65cm}}$} } \restore   &   {}\save[]+<0.5cm,-0.25cm>*\txt<6pc>{%
\fontsize{8pt}{8pt}\selectfont quotient is an abelian variety}  \restore  & G^{\mathrm{gred}}  \ar@{}[r]|-{\text{\scalebox{1.7}{$\le$}}} & G \\
\mathscr{R}_{\mathrm{us},\, k}(G) \ar@{}[r]|-{\text{\scalebox{1.7}{$\lhd$}}} &  \quad \mathscr{R}_{\mathrm{u},\, k}(G)  \ar@{}[r]|-{\text{\scalebox{1.7}{$\lhd$}}}  & G^{\mathrm{sm,\, lin}}   \ar@{}[r]|-{\text{\scalebox{1.7}{$\lhd$}}}  &   G^{\mathrm{lin}} \ar@{}[r]|-{\text{\scalebox{1.7}{$\lhd$}}} & (G^{\mathrm{gred}})^0 \ar@{}[u]|-{\text{\scalebox{1.7}{\rotatebox{90}{$\lhd$}}}} \ar@{}[r]|-{\text{\scalebox{1.7}{$\le$}}} & G^0\ar@{}[u]|-{\text{\scalebox{1.7}{\rotatebox{90}{$\lhd$}}}}
 \\
\vdots \ar@{}[u]|-{\text{\scalebox{1.7}{\rotatebox{90}{$\lhd$}}}} & \vdots \ar@{}[u]|-{\text{\scalebox{1.7}{\rotatebox{90}{$\lhd$}}}} & \ar@{}[u]|-{\text{\scalebox{1.7}{\rotatebox{90}{$\lhd$}}}} G^{\mathrm{uni}} {}\save[]+<2.4cm,0.5cm>*\txt<7.5pc>{%
$\underbrace{ \quad \quad \quad \quad\quad \quad \quad \quad}_{}$  \vspace{-10pt} \\
\fontsize{8pt}{8pt}\selectfont  quotient is a pseudo-(abelian variety) } \restore    &    &  {}\save[]+<-10.1cm,0.5cm>*\txt<5.1pc>{%
$\underbrace{ \quad \quad \quad \quad}_{\substack{\text{quotient is} \\\text{wound unipotent}}}$ } \restore  & {}\save[]+<-1cm,0.5cm>*\txt<7.5pc>{%
$\underbrace{ \quad \quad \quad \quad}_{}$ \vspace{-10pt} \\
\fontsize{8pt}{8pt}\selectfont  has a filtration with quotients $\subset \mathrm{Res}_{k'/k}(\text{finite})$ } \restore &\\
(\mathscr{R}_{\mathrm{us},\, k}(G))_{i} \ar@{}[u]|-{\text{\scalebox{1.7}{\rotatebox{90}{$\lhd$}}}} & (\mathscr{R}_{\mathrm{u},\, k}(G))_{j} \ar@{}[u]|-{\text{\scalebox{1.7}{\rotatebox{90}{$\lhd$}}}} & \ar@{}[u]|-{\text{\scalebox{1.7}{\rotatebox{90}{$\lhd$}}}} G^{\mathrm{tor}} {}\save[]+<3.3cm,3.75cm>*\txt<6pc>{%
\rotatebox{-50}{$\xrightarrow{\hspace{0.6cm}}$} } \restore   &  &&   & \\
\vdots \ar@{}[u]|-{\text{\scalebox{1.7}{\rotatebox{90}{$\lhd$}}}} & \vdots \ar@{}[u]|-{\text{\scalebox{1.7}{\rotatebox{90}{$\lhd$}}}} &&& && \\
 (\mathscr{R}_{\mathrm{us},\, k}(G))_{1}  \ar@{}[u]|-{\text{\scalebox{1.7}{\rotatebox{90}{$\lhd$}}}} & (\mathscr{R}_{\mathrm{u},\, k}(G))_{1} \ar@{}[u]|-{\text{\scalebox{1.7}{\rotatebox{90}{$\lhd$}}}}  &&&&&\\
(\mathscr{R}_{\mathrm{us},\, k}(G))_{0}\ \ \ar@{}[r]|-{\text{\scalebox{1.35}{$=1=$}}} \ar@{}[u]|-{\text{\scalebox{1.7}{\rotatebox{90}{$\lhd$}}}} & \ (\mathscr{R}_{\mathrm{u},\, k}(G))_{0}  \ar@{}[u]|-{\text{\scalebox{1.7}{\rotatebox{90}{$\lhd$}}}}  &&&&&
}
\]
We now define the notation appearing in this diagram, discuss the respective subgroups and their associated subquotients, and give forward references to more detailed further discussions.
\benuma
\m
\label{m:G0} {\bf The identity component $G^0$ of $G$.} This is the connected component of $G$ through which the identity section factors. By \cite{SGA3Inew}*{expos\'e VI${}_{\x{\upshape{A}}}$, th\'eor\`eme 2.6.5}, this $G^0$ is a clopen, geometrically connected, quasi-compact (so of finite type), normal $k$-subgroup of $G$. The formation of $G^0$ commutes with base change to any field extension and with passage to quotients by connected, normal $k$-subgroups of $G$. In particular (or by \cite{SGA3Inew}*{expos\'e~VI${}_{\x{\upshape{A}}}$, proposition 5.5.1}), the quotient $G/G^0$ is \'{e}tale, its base change to $k^s$ becomes a $k^s$-group scheme associated to some abstract group.

\m
\label{m:Ggred}  {\bf The smooth part $G^{\mathrm{gred}}$ of $G$.} This is the largest smooth, closed $k$-subgroup of $G$, equivalently, the largest geometrically reduced, closed $k$-subscheme of $G$ (hence the notation $(-)^{\mathrm{gred}}$; see \cite{SGA3Inew}*{expos\'e VI${}_{\x{\upshape{A}}}$, proposition 1.3.1 (2)} for the equivalence), we review its construction in \S\ref{pp:Xgred} below. The formation of $G^{\mathrm{gred}}$ commutes with base change to any separable field extension and with passage to quotients by smooth, normal, closed $k$-subgroups of $G$. In \S\ref{pp:gred-filtration} below, we recall that the inclusion $(G^{\mathrm{gred}})^0 \le G^0$ may be refined by further closed subgroups whose successive subquotients are subschemes of restrictions of scalars of the form $\mathrm{Res}_{k^\prime/k}(\text{finite }k^\prime\text{-scheme})$ for some finite, purely inseparable field extensions $k^\prime/k$; this tends to be useful for reducing to smooth groups. One reason why smooth $k$-groups are preferable is that, by Grothendieck's theorem \cite{CGP15}*{Lemma C.4.4}, they always have a maximal torus defined over $k$ (the same also holds for commutative groups but not in general, see \cite{SGA3Inew}*{expos\'e XVII, remarque 5.9.1}).

\m
\label{m:Glin}\textbf{The \emph{linear part} $G^{\mathrm{lin}}$ of $G$.} This is the smallest connected, affine, normal $k$-subgroup
\[\label{eqn:Chevalley}
\qq G^{\mathrm{lin}} \lhd (G^{\mathrm{gred}})^0 \quad\text{such that} \quad G^{\mathrm{av}}:= (G^{\mathrm{gred}})^0/G^{\mathrm{lin}} \quad \text{is an abelian variety.}
\]
 The Chevalley theorem \cite{BLR90}*{Section 9.2, Theorem 1} ensures that this $G^{\mathrm{lin}}$ exists, and, in the case when $k$ is perfect, that $G^{\mathrm{lin}}$ is smooth and that its formation commutes with base change to any field extension $k^\prime/k$. In contrast, when $k$ is imperfect, $G^{\mathrm{lin}}$ need not be smooth (see, for instance, \S\ref{pp:pseudo-abelian} below) and, by a limit, spreading out, and Galois descent argument, the base change property only holds for separable field extensions. Nevertheless, regardless of what $k$ is, $G^{\mathrm{lin}}$ has no nontrivial infinitesimal $k$-group quotients: indeed, if the intersection of the kernels of all such quotients was smaller than $G^{\mathrm{lin}}$, then that would contradict the definition of $G^{\mathrm{lin}}$. 

It might be more appropriate to call $G^{\mathrm{lin}}$ the \emph{connected linear part} of $G$, but we prefer brevity.

\m\label{m:Gsmlin}
\textbf{The \emph{smooth linear part} $G^{\mathrm{sm},\, \mathrm{lin}}$ of $G$.} This is the largest connected, smooth, affine $k$-subgroup of $G$, more succinctly, it is simply    
$$
\qq G^{\mathrm{sm},\,\mathrm{lin}} := ((G^{\mathrm{lin}})^{\mathrm{gred}})^0.
$$
It might be more appropriate to call $G^{\mathrm{sm},\, \mathrm{lin}}$ the \emph{connected smooth linear part} of $G$, but we again prefer brevity because this seems unlikely to cause confusion. Indeed, a general $G$ has no largest smooth, affine $k$-subgroup, as the example of the constant $k$-group $\underline{\mathbb{Q}/\mathbb{Z}}$ shows. 

The formation of $G^{\mathrm{sm},\,\mathrm{lin}}$ commutes with base change to any separable field extension and with passage to quotients by connected, smooth, affine normal $k$-subgroups of $G$. Its stability under conjugation by $k^s$-points of $G$ ensures that $G^{\mathrm{sm},\,\mathrm{lin}}$ is normal in $G^{\mathrm{gred}}$. 

Thanks to its smoothness, $G^{\mathrm{sm},\, \mathrm{lin}}$ is more manageable than $G^{\mathrm{lin}}$, although the quotient
$$
\qq G^{\mathrm{pav}}:=(G^{\mathrm{gred}})^0/G^{\mathrm{sm},\,\mathrm{lin}}
$$
is no longer an abelian variety, but only a \emph{pseudo-abelian} variety in the sense of Totaro, that is, it is smooth, connected, and has no nontrivial smooth, connected, affine $k$-subgroups. We review some aspects of pseudo-abelian varieties in \S\ref{sec:pseudo-abelian} below, and in \Cref{prop:PP-groups} we show that they are pseudo-proper (in the sense of \Cref{def:pseudo-proper}). This pseudo-properness is useful for controlling $G^{\mathrm{pav}}$, so also for reducing to connected, smooth, affine $k$-groups in practice. The quotient $G^{\mathrm{lin}}/G^{\mathrm{sm},\, \mathrm{lin}}$ is identified with $(G^{\mathrm{pav}})^{\mathrm{lin}}$, so it is pseudo-finite in the sense that its largest smooth, closed $k$-subgroup $(G^{\mathrm{lin}}/G^{\mathrm{sm},\, \mathrm{lin}})^{\mathrm{gred}}$ is \'etale (see \Cref{def:pseudo-finite}).

\m\label{m:RukG}
\textbf{The \emph{unipotent} $k$-\emph{radical} $\mathscr{R}_{\mathrm{u},\, k}(G)$ of $G$.} This is the largest connected, smooth, unipotent (see below), normal $k$-subgroup of the smooth linear part $G^{\mathrm{sm},\, \mathrm{lin}}$ (equivalently, of $G^{\mathrm{gred}}$), its existence is immediate from the definition, alternatively, one may refer to \cite{SGA3Inew}*{expos\'e~VI${}_{\x{\upshape{B}}}$, corollaire 7.1.1}. The formation of $\mathscr{R}_{\mathrm{u},\, k}(G)$ commutes with base change to any separable field extension and with passage to quotients by connected, smooth, unipotent, normal $k$-subgroups of $G$.  The quotient\footnote{We use the notation $G^{\mathrm{pred}}$ even when $k$ is perfect, for instance, algebraically closed, even though then $G^{\mathrm{pred}}$ is necessarily reductive. We reserve $(-)^{\mathrm{red}}$ for denoting the underlying reduced closed subscheme.}  
    \[\label{eqn:Gpred}
\qq    G^{\mathrm{pred}}:=G^{\mathrm{sm},\,\mathrm{lin}}/\mathscr{R}_{\mathrm{u},\, k}(G)
    \]
    is a \emph{pseudo-reductive} $k$-group in the sense that it is connected, smooth, affine, and has a trivial unipotent $k$-radical, that is, $\mathscr{R}_{\mathrm{u},\, k}(G^{\mathrm{pred}}) = 1$. Pseudo-reductive groups form the most delicate part of the entire diagram above, and analyzing them is subtle. In \S\S\ref{sec:iG}--\ref{sec:pseudo-reductive} below, we present a widely useful and somewhat underappreciated framework for handling them, more precisely, for reducing to reductive groups.

\m \label{m:RuskG}\textbf{The \emph{split unipotent} $k$-\emph{radical} $\mathscr{R}_{\mathrm{us},\, k}(G)$ of $G$.} This is the largest split unipotent (see below), normal $k$-subgroup of the unipotent $k$-radical $\mathscr{R}_{\mathrm{u},\, k}(G)$ (equivalently, of $G^{\mathrm{gred}}$), it exists by, for instance, \cite{CGP15}*{Theorem B.3.4}. Here we recall that a unipotent $k$-group is \emph{split} if it is an iterated extension of the additive group $\mathbb{G}_{a,\, k}$. It is then also connected and smooth, equivalently, a connected, smooth, unipotent $k$-group is split if and only if it admits a dominant $k$-morphism from some $\mathbb{A}^n_k$, in which case it is even isomorphic to $\mathbb{A}^n_k$ as a $k$-scheme, so that every $k$-group quotient of a split unipotent $k$-group is split unipotent, see \cite{Con15b}*{Corollary 3.9}. In contrast, a unipotent $k$-group is \emph{wound} if it has no $\mathbb{G}_{a,\, k}$ as a $k$-subgroup (see \cite{Ros25}*{Proposition A.1} for equivalent characterizations of woundness).\footnote{Some authors require wound unipotent groups to be smooth by definition; our terminology agrees with \cite{BLR90}*{top of p. 174}, except that we do not require wound groups to be connected.} By \cite{CGP15}*{Theorem B.3.4}, the formation of $\mathscr{R}_{\mathrm{us},\, k}(G)$ commutes with base change to separable field extensions and with passage to quotients by split unipotent, normal $k$-subgroups of $G$. 
    
    The quotient $\mathscr{R}_{\mathrm{u},\, k}(G)/\mathscr{R}_{\mathrm{us},\, k}(G)$ is wound unipotent, so
    $$
 \qq   G^{\mathrm{qred}}:=G^{\mathrm{sm},\,\mathrm{lin}}/\mathscr{R}_{\mathrm{us},\, k}(G)
    $$
    is a quasi-reductive $k$-group in the sense that it is connected, smooth, affine, and has a trivial split unipotent $k$-radical, so that $\mathscr{R}_{\mathrm{us},\, k}(G^{\mathrm{qred}}) = 1$, equivalently, so that its unipotent $k$-radical $\mathscr{R}_{\mathrm{u},\, k}(G)$ is wound.

\m
\label{m:cckp} {\bf The \emph{iterated cc$k$p kernels} $(\mathscr{R}_{\mathrm{u},\, k}(G))_i$ and $(\mathscr{R}_{\mathrm{us},\, k}(G))_j$.} These are defined inductively for any smooth, unipotent $k$-group $U$ as follows: $U$ has a unique largest connected, smooth, central, $p$-torsion $k$-subgroup $U_1 \lhd U$, the cc$k$p kernel of $U$ (see \cite{CGP15}*{Definition B.3.1}), and for $i > 1$ one inductively sets $U_i := (U/U_{i - 1})_1$. By \emph{loc.~cit.}~and \cite{CGP15}*{Corollary B.3.3}, the filtration $\{U_i\}_{i \ge 0}$ is exhaustive, its formation commutes with base change to separable field extensions, and if $U$ is wound (resp.,~strongly wound, see \Cref{def:strongly-wound} below), then so is every subquotient $U_{i^\prime}/U_i$ for $i' \ge i$ (see \Cref{lem:sw-cckp} below). This last aspect is remarkable because woundness is most often not inherited by quotients, see \S\ref{pp:wound-structure} below or \cite{Ros25}*{before Definition 1.2, also Proposition 7.7}.

\m
\label{m:pseudo-parabolic} {\bf The pseudo-parabolic subgroups $P_\lambda \le G^{\mathrm{sm},\, \mathrm{lin}}$.} These are certain connected, smooth, affine $k$-subgroups that contain $\mathscr{R}_{u,\, k}(G)$ and are associated to $k$-group homomorphisms $\lambda\colon \mathbb{G}_{m,\, k} \rightarrow G$, we review them in \S\ref{pp:pseudo-parabolic} below. By \cite{CGP15}*{Proposition 3.5.2 (1)}, the pseudo-parabolicity of a smooth $k$-subgroup $P \le G$ may be tested after base change to any separable field extension $k^\prime / k$, and it is also insensitive to base change to any such extension. Each pseudo-parabolic $P_\lambda$ is the preimage of the corresponding pseudo-parabolic of $G^{\mathrm{pred}}$. If the induced $\lambda \colon \mathbb{G}_{m,\, k} \rightarrow G^{\mathrm{pred}}$ is noncentral, then $P_\lambda$ is strictly smaller than $G^{\mathrm{sm},\, \mathrm{lin}}$, so the pseudo-reductive $k$-group $P_\lambda^{\mathrm{pred}}$ is ``smaller'' than $G^{\mathrm{pred}}$, and is even commutative when $(P_\lambda)_{k^s}$ is minimal among the pseudo-parabolics of $G_{k^s}$. Thus, pseudo-parabolics aid the study of the most delicate part $G^{\mathrm{pred}}$ of the diagram above by facilitating passage to ``smaller'' pseudo-reductive groups, all the way up to commutative pseudo-reductive groups, for which nonabelian complexities disappear. To aid this further, we prove in \Cref{thm:pseudo-flag} below that the homogeneous spaces $G^{\mathrm{sm},\, \mathrm{lin}}/P_\lambda$, in fact, already the $(G^{\mathrm{gred}})^0/P_\lambda$, are all pseudo-proper.

\m
\label{m:Guni}{\bf The largest unirational $k$-subgroup $G^{\mathrm{uni}} \le G$.} This is the largest unirational (closed) $k$-subgroup of $G$, it exists by \cite{BLR90}*{bottom of p. 310}. Here we recall that a finite type, reduced $k$-scheme $X$ is \emph{unirational} if there is a dominant rational $k$-morphism $\mathbb{A}^n_k\dashrightarrow X$ (see \cite{Bor91}*{Section AG.13.7}), equivalently, if $X$ is integral with a function field that is a subfield of some purely transcendental extension $k(t_1, \dotsc, t_n)$ over $k$. For every nonempty open $U$ of a unirational $X$ and every semilocal $k$-algebra $A$ with infinite residue fields, we have $U(A) \neq \emptyset$, in particular, if $k$ is infinite, then $X(k)$ is Zariski dense in $X$. Every unirational $X$ is generically smooth, in particular, our $G^{\mathrm{uni}}$ is smooth and connected, and every quotient of $G^{\mathrm{uni}}$ is again unirational. All maps from $\mathbb{A}^1_k$ to abelian varieties are constant, so $G^{\mathrm{uni}} \le G^{\mathrm{sm},\, \mathrm{lin}}$, in other words, $G^{\mathrm{uni}}$ is affine.

By a recent result of Rosengarten \cite{Ros24}*{Theorem 1.6}, the formation of $G^{\mathrm{uni}}$ commutes with base change to separable field extensions, in particular, $G^{\mathrm{uni}}$ is normal in $G^{\mathrm{gred}}$. By considering function fields, we see that the formation of $G^{\mathrm{uni}}$ also commutes with passage to quotients by split unipotent, normal $k$-subgroups of $G$. 
For a smooth $k$-group $G$ we may keep iteratively forming quotients by $(-)^{\mathrm{uni}}$ to eventually reduce to the case when $G^{\mathrm{uni}} = 1$. Bosch--L\"utkebohmert--Raynaud have characterized groups at which this process stops in the commutative case, more precisely, by \cite{BLR90}*{Section 10.3, Theorem 1}, for a connected, smooth, \emph{commutative} $k$-group $G$, the following are equivalent:
\benumr
    \item $G^{\mathrm{uni}} = 0$;
    \item $G(S) \xrightarrow{\sim} G(U)$ for every dense open immersion $U \subset S$ of smooth $k$-schemes;
\end{enumerate}
and, granted that $G$ is a dense open in a proper, regular $k$-scheme $\overline{G}$, these are equivalent to
\begin{enumerate}
    \item[(iii)] $\overline{G}^{\mathrm{sm}} = G$, that is, $G$ is precisely the $k$-smooth locus of $\overline{G}$.
\end{enumerate}
By the resolution of singularities conjecture, such a $\overline{G}$ ought to exist, and in practice one may sometimes build it, for instance, for groups given by the vanishing of $p$-polynomials as in \S\ref{pp:wound-structure} below, compare with \cite{Tot13}*{Example 9.7}. For a version of the criterion above beyond commutative groups, see \Cref{prop:PP-groups}~\ref{m:PPG-i} below.

\m
\label{m:Gtor} {\bf The $k$-subgroup $G^{\mathrm{tor}} \le G$ generated by the $k$-tori.} This is the $k$-subgroup of $G$, equivalently, of $G^{\mathrm{sm},\, \mathrm{lin}}$, generated by the $k$-tori of $G$. It is connected, smooth, affine, and unirational over $k$ because $k$-tori are unirational (see \cite{Bor91}*{Chapter III, Example 8.13 (2)}), so that $G^{\mathrm{tor}}$ lies in $G^{\mathrm{uni}}$. By \cite{CGP15}*{Proposition A.2.11} (and \ref{m:Gsmlin} above), the formation of $G^{\mathrm{tor}}$ commutes with base change to any separable field extension (and even to any field extension if $G$ is smooth), $G^{\mathrm{tor}}$ is normal in $G^{\mathrm{gred}}$, and $G^{\mathrm{sm},\, \mathrm{lin}}/G^{\mathrm{tor}}$ is unipotent. If $G$ is quasi-reductive (resp.,~pseudo-reductive), then, by normality, so are $G^{\mathrm{tor}}$ and $G^{\mathrm{uni}}$. Conversely, $G^{\mathrm{tor}} = 1$ if and only if $G^{\mathrm{sm},\, \mathrm{lin}}$ is unipotent, that is, if and only if $G^{\mathrm{pred}} = 1$; whereas $G^{\mathrm{uni}} = 1$ if and only if $G^{\mathrm{sm},\, \mathrm{lin}}$ is strongly wound unipotent (see \Cref{def:strongly-wound} below).
\eenum
\epp

The formation of the subgroups \ref{m:G0}--\ref{m:Gtor}, so of the entire diagram above, is functorial in $G$ and commutes with products and with base change to separable field extensions, so also with restrictions of scalars along finite separable field extensions. In particular, being a pseudo-abelian variety (resp., a pseudo-reductive group; resp., a quasi-reductive group; resp., a wound unipotent group) commutes with and may be tested after such a base change, and if a smooth $k$-group acts by group automorphisms on $G$, or even merely on $G^{\mathrm{gred}}$, then this action preserves the $k$-subgroups $G^{\mathrm{gred}}$, $(G^{\mathrm{gred}})^0$, $G^{\mathrm{lin}}$, $G^{\mathrm{sm},\,\mathrm{lin}}$, $G^{\mathrm{uni}}$, $G^{\mathrm{tor}}$, $\mathscr{R}_{\mathrm{u},\, k}(G)$, $\mathscr{R}_{\mathrm{us},\, k}(G)$, $(\mathscr{R}_{\mathrm{u},\, k}(G))_i$, and $(\mathscr{R}_{\mathrm{us},\, k}(G))_j$. By letting this action be that of $G^{\mathrm{gred}}$ on itself by conjugation, we see that these subgroups are normal even in~$G^{\mathrm{gred}}$.

Granted that we discard $G^0$, $G^{\mathrm{lin}}$, $G^{\mathrm{uni}}$, and $G^{\mathrm{tor}}$, the formation of the diagram displayed above also commutes with restrictions of scalars along arbitrary finite field extensions: for this, the settled separable case reduces us to purely inseparable extensions, then we note the preservation of $G^{\mathrm{gred}}$ (see \S\ref{pp:conv} and \S\ref{pp:Xgred}), then also of $(G^{\mathrm{gred}})^0$ (also use \cite{CGP15}*{Proposition A.5.9}), then also of $G^{\mathrm{sm},\, \mathrm{lin}}$, $\mathscr{R}_{\mathrm{u},\, k}(G)$, $\mathscr{R}_{\mathrm{us},\, k}(G)$ (see also \cite{CGP15}*{proof of Proposition 1.1.10}), then also of the pseudo-parabolic subgroups (see \cite{CGP15}*{Proposition 2.2.13}), and finally of the cc$k$p filtrations (combine \cite{CGP15}*{Proposition A.5.15 (1)} with the previous steps). In particular, pseudo-abelian varieties (resp.,~pseudo-reductive groups; resp.,~quasi-reductive groups; resp.,~smooth, wound unipotent groups) are stable under restrictions of scalars, and they are also stable under extensions and under passage to connected, smooth, normal $k$-subgroups.

\bpp[The case when $k$ is perfect] \label{pp:perfect-case}
We will be especially interested in the case when the field $k$ is imperfect: this is when the geometric phenomena are particularly rich and when all the inclusions of subgroups in the diagram displayed above are in general strict. In contrast, in the case when $k$ is perfect, the structure theory above simplifies as follows.
\epp
\benumr
    \item \label{m:PC-i}
We have $G^{\mathrm{lin}} = G^{\mathrm{sm},\, \mathrm{lin}}$ and the pseudo-abelian variety $G^{\mathrm{pav}}$ is an abelian variety, see \ref{m:Glin}. In fact, the following inclusion is an equality precisely over perfect fields:
    $$
  \qq  \{\text{abelian varieties}\} \subseteq \{\text{pseudo-abelian varieties}\}
    $$
    (to see the strictness of the inclusion over every imperfect field, consider restrictions of scalars of abelian varieties from purely inseparable extensions and see \cite{CGP15}*{Example A.5.6}; for more interesting examples, see \cite{Tot13}*{Corollaries 6.5 and 7.3}).
    
    \item \label{m:PC-ii} We have $\mathscr{R}_{\mathrm{us},\, k}(G) = \mathscr{R}_{\mathrm{u},\, k}(G)$, both of them descend $\mathscr{R}_{u,\, \overline{k}}(G_{\overline{k}})$ to $k$, and both $G^{\mathrm{pred}}$ and $G^{\mathrm{qred}}$ are reductive, see \ref{m:RukG} and \ref{m:RuskG}. In fact, each of the following inclusions is an equality precisely over perfect fields:
    $$
\qq    \{\text{reductive groups}\} \subseteq \{\text{pseudo-reductive groups}\} \subseteq \{\text{quasi-reductive groups}\}
    $$
    (for strictness of the inclusions over every imperfect field, again consider restrictions of scalars). Moreover, over perfect fields every connected, smooth unipotent group is split, that is, the woundness phenomenon is specific to imperfect fields.
\eenum
    In the case when $k$ is even of characteristic $0$, by the Cartier theorem \cite{SGA3Inew}*{expos\'e VI${}_{\x{\upshape{B}}}$, corollaire 1.6.1}, every locally of finite type $k$-group is smooth, that is,
    \benum
    \item[(iii)] \label{m:PC-iii} We have $G^{\mathrm{gred}} = G$.
    \eenum
    Over a general field $k$, the Cartier theorem has a useful generalization: by \cite{SGA3Inew}*{expos\'e~VII${}_{\x{\upshape{A}}}$, proposition~8.3}, for any $k$-group scheme $G$ locally of finite type, there is an infinitesimal, normal $k$-subgroup $\alpha_G \lhd G$ that may be chosen to be any sufficiently large Frobenius kernel of $G$ such that the $k$-group $G/\alpha_G$ is smooth.

\bpp[The derived subgroup $\mathscr D(G)$] \label{pp:DG}
In addition to \S\ref{pp:fundamental-filtration}, every smooth $k$-group scheme $G$ has a derived subgroup $\mathscr{D}(G)$ defined as the smallest closed $k$-subgroup containing the image of the commutator map $G \times_k G \xrightarrow{(g,\,h)\, \mapsto\, ghg^{-1}h^{-1}} G$, see \cite{SGA3Inew}*{expos\'e VI${}_{\x{\upshape{B}}}$, d\'efinition 7.2.2 b)}. By \cite{SGA3Inew}*{expos\'e VI${}_{\x{\upshape{B}}}$, proposition 7.1, corollaire 7.3}, this $\mathscr{D}(G)$ exists, is smooth and normal in $G$, its formation commutes with base change to arbitrary field extensions, $G/\mathscr{D}(G)$ is commutative, and $G$ is commutative if and only if $\mathscr{D}(G) = 1$. By \cite{SGA3Inew}*{expos\'e VI${}_{\x{\upshape{B}}}$, corollaire 7.2.1, proposition~7.8, corollaire 7.10}, if $G$ is of finite type, then $\mathscr{D}(G)$ is the image (as fpqc sheaves) of the commutator map above and $\mathscr{D}(G)(\overline{k})$ is the derived subgroup of $G(\overline{k})$, and if $G$ is connected, then so is $\mathscr{D}(G)$.

If a smooth $k$-group $G$ is affine (resp.,~pseudo-reductive; resp.,~quasi-reductive; resp.,~unipotent; resp.,~wound unipotent; resp.,~unirational), then so is $\mathscr{D}(G)$. If $G = G^{\mathrm{tor}}$, then \cite{CGP15}*{Propositions~A.2.8, A.2.10} ensure that $G/\mathscr{D}(G)$ is a torus (it is connected, smooth, commutative, and generated by its $k$-tori) that is a quotient of any maximal $k$-torus of $G$. 

A smooth $k$-group scheme $G$ is \textit{perfect} if $\mathscr{D}(G) = G$. A \textit{pseudo-semisimple} (resp.,~\textit{quasi-semisimple}) $k$-group is a pseudo-reductive (resp.,~quasi-reductive) $k$-group that is perfect, so a reductive $k$-group is quasi-semisimple (resp.,~pseudo-semisimple) if and only if it is semisimple.
\epp


\csub[Passage to smooth groups and pseudo-finiteness]

\label{sec:pass-to-smooth}

Other than the identity component $G^0$ discussed in \S\ref{pp:fundamental-filtration}~\ref{m:G0}, the first piece of the fundamental filtration of $G$ is the largest smooth, closed $k$-subgroup $G^{\mathrm{gred}}$. We review the construction of $G^{\mathrm{gred}}$ and recall a useful technique for passing from $G$ to $G^{\mathrm{gred}}$.

\bpp[The underlying geometrically reduced subspace $X^{\mathrm{gred}} \subset X$] \label{pp:Xgred}
Each locally of finite type $k$-algebraic space $X$ has the largest geometrically reduced, closed $k$-subspace
\[
X^{\mathrm{gred}} \subseteq X^{\mathrm{red}} \subseteq X
\]
defined as follows. One defines $X^{\mathrm{gred}}$ by Galois descent by declaring $(X^{\mathrm{gred}})_{k^s}$ to be the schematic (or merely Zariski) closure of $X(k^s)$ in $X_{k^s}$, see \cite{SP}*{Lemmas~\href{https://stacks.math.columbia.edu/tag/082X}{082X} and \href{https://stacks.math.columbia.edu/tag/0830}{0830}}. By \emph{loc.~cit.},~the definition is compatible with pullback under any \'etale $k$-morphism $X^\prime \rightarrow X$, so, by considering an atlas, the claim that this $X^{\mathrm{gred}}$ is the largest geometrically reduced, closed subspace reduces to the scheme case that was settled in \cite{CGP15}*{Lemma~C.4.1}. Of course, if $k$ is perfect, then we simply have $X^{\mathrm{gred}} = X^{\mathrm{red}}$, but this fails for imperfect $k$: for instance, if $X = \mathrm{Spec}(k^\prime)$ for a nontrivial, finite, purely inseparable extension $k^\prime/k$, then $X^{\mathrm{gred}} = \emptyset$.

By definition, $X^{\mathrm{gred}}$ is functorial in $X$ and commutes with pullback along any smooth $k$-morphism $X^\prime \rightarrow X$. Moreover, as in \emph{loc.~cit.},~$X^{\mathrm{gred}}(k^\prime) = X(k^\prime)$ for every separable field extension $k^\prime/k$ and the formation of $X^{\mathrm{gred}}$ commutes with products and with base change to any such $k^\prime$ (in particular, if $X(k^\prime) \neq \emptyset$, then already $X(k^s) \neq \emptyset$). Therefore, the inclusion 
\[\label{eqn:reduce-to-gred}
X^{\mathrm{gred}}(S) \subseteq X(S)
\]
 is an equality for every integral $k$-algebraic space $S$ whose function field is a separable extension of $k$ (the pullback of $X^{\mathrm{gred}}$ under any $S$-point of $X$ is a closed subspace of $S$ that contains the generic point), see \cite{SP}*{Definition~\href{https://stacks.math.columbia.edu/tag/0ENE}{0ENE}}.

The compatibility with products and the functoriality of $X^{\mathrm{gred}}$ imply that if $X$ is a $k$-group, then $X^{\mathrm{gred}}$ is its closed $k$-subgroup, the largest smooth (see \S\ref{pp:fundamental-filtration}~\ref{m:Ggred}), closed $k$-subgroup of $X$. Similarly, if $X$ is a $k$-group and $E$ is an $X$-torsor over $k$, then $X^{\mathrm{gred}}$ is either empty or a $X^{\mathrm{gred}}$-torsor whose induced $X$-torsor is $E$: in fact, $X^{\mathrm{gred}} \neq \emptyset$ exactly when $E$ trivializes over some separable extension of $k$, equivalently, over $k^s$.
\epp

\begin{definition} \label{def:pseudo-finite}
A $k$-algebraic space $X$ is
\begin{itemize}
    \item \emph{pseudo-\up{locally quasi-finite}} if it is locally of finite type and $X^{\mathrm{gred}}$ is locally quasi-finite (equivalently, \'etale; equivalently, $X$ has no positive-dimensional $k$-smooth subspaces; equivalently, each quasi-compact open subscheme of $X$ has only finitely many $k^s$-points);
    \item \emph{pseudo-finite} if it is separated, of finite type, and $X^{\mathrm{gred}}$ is finite (equivalently, finite \'etale; equivalently, $X$ has only finitely many $k^s$-points).
\end{itemize}
\end{definition}

Being pseudo-finite is what occurs in practice, although being pseudo-(locally quasi-finite) is useful for making some statements sharp, without unnecessary hypotheses. In \S\ref{pp:iG}, \S\ref{pp:pav-iG}, and \ref{pp:pred-iG} below, we will see many pseudo-finite, positive-dimensional $k$-groups relevant for the study of pseudo-abelian varieties and pseudo-reductive groups. 

\begin{remark} \label{rem:Xgred-etale}
By \S\ref{pp:Xgred}, being pseudo-(locally quasi-finite) (resp.,~pseudo-finite) is stable under and may be tested after base change to any separable field extension $k^\prime/k$. Moreover, if $X$ is pseudo-(locally quasi-finite), then $X^{\mathrm{gred}}$ is a separated (\'etale) $k$-scheme, see \cite{SP}*{Lemma~\href{https://stacks.math.columbia.edu/tag/06LZ}{06LZ}}.
\end{remark}

\begin{remark}
Pseudo-(locally quasi-finite) (resp.,~pseudo-finite) $k$-group schemes are stable under passing to $k$-subgroups and extensions, but are not stable under quotients.
\end{remark}

\begin{remark} \label{rem:pseudo-finite-perfect}
For perfect $k$, an $X$ is pseudo-(locally quasi-finite) (resp.,~pseudo-finite) if and only if it is locally quasi-finite (resp.,~finite). Over imperfect $k$, however, pseudo-finite yet positive-dimensional $k$-schemes are pervasive, as the following examples show.
\end{remark}

\beg\label{eg:pseudo-etale}
If $X$ is pseudo-finite, then so is every $k$-subspace of a quasi-finite, separated $X$-algebraic space $X^\prime$. For instance, every finite $k$-scheme is pseudo-finite and, more interestingly, for every finite extension $k^\prime/k$ and every finite (or merely pseudo-finite) $k^\prime$-scheme $X^\prime$, every subspace of $X := \mathrm{Res}_{k^\prime/k}X^\prime$ is pseudo-finite. When $k^\prime/k$ is not separable, such $X$ may have an arbitrarily large dimension, as happens already with $X^\prime = \alpha_p$.
\eeg

\beg \label{eg:pseudo-finite-G}
For a $k$-group scheme $G$ locally of finite type, the quotient $G/G^{\mathrm{gred}}$ is pseudo-(locally quasi-finite) (resp.,~even pseudo-finite if $G$ is of finite type): indeed, its only $k^s$-point is the identity section. Similarly, for a $G$-torsor $E$ over $k$, the quotient $E/G^{\mathrm{gred}}$ is pseudo-(locally quasi-finite) and $(E/G^{\mathrm{gred}})^{\mathrm{gred}} \neq \emptyset$ if and only if $E$ trivializes over $k^s$.
\eeg

To study pseudo-finite groups effectively, and also to pass from $G$ to $G^{\mathrm{gred}}$ for general $k$-group schemes of finite type, we use the following method of going deeper into $G$ until we reach $G^{\mathrm{gred}}$.

\bpp[Filtration by stabilizers (compare with \cite{GGMB14}*{{\bf section 5.4}})] \label{pp:gred-filtration}
As we will see, any $k$-group scheme $G$ of finite type has a filtration by $k$-subgroups
\[\label{eqn:gred-filtration}
G^{\mathrm{gred}} = G^{(n)} \lneq G^{(n-1)} \lneq \cdots \lneq G^{(0)} = G \quad (\text{so that } G^{\mathrm{gred}} = (G^{(i)})^{\mathrm{gred}} \text{ for all } i)
\]
such that each $G^{(i)}$ acts on the affine $k$-scheme $\mathrm{Res}_{k_i/k}(Q_i)$ for some finite, purely inseparable field extension $k_i/k$ and $k_i$-group quotient $Q_i \colonequals G^{(i)}_{k_i}/(G^{(i)}_{k_i})^{\mathrm{red}}$, and $G^{(i+1)}$ is the stabilizer of the unique $k$-point of $\mathrm{Res}_{k_i/k}(Q_i)$. Indeed, inductively on $i$,
\begin{itemize}
    \item we choose $k_i$ to be a field of definition of the $k^{\mathrm{perf}}$-subgroup $(G^{(i)}_{k^{\mathrm{perf}}})^{\mathrm{red}} \le G^{(i)}_{k^{\mathrm{perf}}}$, so that $(G^{(i)}_{k_i})^{\mathrm{gred}}$ is a $k_i$-subgroup of $G^{(i)}_{k_i}$ that descends $(G^{(i)}_{k^{\mathrm{perf}}})^{\mathrm{red}}$ and $(G^{(i)}_{k_i})^{\mathrm{gred}} \cong (G^{(i)}_{k_i})^{\mathrm{red}}$;
    \item we set $Q_i := G^{(i)}_{k_i}/(G^{(i)}_{k_i})^{\mathrm{gred}}$, so that $Q_i$ is $k_i$-finite, connected, and $Q_i(k_i)$ and $Q_i(k_i^s)$ are both singletons (see \Cref{eg:pseudo-finite-G});
    \item we note that the $G^{(i)}_{k_i}$-action on $Q_i$ gives rise to a $\mathrm{Res}_{k_i/k}(G^{(i)}_{k_i})$-action on $\mathrm{Res}_{k_i/k}(Q_i)$, so, by restricting to $G^{(i)} \le \mathrm{Res}_{k_i/k}(G^{(i)}_{k_i})$, also to a $G^{(i)}$-action on $\mathrm{Res}_{k_i/k}(Q_i)$;
    \item we let $G^{(i + 1)} \le G^{(i)}$ be the stabilizer of the unique $k$-point of $\mathrm{Res}_{k_i/k}(Q_i)$;
    \item we note that $G^{(i + 1)}$ does not depend on the choice of $k_i$ because enlarging the latter has the effect of $G^{(i)}$-equivariantly embedding $\mathrm{Res}_{k_i/k}(Q_i)$ into a larger such restriction of scalars.
\end{itemize}
If $Q_i$ is not reduced, then $G^{(i + 1)} \lneq G^{(i)}$: indeed, the counit $(\mathrm{Res}_{k_i/k}(-))_{k_i} \rightarrow (-)$ of the adjunction is functorial and commutes with products, so $(\mathrm{Res}_{k_i/k}(Q_i))_{k_i} \rightarrow Q_i$ is $G^{(i)}_{k_i}$-equivariant, to the effect that if we had $G^{(i + 1)} = G^{(i)}$, then $G^{(i)}_{k_i}$ could not act transitively on $Q_i$. Moreover, $G^{(i + 1)}$ contains $(G^{(i)})^{\mathrm{gred}} \cong G^{\mathrm{gred}}$ because $(\mathrm{Res}_{k_i/k}(Q_i))(k^s)$ is a singleton and $G^{(i)}(k^s)$ is schematically dense in $(G^{(i)})^{\mathrm{gred}}$. By Noetherian induction, since $G$ is of finite type and we have $G^{(i + 1)} \lneq G^{(i)}$ whenever $G^{\mathrm{gred}} \lneq G^{(i)}$, the filtration eventually stabilizes at $G^{\mathrm{gred}}$. 

By construction, the filtration \eqref{eqn:gred-filtration} is functorial in $G$ and commutes with products and with base change along separable field extensions, so also with restriction of scalars along finite separable field extensions. In particular, if a smooth $k$-group acts by $k$-group automorphisms on $G$, then this action preserves the $k$-subgroups $G^{(i)}$.
\epp

\begin{lemma} \label{cor:pfinite-affine}
In \uS\uref{pp:gred-filtration}, for every $k$-group $G$ of finite type, the $G^{(i)}/G^{(i')}$ for $i' \ge i$ are quasi-affine. In particular, a pseudo-finite $k$-group scheme is affine.
\end{lemma}

\begin{proof}
Each $G^{(i)}/G^{(i + 1)}$ is quasi-affine: by \cite{GGMB14}*{remarques~2.1.4~(i)}, it is a subscheme~of~the affine $k$-scheme $\mathrm{Res}_{k_i/k}(Q_i)$. Thus, by induction and descent for quasi-affine maps \cite{SP}*{Lemma~\href{https://stacks.math.columbia.edu/tag/02L7}{02L7}}, each $G^{(i)}/G^{(i')}$ is quasi-affine. If $G$ is pseudo-finite, then $G^{\mathrm{gred}}$ is finite \'etale, so it follows that $G$ is quasi-affine. Then \cite{SGA3Inew}*{expos\'e VI${}_{\x{\upshape{B}}}$, proposition 11.11} shows that $G$ is affine. 
\end{proof}

\begin{remark}
In \cite{Tot13}*{Lemma 6.3}, assuming that $k$ is imperfect, Totaro constructs pseudo-finite, commutative extensions of connected, smooth, commutative, $p$-torsion $k$-groups $U$ by $\alpha_p$. \Cref{cor:pfinite-affine} implies that such extensions do not exist if instead, for instance, $U$ is a nonzero abelian variety, although this also follows already from \cite{Gro62}*{proposition 3.1}.
\end{remark}

\begin{lemma} \label{lem:GGMB}
Let $G$ be a $k$-group scheme of finite type and let $X$ be a quasi-affine $k$-scheme of finite type equipped with a $G$-action.

\benum
    \item \label{m:GGMB-i} Let $G_x \subset G$ be the schematic stabilizer of an $x \in X(k)$. If every $\overline{k}$-torus of $G_{\overline{k}}$ stabilizes $x$ \up{that is, lies in $(G_x)_{\overline{k}}$}, then the orbit $G \cdot x \cong G / G_x$ is closed in $X$.
    \item \label{m:GGMB-ii} If every $\overline{k}$-torus of $G_{\overline{k}}$ lies in $(G^{\mathrm{gred}})_{\overline{k}}$ \up{see Remark~\uref{rem:condition-star}}, then each $G^{(i)} / G^{(i + 1)} \subseteq \mathrm{Res}_{k_i / k}(Q_i)$ in \uS\uref{pp:gred-filtration} is a closed immersion, the $G^{(i)} / G^{(i')}$ for $i' \ge i$ are all affine, and $G / G^{\mathrm{gred}}$ is affine.
\eenum
\end{lemma}

\begin{proof}
The claim \ref{m:GGMB-ii} follows from \ref{m:GGMB-i} and the affineness of $\mathrm{Res}_{k_i / k}(Q_i)$, so we focus on \ref{m:GGMB-i}. In the latter, the orbit $G \cdot x$ is a subscheme of $X$, see, for instance, \cite{GGMB14}*{remarques~2.1.4~(i)}. To show that it is closed, we proceed similarly to \cite{GGMB14}*{lemme~2.4.7} that treated the case when $G$ is affine. We first replace $k$ by $\overline{k}$ and $G$ by $G^{\mathrm{red}}$ to reduce to when $k = \overline{k}$ and $G$ is smooth. We then consider finitely many translates of $G^0 \cdot x$ to reduce further to when $G$ is also connected. Moreover, we replace $X$ by the schematic closure of $G \cdot x$ to assume that $G \cdot x$ is schematically dense in $X$. 

By the anti-Chevalley theorem \cite{CGP15}*{Theorem A.3.9}, then $G$ is an extension of a smooth affine $k$-group $G^{\mathrm{aff}}$ by a semi-abelian variety $A$. By assumption, the toral part $T$ of $A$ fixes $x$, so, since it is also normal in $G$, it fixes every $k$-point of $G \cdot x$. Since $k$-points of $G \cdot x$ are schematically dense in $X$, we find that $T$ acts trivially on $X$, in other words, the action of $G$ factors through $G / T$. Moreover, by \cite{CGP15}*{Proposition C.4.5 (2)}, the maximal $k$-tori of $G / T$ are precisely the images of the maximal tori of $G$. Thus, by replacing $G$ by $G / T$, we retain our assumption about $G_x$, and hence reduce to the case when $A$ is an abelian variety. However, $X$ is quasi-affine and $A$ is normal in $G$, so we likewise find that $A$ fixes every $k$-point of $G \cdot x$. It then acts trivially on $X$, so we may replace $G$ by $G^{\mathrm{aff}}$ to reduce to when $G$ is affine. 

In the affine case, the $k$-subgroup $G^{\mathrm{tor}} \le G$ generated by the $k$-tori of $G$ is connected, smooth, and normal with $G / G^{\mathrm{tor}}$ unipotent (see \cite{CGP15}*{Proposition A.2.11}). By our assumption on $G_x$ and the argument with $k$-points as above, $G^{\mathrm{tor}}$ acts trivially on $X$. We may therefore pass to $G / G^{\mathrm{tor}}$ to reduce to the case when $G$ is unipotent. Then, however, the desired closedness of $G \cdot x \subset X$ is the Rosenlicht lemma \cite{SGA3II}*{expos\'e XVII, lemme 5.7.3}. 
\end{proof}

\begin{remark}
\Cref{lem:GGMB} fails if $G$ is merely locally of finite type (and not quasi-compact). For instance, for any $a \in k^\times$, the constant $k$-group $\mathbb{Z}$ acts on $\mathbb{A}^1_k$ by $(n, t) \mapsto a^n t$, where $t$ is the coordinate of $\mathbb{A}^1_k$; if $a$ is not a root of unity, then the orbit of the $k$-point $t = 1$ is not closed in $\mathbb{A}^1_k$.
\end{remark}

\begin{remark}\label{rem:condition-star}
By \cite{GGMB14}*{lemme 2.4.5}, \cite{SGA3IIInew}*{expos\'e XVII, th\'eor\`eme 7.3.1} (applied to $(G^{\mathrm{lin}})^0$), and \cite{SGA3II}*{expos\'e XII, proposition 1.12}, for a $k$-group scheme $G$ locally of finite type, every $\overline{k}$-torus of $G_{\overline{k}}$ lies in $(G^{\mathrm{gred}})_{\overline{k}}$, that is, $(G_{\overline{k}})^{\mathrm{tor}} \le (G^{\mathrm{gred}})_{\overline{k}}$,

\benumr
    \item \label{m:CS-i} if $G^0$ is a normal subgroup of smooth $k$-group scheme (in particular, if $G$ is smooth); or
    \item \label{m:CS-ii} if $G^0$ is nilpotent (for instance, either commutative or unipotent, see \cite{SGA3IIInew}*{expos\'e~XVII, corollaire 3.7}).
\eenum

For a $k$-group scheme locally of finite type, the condition that every $\overline{k}$-torus of $G_{\overline{k}}$ lies in $(G^{\mathrm{gred}})_{\overline{k}}$ is stable under base change to any separable field extension $k^{\prime}/k$ (see \S\ref{pp:fundamental-filtration}~\ref{m:Gtor} and \S\ref{pp:Xgred}) but is not in general inherited by inner forms, see \cite{FG21}*{Example 7.2 (b) and Remark 7.3 (b)}. As an example, if a pseudo-finite $k$-group $G$ satisfies this condition, then \S\ref{pp:fundamental-filtration}~\ref{m:Gtor} implies that $(G^0_{\overline{k}})^{\mathrm{red}}$ is unipotent (compare with \Cref{cor:pfinite-affine} above).
\end{remark}

\csub[The comparison map $i_{G,\, \overline{G}}$] \label{sec:iG}

We approach pseudo-reductive groups (resp.,~and pseudo-abelian varieties to some extent) via a comparison map $i_{G,\, \overline{G}}$ that relates them to restrictions of scalars of reductive groups (resp.,~of abelian varieties). We now analyze this map in an abstract setting and establish \Cref{prop:ker-pseudo-finite,prop:Brauer-affine} that allow us to control its kernel and ``cokernel.'' Although the map $i_{G,\, \overline{G}}$ is known to play a prominent role in the theory of pseudo-reductive groups, even in this case the affineness of its ``cokernel'' established in \Cref{prop:Brauer-affine} seems new.

\bpp[The comparison map $i_{G,\, \overline{G}}$ and its kernel] \label{pp:iG}
Let $G$ be a $k$-group scheme locally of finite type, let $k^\prime/k$ be a finite field extension, and let $q\colon G_{k^\prime} \twoheadrightarrow \overline{G}$ be a $k^\prime$-group scheme quotient. By the universal property, $q$ corresponds to a $k$-group scheme homomorphism
$$
i_{G,\, \overline{G}} \colon G \rightarrow \mathrm{Res}_{k^\prime/k}(\overline{G})
$$
(see \S\ref{pp:conv} and \S\ref{pp:k-groups-lft}). The corresponding $k^\prime$-group homomorphism $\mathrm{Ker}(i_{G, \,\overline{G}})_{k^\prime} \rightarrow \overline{G}$ is trivial, so $\mathrm{Ker}(i_{G, \,\overline{G}})_{k^\prime} \lhd \mathrm{Ker}(q)$. Even if the choice of $k^\prime$ appears to be noncanonical, this does not matter: enlarging $k^\prime$ and using the corresponding base change of $\overline{G}$ amounts to postcomposing $i_{G,\,\overline{G}}$ with a closed immersion, see \cite{BLR90}*{bottom of p.~197}.
\epp

To analyze the kernel of $i_{G,\, \overline{G}}$ in \Cref{prop:ker-pseudo-finite}, we use the following lemma.

\begin{lemma} \label{lem:res-unipotent}
For a finite, purely inseparable field extension $k^\prime/k$ and a $k^\prime$-group scheme $G$ locally of finite type, the counit map $(\mathrm{Res}_{k^\prime/k}(G))_{k^\prime} \twoheadrightarrow G$ has unipotent kernel that is split if $G$ is smooth.
\end{lemma}

\begin{proof}
The counit map is given by the functoriality of the restriction of scalars relative to the diagonal map $k^\prime \otimes_k k^\prime \twoheadrightarrow k^\prime$. Since $k^\prime/k$ is purely inseparable, this diagonal map is simply the quotient of a local Artinian $k^\prime$-algebra $A \colonequals k^\prime \otimes_k k^\prime$ by its maximal ideal $\mathfrak{m}$. Thus, by filtering by powers of $\mathfrak{m}$ and using \S\ref{pp:conv} and formal smoothness, it suffices to show that the kernel of the map $\mathrm{Res}_{(A/\mathfrak{m}^{n + 1})/k^\prime}(G_{A/\mathfrak{m}^{n + 1}}) \rightarrow \mathrm{Res}_{(A/\mathfrak{m}^n)/k^\prime}(G_{A/\mathfrak{m}^n})$ is a power of $\mathbb{G}_{a,\, k^\prime}$. This, however, follows from deformation theory \cite{SGA3Inew}*{expos\'e III, th\'eor\`eme 0.1.8}. 
\end{proof}

\begin{proposition} \label{prop:ker-pseudo-finite}
In the setting of \uS\uref{pp:iG}, 
\benum
    \item \label{m:KPF-a} If $\mathrm{Ker}(q)$ is quasi-compact \up{e.g.,~if so is $G$}, then both $i_{G,\, \overline{G}}$ and $\mathrm{Ker}(i_{G,\, \overline{G}})$ are quasi-compact;
    \item \label{m:KPF-b} If $\mathrm{Ker}(q)^0$ \up{resp.,~$\mathrm{Ker}(q)$} is unipotent, then so is $\mathrm{Ker}(i_{G,\, \overline{G}})^0$ \up{resp., $\mathrm{Ker}(i_{G,\, \overline{G}})$};
    \item \label{m:KPF-c} If $\mathrm{Ker}(q)^0$ is unipotent \up{resp.,~and $\mathrm{Ker}(q)$ is quasi-compact} and $G^{\mathrm{gred}}$ has no nontrivial, connected, smooth, unipotent, normal $k$-subgroups, then $\mathrm{Ker}(i_{G,\, \overline{G}})$ is pseudo-\up{locally quasi-finite} \up{resp.,~pseudo-finite}.
\end{enumerate}
\end{proposition}

\begin{proof} \hfill
\benum
    \item By \Cref{lem:res-unipotent}, the kernels of the horizontal maps in the commutative square  
\[
\qq \xymatrix{
(\mathrm{Res}_{k^\prime/k}(G_{k^\prime}))_{k^\prime} \ar@{->>}[r] \ar[d] & G_{k^\prime} \ar@{->>}[d]^-{q} \\
(\mathrm{Res}_{k^\prime/k}(\overline{G}))_{k^\prime} \ar@{->>}[r] & \overline{G}
}
\]
    are unipotent, so these maps are affine. The quasi-compactness of $\mathrm{Ker}(q)$ then ensures that the vertical maps are quasi-compact. Thus, since $G \hookrightarrow \mathrm{Res}_{k^\prime/k}(G_{k^\prime})$ is a closed immersion by \cite{BLR90}*{bottom of p.~197}, we conclude that $i_{G,\, \overline{G}}$ is quasi-compact. By base change, then $\mathrm{Ker}(i_{G,\, \overline{G}})$ is also quasi-compact.
    
    \item It suffices to recall that $\mathrm{Ker}(i_{G, \,\overline{G}})_{k^\prime} \lhd \mathrm{Ker}(q)$ and to review \S\ref{pp:conv}.
    
    \item Since $\mathrm{Ker}(i_{G,\, \overline{G}})$ is normal in $G$ and of formation compatible with base change to separable extensions, $G(k^s)$-conjugation preserves $\mathrm{Ker}(i_{G,\, \overline{G}})^{\mathrm{gred}}$. Thus, the latter is normal in $G^{\mathrm{gred}}$ (see \S\ref{pp:fundamental-filtration}~\ref{m:Ggred}). Our assumption on $G^{\mathrm{gred}}$ and the unipotence of $\mathrm{Ker}(i_{G,\, \overline{G}})^0$ given by \ref{m:KPF-b} now ensure that $(\mathrm{Ker}(i_{G,\, \overline{G}})^{\mathrm{gred}})^0$ is trivial, to the effect that $\mathrm{Ker}(i_{G,\, \overline{G}})$ is pseudo-(locally quasi-finite). If $\mathrm{Ker}(q)$ is quasi-compact, then, by \ref{m:KPF-a}, so is $\mathrm{Ker}(i_{G,\, \overline{G}})$, so that it is even pseudo-finite (see \Cref{def:pseudo-finite}). \qedhere
\end{enumerate}
\end{proof}

To analyze the ``cokernel'' of $i_{G,\, \overline{G}}$ in \Cref{prop:Brauer-affine}, we use the following lemma. 

\begin{lemma} \label{lem:unip-quotient}
Let $G$ be a $k$-group scheme locally of finite type with $G^0$ an extension of a solvable, affine $k$-group by an infinitesimal $k$-group. For any closed $k$-subgroup $H \le G$, the connected components of $G/H$ are clopen and affine, in particular, if $G$ is of finite type, then $G/H$ is affine.
\end{lemma}

\begin{proof}
The key aspect is the affineness: indeed, $G/H$ is a separated $k$-scheme locally of finite type by the results reviewed in \S\ref{pp:k-groups-lft}, and so its connected components are clopen by \cite{SP}*{Lemmas \href{https://stacks.math.columbia.edu/tag/04MF}{04MF} and \href{https://stacks.math.columbia.edu/tag/04ME}{04ME}}. By the end of \S\ref{pp:perfect-case}, for a sufficiently large Frobenius kernel $\alpha_G \lhd G$, both $\overline{G} := G/\alpha_G$ and $\overline{H} := H/(\alpha_G \cap H)$ are smooth $k$-groups and, by our assumption, $\overline{G}^0$ is solvable. Moreover, the flat, separated surjection $G/H \twoheadrightarrow \overline{G}/\overline{H}$ is finite: indeed, its base change along itself admits a surjection from $\alpha_G \times_k G/H$, so \cite{SP}*{Lemmas \href{https://stacks.math.columbia.edu/tag/0AH6}{0AH6} and \href{https://stacks.math.columbia.edu/tag/02LS}{02LS}} give the claimed finiteness. Thus, it suffices to show that every connected component of $\overline{G}/\overline{H}$ is affine. In effect, we may assume that $G^0$ is solvable, affine, and that $G$ and $H$ are smooth, so that so is $G/H$. 

The map $G/H^0 \rightarrow G/H$ is an $H/H^0$-torsor, so, by \S\ref{pp:fundamental-filtration}~\ref{m:G0} and \cite{SGA3II}*{expos\'e X, corollaire~5.14}, the connected components of $G/H^0$ are finite \'etale over those of $G/H$. By \cite{SP}*{Lemma \href{https://stacks.math.columbia.edu/tag/01YQ}{01YQ}}, it then suffices to argue that every connected component of $G/H^0$ is affine, that is, we may assume that $H$ is connected. Once $H$ is connected, the connected components of $G/H$ are precisely the images of the connected components of $G$. Moreover, each component of $G$ becomes a finite union of $G^0$-torsors after base change to some finite separable extension of $k$. Thus, effectivity of descent for affine schemes \cite{SP}*{Lemma \href{https://stacks.math.columbia.edu/tag/0245}{0245}} allows us to replace $G$ by $G^0$ to reduce to the case when $G$ is also connected. Once both $G$ and $H$ are smooth and connected, since we may base change to $\overline{k}$ (see \cite{SP}*{\href{https://stacks.math.columbia.edu/tag/02L5}{Lemma~02L5}}), the affineness is a special case of \cite{Con15b}*{Example 5.5} and \cite{Bri21}*{Theorem~1}. 
\end{proof}

\begin{proposition} \label{prop:Brauer-affine}
In the setting of \uS\uref{pp:iG}, suppose that $k^\prime/k$ is purely inseparable and that $\mathrm{Ker}(q)$ is quasi-compact. Then
$$
Q_{G,\,\overline{G}} := \mathrm{Res}_{k^\prime/k}(\overline{G})/i_{G,\,\overline{G}}(G)
$$
is an affine scheme, and if $\overline{G}$ is commutative, then $Q_{G,\,\overline{G}}$ is a commutative, unipotent $k$-group.
\end{proposition}

\begin{proof}
When $q$ is an isomorphism, a similar result was established over any base in \cite{brauer-purity}*{Lemma~2.1}, so we will build on its method. Since $i_{G,\, \overline{G}}$ is quasi-compact by \Cref{prop:ker-pseudo-finite}~\ref{m:KPF-a}, arguments as in the proof of \Cref{lem:unip-quotient} ensure that $Q_{G,\,\overline{G}}$ is a separated $k$-scheme locally of finite type. 

By descent for affineness \cite{SP}*{Lemma~\href{https://stacks.math.columbia.edu/tag/02L5}{02L5}} and \S\ref{pp:conv}, it is enough to show that $(Q_{G,\,\overline{G}})_{k^\prime}$ is affine (resp.,~and is unipotent when $\overline{G}$ is commutative). Over $k^\prime$ we have the counit map 
$$
(\mathrm{Res}_{k^\prime/k}(\overline{G}))_{k^\prime} \cong \mathrm{Res}_{k^\prime\otimes_k k^\prime/k^\prime}(\overline{G}_{k^\prime\otimes_k k^\prime}) \xrightarrow{h} \mathrm{Res}_{k^\prime/k^\prime}(\overline{G}) \cong \overline{G},
$$
as well as a diagram of $k^\prime$-group homomorphisms
\[
\begin{xy}
\xymatrix{
1 \ar[r] & \mathrm{Ker}(q) \ar[r] \ar[d]^{i} & G_{k^\prime} \ar[r]^{q} \ar[d]^{(i_{G,\,\overline{G}})_{k^\prime}} & \overline{G} \ar[r] \ar@{=}[d] & 1 \\
1 \ar[r] & \mathrm{Ker}(h) \ar[r] & (\mathrm{Res}_{k^\prime/k}(\overline{G}))_{k^\prime} \ar[r]^-{h} & \overline{G} \ar[r] & 1,
}
\end{xy}
\]
whose commutativity follows from the case when $q$ is the identity. We conclude that $h$ is surjective and, by a diagram chase, that
\begin{equation} \label{eqn:snake}
\mathrm{Ker}(h)/i(\mathrm{Ker}(q)) \cong (Q_{G,\,\overline{G}})_{k^\prime}.
\end{equation}
Moreover, $\mathrm{Ker}(h)$ is unipotent by \Cref{lem:res-unipotent} and $i$ inherits quasi-compactness from $i_{G,\,\overline{G}}$. At this point, by \eqref{eqn:snake} and \S\ref{pp:conv}, \Cref{lem:unip-quotient} gives the claim.
\end{proof}

\csub[From pseudo-abelian varieties to abelian varieties]\label{sec:pseudo-abelian}

Pseudo-abelian varieties were introduced by Totaro in \cite{Tot13}. In this section, we review the aspects of their theory relevant for us, stressing the comparison map $i_G$ of \S\ref{sec:iG} because this highlights the analogy with the theory of pseudo-reductive groups.

\bpp[Pseudo-abelian varieties] \label{pp:pseudo-abelian}
A \emph{pseudo-abelian variety} is a smooth, connected $k$-group scheme $G$ that has no nontrivial connected, smooth, affine, normal $k$-subgroups, see \cite{Tot13}*{Definition 0.1}. By \cite{Tot13}*{Theorem 2.1}, a pseudo-abelian variety $G$ is automatically commutative, in particular, the normality assumption in the definition could be dropped: $G$ has no nontrivial connected, smooth, affine $k$-subgroup, in other words, $G^{\mathrm{sm},\, \mathrm{lin}} = 1$. In addition, by \emph{loc.\ cit.}, a pseudo-abelian variety $G$ is in a unique way an extension
\begin{equation}
0 \rightarrow G^{\mathrm{ab}} \rightarrow G \rightarrow G^{\mathrm{unip}} \rightarrow 0
\label{eqn:pseudo-abelian}
\end{equation}
of a connected, smooth, commutative, unipotent $k$-group $G^{\mathrm{unip}}$ by an abelian variety $G^{\mathrm{ab}}$. In particular, $G$ is proper if and only if it is an abelian variety.

On the other hand, \S\ref{pp:fundamental-filtration} gives us an extension in the other order:
\begin{equation}
0 \rightarrow G^{\mathrm{lin}} \rightarrow G \rightarrow G^{\mathrm{av}} \rightarrow 0.
\label{eqn:Chevalley}
\end{equation}

The definition and the commutativity of pseudo-abelian varieties imply that the connected, affine $k$-group $G^{\mathrm{lin}}$ has no positive-dimensional smooth $k$-subgroups, so $G^{\mathrm{lin}}$ is pseudo-finite, see \S\ref{pp:Xgred}. Since there are no nontrivial homomorphisms from a torus to a unipotent group, nor to an abelian variety, \eqref{eqn:pseudo-abelian} and \S\ref{pp:fundamental-filtration} imply that the connected, smooth, affine, commutative $\overline{k}$-group $((G^{\mathrm{lin}})_{\overline{k}})^{\mathrm{red}}$ is unipotent and is identified with $(G_{\overline{k}})^{\mathrm{lin}}$.

Both extensions \eqref{eqn:pseudo-abelian} and \eqref{eqn:Chevalley} help in practice, and it is useful to keep in mind the relationship between them. Namely, the map of abelian varieties
\[
G^{\mathrm{ab}} \twoheadrightarrow G^{\mathrm{av}}
\]
is an isogeny: it is surjective because the affine group $G^{\mathrm{unip}}$ cannot surject on a nonzero abelian variety, and its kernel $G^{\mathrm{ab}} \cap G^{\mathrm{lin}}$ is finite because it is both proper and affine. Similarly,
\[
G^{\mathrm{lin}} \twoheadrightarrow G^{\mathrm{unip}}
\]
is an isogeny: it is surjective because the abelian variety $G^{\mathrm{av}}$ cannot surject on a positive-dimensional affine group, and its kernel $G^{\mathrm{ab}} \cap G^{\mathrm{lin}}$ is finite. The order of the common kernel $G^{\mathrm{ab}} \cap G^{\mathrm{lin}}$, and hence also the common degree, of these two isogenies is a power of the characteristic exponent $p$ of $k$ because the component group of $(G^{\mathrm{ab}} \cap G^{\mathrm{lin}})_{\overline{k}}$ inherits unipotence from $((G^{\mathrm{lin}})_{\overline{k}})^{\mathrm{red}}$, see \S\ref{pp:conv}.

By \cite{Tot13}*{Corollary 6.5}, if $k$ is imperfect of characteristic $p$, then every connected, smooth, commutative, unipotent $k$-group of exponent $p$ (such as $\mathbb{G}_a$) occurs as $G^{\mathrm{unip}}$ for some pseudo-abelian variety $G$ over $k$. In particular, the isogenies $G^{\mathrm{lin}} \twoheadrightarrow G^{\mathrm{unip}}$ exhibit many pseudo-finite, commutative covers of split unipotent $k$-groups, see also \cite{Tot13}*{Lemma 6.3}.
\epp

\bpp[The comparison map $i_G$ for pseudo-abelian varieties] \label{pp:pav-iG} 
Let $G$ be a pseudo-abelian variety over $k$. Since $(-)^{\mathrm{sm},\, \mathrm{lin}}$ and $(-)^{\mathrm{lin}}$ agree over perfect fields (see \S\ref{pp:perfect-case}~\ref{m:PC-i}), by \cite{EGAIV2}*{corollaire~4.8.11}, there is the smallest finite, purely inseparable field extension $k^\prime/k$ such that $(G_{k^\prime})^{\mathrm{lin}}$ is smooth, in other words, such that it is a $k^\prime$-descent of the smooth, unipotent $\overline{k}$-subgroup $(G_{\overline{k}})^{\mathrm{lin}}$. Similarly to \eqref{eqn:Chevalley}, we set $\overline{G} \colonequals (G_{k^\prime})^{\mathrm{av}} \cong G_{k^\prime}/(G_{k^\prime})^{\mathrm{lin}}$ and consider the resulting comparison map
\[
i_G\colon G \rightarrow \mathrm{Res}_{k^\prime/k}(\overline{G}),
\]
whose formation commutes with base change to any separable extension of $k$. The image $i_G(G)$ is smooth and connected, so it is a pseudo-abelian subvariety of the pseudo-abelian variety $\mathrm{Res}_{k^\prime/k}(\overline{G})$, see the final part of \S\ref{pp:fundamental-filtration}. The $k^\prime$ and $\overline{G}$ remain the same for $i(G)$ and the map $i_{i(G)}$ is the inclusion $i_G(G) \subset \mathrm{Res}_{k^\prime/k}(\overline{G})$, compare with the proof of \cite{CGP15}*{Lemma 9.2.1}. \Cref{prop:ker-pseudo-finite,prop:Brauer-affine} ensure that $\mathrm{Ker}(i_G)$ is pseudo-finite and unipotent, while $\mathrm{Coker}(i_G)$ is connected, smooth, and unipotent. With this control of the (co)kernel, it is fruitful to study the pseudo-abelian variety $G$ by reducing to the pseudo-abelian variety $i_G(G)$ and then to the abelian variety $\overline{G}$.
\epp

\csub[From quasi-reductive groups to reductive groups]\label{sec:pseudo-reductive}

Pseudo-reductive groups were introduced by Tits \cite{Tit13}*{cours 1991--1992} and have been studied extensively by him and, more recently, by Conrad--Gabber--Prasad \cite{CGP15} and Conrad--Prasad \cite{CP16}, \cite{CP17}. We now frame the pseudo-reductive theory around the comparison map of \S\ref{sec:iG} and review some aspects needed for our later arguments.

\bpp[Notation] 
Throughout this section, we fix a connected, smooth, affine $k$-group $G$.
\epp

\bpp[The comparison map $i_G$] \label{pp:pred-iG}
We may choose the smallest finite, purely inseparable field extension $k^\prime/k$ such that $\mathscr{R}_{\mathrm{u},\, k^\prime}(G_{k^\prime})$ is a $k^\prime$-descent of $\mathscr{R}_{\mathrm{u},\,\overline{k}}(G_{\overline{k}})$ (see \cite{EGAIV2}*{corollaire 4.8.11} for the existence of \emph{the} smallest such $k^\prime$), in other words, such that the $k^\prime$-group $\overline{G} \colonequals (G_{k^\prime})^{\mathrm{pred}}$ of \eqref{eqn:Gpred} is reductive (recall from \S\ref{pp:perfect-case}~\ref{m:PC-ii} that the unipotent radical over a perfect field descends the geometric unipotent radical). We say that $G$ is \emph{simply connected} (resp.,~\emph{adjoint}) if $\overline{G}$ is semisimple and simply connected (resp.,~adjoint). In general, we consider the resulting comparison map
\[
i_G\colon G \rightarrow \mathrm{Res}_{k^\prime/k}(\overline{G}),
\]
whose formation commutes with base change to any separable extension of $k$. \Cref{prop:ker-pseudo-finite} ensures that $\mathrm{Ker}(i_G)$ is unipotent and that it is pseudo-finite if and only if $G$ is pseudo-reductive. In fact, for pseudo-reductive $G$, as we will review in \S\ref{pp:pseudo-minimal} below, $\mathrm{Ker}(i_G)$ is a central extension of \emph{commutative}, pseudo-finite, unipotent $k$-groups (and is even itself commutative if $\mathrm{char}(k) \neq 2$). Even if $G$ is pseudo-reductive, $\mathrm{Ker}(i_G)$ need not be connected or smooth, see \cite{CGP15}*{Example~1.6.3, Remark 9.1.11, Theorem 9.8.1 (3)--(4)}.

By \cite{CGP15}*{Lemma 9.2.1}, the image $i_G(G)$ is pseudo-reductive, the $k^\prime$ and $\overline{G}$ remain the same for $i_G(G)$, and the map $i_{i_G(G)}$ is the inclusion $i_G(G) \le \mathrm{Res}_{k^\prime/k}(\overline{G})$. By \Cref{prop:Brauer-affine}, the homogeneous space $\mathrm{Res}_{k^\prime/k}(\overline{G})/i_G(G)$ is connected, smooth, and affine.

With this control of the kernel and the ``cokernel,'' it is often fruitful to approach the study of $G$ by first reducing to the pseudo-reductive group $i_G(G)$ and then to the reductive group $\overline{G}$.
\epp

\bpp[Passage to pseudo-reductive groups of minimal type]\label{pp:pseudo-minimal}
Suppose now that $G$ is pseudo-reductive. By \cite{CGP15}*{Proposition 9.4.2~(i)}, for a maximal $k$-torus $T \subset G$ (which exists, see~\S\ref{pp:fundamental-filtration}~\ref{m:Ggred}), the subgroup
\[
\mathscr{C}_G := \mathrm{Ker}(i_G) \cap Z_G(T) \subset \mathrm{Ker}(i_G)
\]
is central in $G$ and does not depend on the choice of $T$. Of course, $\mathscr{C}_G$ is also pseudo-finite and unipotent because it inherits these properties from $\mathrm{Ker}(i_G)$. By \cite{CGP15}*{Proposition 9.4.2~(i)--(ii)} and \cite{CGP15}*{Corollary 9.4.3}, the quotient $\mathscr{G} := G/\mathscr{C}_G$ is pseudo-reductive with $\mathscr{C}_{\mathscr{G}} = 1$, and the map $i_G$ is the composition of $i_{\mathscr{G}}$ and the quotient $G \twoheadrightarrow \mathscr{G}$. A pseudo-reductive $k$-group $G$ is
\begin{itemize}
    \item \emph{of minimal type} if $\mathscr{C}_G = 1$ (see \cite{CGP15}*{Definition 9.4.4}), equivalently, if $G \cong \mathscr{G}$;
    \item \emph{ultraminimal} if $i_G$ is injective, equivalently, if $G$ is of minimal type and $G_{k^s}$ has a reduced root system (see \cite{CGP15}*{Theorem 9.4.7}).
\end{itemize}

For instance, by the above, $i_G(G)$ is an ultraminimal pseudo-reductive $k$-group. Pseudo-reductive groups of minimal type that are not ultraminimal exist only in characteristic $2$, and they cause by far the most complications in the general theory; for instance, they are the \emph{raison d'\^etre} for the book \cite{CP16}. In practice, since $\mathscr{C}_G$ is central, it is often straightforward to reduce to the minimal type case, but passing to ultraminimal groups tends to be more delicate.

One may control $\mathrm{Ker}(i_G)$ even if $G$ is of minimal type but not ultraminimal. Namely, in this case, by \cite{CGP15}*{Theorem 9.4.7}, the pseudo-finite, unipotent $k$-group $\mathrm{Ker}(i_G)$ is connected, commutative, but not central, and if $G$ contains a \emph{split} maximal $k$-torus $T$, then $\mathrm{Ker}(i_G)$ is the direct product of its intersections with the root groups of the multipliable roots of $T$.
\epp

\csub[The control of pseudo-parabolic subgroups]\label{sec:pseudo-parabolic}

The study of connected, smooth, affine $k$-groups benefits from the theory of pseudo-parabolic subgroups introduced by Borel--Tits in \cite{BT78}, as these are means for passing to smaller pseudo-reductive groups, see \S\ref{pp:fundamental-filtration}~\ref{m:pseudo-parabolic}. We review some aspects of the pseudo-parabolic subgroup theory for use in subsequent chapters.

\bpp[Subgroups associated to $\mathbb{G}_m$-actions]\label{pp:Gm-actions}
Let $A$ be a ring and let $G$ be a finitely presented, affine $A$-group scheme equipped with a left action of $\mathbb{G}_{m,\, A}$ over $A$:
\[
\lambda\colon \mathbb{G}_{m,\,A} \times G \rightarrow G.
\]
A common case is when the action is conjugation composed with an $A$-homomorphism $\mathbb{G}_{m,\, A} \rightarrow G$ (often also denoted $\lambda$ for simplicity). By \cite{CGP15}*{Remark 2.1.11 (with Lemma 2.1.5)}, the $\mathbb{G}_{m,\,A}$-action gives rise to the finitely presented, closed $A$-subgroups:
\benuma
    \item $P_G(\lambda) \le G$, the \emph{attractor} subgroup, i.e., the subfunctor parametrizing those $A$-algebra valued points $g$ of $G$ for which the action map $t \mapsto t \cdot g$ extends to a map $\mathbb{A}^1 \rightarrow G$;
    \item $U_G(\lambda) \lhd P_G(\lambda)$, the subfunctor parametrizing those such $g$ for which the resulting map $\mathbb{A}^1 \rightarrow G$ sends the zero section of $\mathbb{A}^1$ to the identity section of $G$;
    \item $Z_G(\lambda) \le P_G(\lambda)$, the subfunctor parametrizing the $\lambda$-stable sections of $G$.
\end{enumerate}

The formation of $P_G(\lambda)$, $U_G(\lambda)$, and $Z_G(\lambda)$ commutes with arbitrary base change in $A$, as well as with intersections with $\lambda$-stable, closed, finitely presented $A$-subgroups of $G$, and, by \cite{CGP15}*{Remark~2.1.11 (with Proposition 2.1.8 (2))}, we have
\[
P_G(\lambda) \cong U_G(\lambda) \rtimes Z_G(\lambda).
\]
By \cite{CGP15}*{Remark 2.1.11 (with Lemma 2.1.5 and Proposition 2.1.8 (4))}, the $A$-group $U_G(\lambda)$ has connected, unipotent $A$-fibers, and if the $A$-fibers of $G$ are connected, then so are those of $P_G(\lambda)$ and $Z_G(\lambda)$. By \cite{CGP15}*{Remark 2.1.11 (with Propositions 2.1.8 (3) and 2.1.10)}, if $G$ is $A$-smooth, then so are $P_G(\lambda)$, $U_G(\lambda)$, and $Z_G(\lambda)$, and the $A$-fibers of $U_G(\lambda)$ are even split unipotent.

Letting $-\lambda$ denote the $\mathbb{G}_{m,\, A}$-action on $G$ obtained from $\lambda$ by precomposing with the inversion of $\mathbb{G}_{m,\, A}$, by \emph{loc.~cit.} and \cite{EGAIV4}*{corollaire 17.9.5}, the multiplication map
\be\label{eqn:open}
U_G(-\lambda) \times P_G(\lambda) \to G
\ee
is an open immersion provided that its source is $A$-flat (which holds if either $G$ is $A$-smooth or $A$ is a field). If \eqref{eqn:open} is an open immersion and the reduced geometric $A$-fibers of $G$ are connected and solvable, then \eqref{eqn:open} is even an isomorphism: indeed, we may assume that $A$ is a field and $G$ is smooth, note that $G/P_G(\lambda)$ is affine by \cite{Con15b}*{Example 5.5} and \cite{Bri21}*{Theorem 1}, and then conclude from \Cref{lem:GGMB}~\ref{m:GGMB-i} that the open $U_G(-\lambda) \subset G/P_G(\lambda)$ is simultaneously closed (this argument reproves \cite{CGP15}*{Proposition 2.1.12 (1)}).

The formation of the open immersion \eqref{eqn:open} commutes with intersections with $\lambda$-stable, closed, finitely presented $A$-subgroups $G' \le G$ for which $U_{G'}(-\lambda) \times P_{G'}(\lambda)$ is $A$-flat. Similarly, the $P_G(\lambda)$, $U_G(\lambda)$, and $Z_G(\lambda)$ are compatible with flat surjections as follows: by \cite{CGP15}*{Remark 2.1.11 (with Corollary 2.1.9)} and \cite{EGAIV3}*{corollaire 11.3.11}, if $G$ has connected $A$-fibers, $P_G(\lambda)$ (resp., $U_G(\lambda)$; resp., $Z_G(\lambda)$) is $A$-flat, and $\pi\colon G \twoheadrightarrow G'$ is a faithfully flat, $\mathbb{G}_{m,\, A}$-equivariant surjection onto a finitely presented, affine $A$-group $G'$ with a $\mathbb{G}_{m,\, A}$-action $\lambda'$, then $\pi$ induces a flat surjection $P_G(\lambda) \twoheadrightarrow P_{G'}(\lambda')$ (resp., $U_G(\lambda) \twoheadrightarrow U_{G'}(\lambda')$; resp., $Z_G(\lambda) \twoheadrightarrow Z_{G'}(\lambda')$) whose target is also $A$-flat.
\epp

\bpp[Pseudo-parabolic subgroups]\label{pp:pseudo-parabolic}
For a connected, smooth, affine $k$-group scheme $G$, a $k$-subgroup $P \le G$ is \emph{pseudo-parabolic} if it has the form
\[
P = P_\lambda \colonequals P_G(\lambda)\mathscr{R}_{\mathrm{u},\,k}(G), 
\quad\quad \text{equivalently,} \quad\quad 
P = P_\lambda \colonequals P_G(\lambda)\mathscr{R}_{\mathrm{us},\,k}(G)
\]
for a $k$-homomorphism $\lambda\colon\mathbb{G}_{m,\,k}\to G$, see \cite{CGP15}*{Definition 2.2.1 and Proposition 2.2.4}. 

By \S\ref{pp:Gm-actions}, pseudo-parabolic $k$-subgroups are connected, smooth, and affine. Moreover, by \cite{CGP15}*{Corollary 2.2.5}, we have
\[
\mathscr{R}_{\mathrm{u},\,k}(P_\lambda) = U_G(\lambda)\mathscr{R}_{\mathrm{u},\,k}(G) 
\quad \text{and} \quad 
\mathscr{R}_{\mathrm{us},\,k}(P_\lambda) = U_G(\lambda)\mathscr{R}_{\mathrm{us},\,k}(G),
\]
so that $\mathscr{R}_{\mathrm{us},\,k}(P_\lambda) = U_G(\lambda)$ when $G$ is quasi-reductive, and $\mathscr{R}_{\mathrm{u},\,k}(P_\lambda)$ is split unipotent when $G$ is pseudo-reductive. By \cite{CGP15}*{Proposition 3.5.2 (1)} (with the final aspect of \S\ref{pp:fundamental-filtration}~\ref{m:Ggred}), a $k$-subgroup $P \le G$ is pseudo-parabolic if and only if $P_{k'} \le G_{k'}$ is pseudo-parabolic for some (equivalently, any) separable extension $k'/k$. A pseudo-parabolic $k$-subgroup $P \le G$ is a \emph{pseudo-Borel} if $P_{k^s} \le G_{k^s}$ is a minimal pseudo-parabolic of $G_{k^s}$, equivalently, by \cite{CGP15}*{Proposition 3.5.4}, if the image of $P_{\overline{k}}$ is a Borel subgroup of $(G_{\overline{k}})^{\mathrm{pred}}$. By \cite{CGP15}*{Proposition 2.2.9}, the pseudo-parabolic (resp., pseudo-Borel) $k$-subgroups of a reductive $k$-group are precisely its parabolic (resp., Borel) $k$-subgroups. By \cite{CP17}*{Corollary 4.3.5} (resp., \cite{CGP15}*{Proposition 3.5.8}), pseudo-parabolic $k$-subgroups of a pseudo-parabolic $P \le G$ (resp., smooth $k$-subgroups of $G$ containing $P$) are pseudo-parabolic. 

By \cite{CGP15}*{Proposition 2.2.10}, the pseudo-parabolic $k$-subgroups of $G$ are the preimages of the pseudo-parabolic $k$-subgroups of $G^{\mathrm{pred}}$ (so also of those of $G^{\mathrm{qred}}$). We may even replace $G^{\mathrm{pred}}$ by its pseudo-reductive quotient of minimal type (see \S\ref{pp:pseudo-minimal}): by \cite{CGP15}*{Proposition 2.2.12 (3)}, for every quotient $G \twoheadrightarrow G'$ of pseudo-reductive $k$-groups with a central kernel, the pseudo-parabolic $k$-subgroups of $G$ are the preimages of the pseudo-parabolic $k$-subgroups of $G'$. Of course, by definition, the pseudo-parabolic $k$-subgroups of a quasi-reductive $G$ are precisely the $P_G(\lambda)$ for $k$-homomorphisms $\lambda\colon\mathbb{G}_{m,\,k}\to G$. 

In the case when $G = \mathrm{Res}_{k'/k}(G')$ for a nonzero, finite, reduced $k$-algebra $k'$ and a $k'$-group $G'$ with pseudo-reductive $k'$-fibers, by \cite{CGP15}*{Proposition 2.2.13}, the pseudo-parabolic $k$-subgroups $P \le G$ correspond to the fiberwise pseudo-parabolic $k'$-subgroups $P' \le G'$ via the inverse bijections $P' \mapsto \mathrm{Res}_{k'/k}(P')$ and $P \mapsto \pi(P_{k'})$, where $\pi\colon G_{k'} \to G'$ is the counit of the adjunction. In this setting, by \cite{CGP15}*{Corollary A.5.4 (3)} (with \S\ref{pp:conv}), for the corresponding $P$ and $P'$, we have $G/P \cong \mathrm{Res}_{k'/k}(G'/P')$. To reach such restrictions of scalars in practice, we  use the following lemmas.
\epp

The following lemma pointed out to us by Gabber controls pseudo-parabolic $k$-subgroups when passing from a pseudo-reductive $k$-group of minimal type to its ultraminimal quotient. 

\begin{lemma} \label{lem:minimal-to-ultra}
Let $G$ be a pseudo-reductive $k$-group of minimal type, let $i_G \colon G \rightarrow \mathrm{Res}_{k^\prime/k}(\overline{G})$ be the map as in \uS\uref{pp:pred-iG}, and let $P \le G$ be a pseudo-parabolic $k$-subgroup. The $k$-subgroup $i_G(P) \le i_G(G)$ is pseudo-parabolic and, for a $k$-scheme $S$ and an $S$-point $y$ of $i_G(G)/i_G(P)$, the $y$-fiber of the map
\[
G/P \rightarrow i_G(G)/i_G(P)
\]
is a torsor under an $S$-form of the $k$-group $\mathrm{Ker}(i_G)/(\mathrm{Ker}(i_G) \cap P)$, and this $S$-form is trivial whenever $y$ lifts to an $S$-point of $i_G(G)$. Moreover, we have
\be\label{eqn:KeriG-cap-P}
\mathrm{Ker}(i_G) \simeq (\mathrm{Ker}(i_G) \cap P) \times \mathrm{Ker}(i_G)/(\mathrm{Ker}(i_G) \cap P).
\ee
\end{lemma}

\begin{proof}
By \S\S\ref{pp:Gm-actions}--\ref{pp:pseudo-parabolic}, the $k$-subgroup $i_G(P) \le i_G(G)$ is pseudo-parabolic. The $k$-group $\mathrm{Ker}(i_G)$ (or even $G$) acts on $G/P$, and if an $x \in G(S)$ lifts $y$, then the $y$-fiber of the map $G/P \rightarrow i_G(G)/i_G(P)$ is $\mathrm{Ker}(i_G)$-equivariantly isomorphic to
\[
\mathrm{Ker}(i_G)_S/(\mathrm{Ker}(i_G)_S \cap xP_Sx^{-1}) \cong x(\mathrm{Ker}(i_G)/(\mathrm{Ker}(i_G) \cap P))_Sx^{-1},
\]
which is a base change of $\mathrm{Ker}(i_G)/(\mathrm{Ker}(i_G) \cap P)$. Since $G$ is of minimal type, $\mathrm{Ker}(i_G)$ is commutative (see \S\ref{pp:pseudo-minimal}), so the conjugation action of $G$ on $\mathrm{Ker}(i_G)$ factors through $i_G(G)$. By fppf descent, we conclude that if $y$ merely lifts to an $S$-point of $i_G(G)$, then the aforementioned $y$-fiber is a $\mathrm{Ker}(i_G)/(\mathrm{Ker}(i_G)\cap P)$-torsor over $S$, and that in general the $y$-fiber is a torsor under an $S$-form of $\mathrm{Ker}(i_G)/(\mathrm{Ker}(i_G)\cap P)$, as claimed.

Let $\lambda \colon \mathbb{G}_{m,\, k} \rightarrow G$ be a $k$-homomorphism with $P = P_G(\lambda)$ and let $T \subset G$ be a maximal $k$-torus through which $\lambda$ factors. As we reviewed in \S\ref{pp:pseudo-minimal}, the base change $\mathrm{Ker}(i_G)_{k^s}$ is a direct product of its intersections $\mathrm{Ker}(i_G)_{k^s} \cap U_a$ with the root groups $U_a$ for the multipliable roots $a$ of $T_{k^s}$. In particular, \cite{CGP15}*{Proposition 2.1.8 (2)-(3) and Corollary 3.3.12} ensure that $(\mathrm{Ker}(i_G) \cap P)_{k^s}$ (resp., $(\mathrm{Ker}(i_G)/(\mathrm{Ker}(i_G) \cap P))_{k^s}$) is isomorphic to the direct product of the intersections $\mathrm{Ker}(i_G)_{k^s} \cap U_a$ for those multipliable roots $a$ of $T_{k^s}$ for which $\langle a, \lambda\rangle \ge 0$ (resp., $\langle a, \lambda\rangle < 0$). In particular, $(\mathrm{Ker}(i_G)/(\mathrm{Ker}(i_G) \cap P))_{k^s}$ is $\mathrm{Gal}(k^s/k)$-equivariantly identified with a direct factor of $\mathrm{Ker}(i_G)$ complementary to $(\mathrm{Ker}(i_G) \cap P)_{k^s}$, and Galois descent gives the claimed decomposition
\[
\mathrm{Ker}(i_G) \simeq (\mathrm{Ker}(i_G) \cap P) \times \mathrm{Ker}(i_G)/(\mathrm{Ker}(i_G) \cap P). \qedhere
\]
\end{proof}

For ultraminimal pseudo-reductive $G$, we control its pseudo-parabolic $k$-subgroups as follows. 

\begin{lemma} \label{lem:ultra-to-res}
Let $G$ be an ultraminimal pseudo-reductive $k$-group, let $i_G \colon G \hookrightarrow \mathrm{Res}_{k^\prime/k}(\overline{G})$ be its comparison map as in \uS\uref{pp:pred-iG}, let $P \le G$ be a pseudo-parabolic $k$-subgroup, and let $\overline{P} \le \overline{G}$ be the image of $P_{k^\prime}$. Then $i_G$ induces a closed immersion
\[
G/P \hookrightarrow \mathrm{Res}_{k^\prime/k}(\overline{G}/\overline{P}).
\]
\end{lemma}

\begin{proof}
Base change to $k^s$ allows us to assume that $k$ is separably closed (see \S\ref{pp:pred-iG} and \S\ref{pp:pseudo-parabolic}). Let $\lambda \colon \mathbb{G}_{m,\, k} \to G$ be a $k$-homomorphism with $P = P_G(\lambda)$ and let $\overline{\lambda} \colon \mathbb{G}_{m,\, k^\prime} \to \overline{G}$ be the composition of $\lambda_{k^\prime}$ and $G_{k^\prime} \twoheadrightarrow \overline{G}$. Then $\overline{\lambda}$ is the $k^\prime$-homomorphism corresponding to $\lambda \circ i_G$ and $\overline{P} = P_{\overline{G}}(\overline{\lambda})$ (see \S\ref{pp:Gm-actions}). Since $i_G$ is injective, we shorten $\lambda \circ i_G$ to $\lambda$ and note that, by \cite{CGP15}*{Propositions 2.1.13},
\[
\mathrm{Res}_{k^\prime/k}(\overline{P}) = P_{\mathrm{Res}_{k^\prime/k}(\overline{G})}(\lambda) \qquad \text{and} \qquad \mathrm{Res}_{k^\prime/k}(U_{\overline{G}}(-\overline{\lambda})) = U_{\mathrm{Res}_{k^\prime/k}(\overline{G})}(-\lambda).
\]
Therefore, by \S\ref{pp:Gm-actions}, the following square is Cartesian
$$
\xymatrix{
U_G(-\lambda) \times P_G(\lambda) \ar@{^(->}[r] \ar@{^(->}[d]  & \mathrm{Res}_{k^\prime/k}(U_{\overline{G}}(-\overline{\lambda})) \times \mathrm{Res}_{k^\prime/k}(\overline{P}) \ar@{^(->}[d]  \\
G \ar@{^(->}[r]^-{i_G}  & \mathrm{Res}_{k^\prime/k}(\overline{G}) 
}
$$
and its vertical (resp.,~horizontal) maps are open (resp.,~closed) immersions. In particular, since \( \mathrm{Res}_{k^\prime/k}(\overline{G})/\mathrm{Res}_{k^\prime/k}(\overline{P}) \cong \mathrm{Res}_{k^\prime/k}(\overline{G}/\overline{P}) \) (see \S\ref{pp:pseudo-parabolic}), the vertical (resp.,~the top horizontal) maps of 
$$
\xymatrix{
U_G(-\lambda) \ar@{^(->}[r] \ar@{^(->}[d] & \mathrm{Res}_{k^\prime/k}(U_{\overline{G}}(-\overline{\lambda})) \ar@{^(->}[d] \\
G/P \ar@{^(->}[r] & \mathrm{Res}_{k^\prime/k}(\overline{G}/\overline{P}) 
}
$$
are open (resp.,~closed) immersions, and this square is also Cartesian. To conclude that the bottom horizontal map is also a closed immersion, by \( G(k) \)-translation, all that remains is to show that the \( G(k) \)-translates of \( \mathrm{Res}_{k^\prime/k}(U_{\overline{G}}(-\overline{\lambda})) \) cover \( \mathrm{Res}_{k^\prime/k}(\overline{G}/\overline{P}) \). For this, since \( k^\prime/k \) is purely inseparable, by \cite{CGP15}*{Corollary A.5.4 (2)}, it suffices to show that the \( G(k) \)-translates of \( U_{\overline{G}}(-\overline{\lambda}) \) cover \( \overline{G}/\overline{P} \), which follows from \cite{CGP15}*{top of p. 587}. 
\end{proof}

\csub[The control of unipotent groups]
\label{sec:wound-control}

We conclude discussing the diagram of \S\ref{pp:fundamental-filtration} by reviewing the structure of wound unipotent groups. 

\bpp[Structure of smooth, connected, unipotent $k$-groups]\label{pp:wound-structure}
In general, to handle the subquotients of the cc$k$p filtration, we need to understand connected, smooth, commutative, $p$-torsion, unipotent $k$-groups $G$. These can be made explicit: by a result of Tits \cite{BLR90}*{Section 10.2, Proposition 10} (or \cite{CGP15}*{Proposition B.1.13}), every such $G$ is a subgroup of $\mathbb{G}_a^{N + 1}$ given by the vanishing locus of some $p$-polynomial 
$$
\textstyle F = \sum_{i = 0}^{N}\sum_{j = 0}^{n_i} f_{ij}t_i^{p^j} \in k[t_0, \dotsc, t_N] \quad \text{with} \quad f_{in_i} \neq 0 \quad \text{when} \quad n_{i} \ge 0,
$$
and the smoothness condition amounts to there being a nonzero linear term. This gives a concrete way to understand cohomology: functorially in a $k$-algebra $R$, we have
\be \label{eqn:unipotent-coho}
H^0(R, G) \cong (R^{N +1})^{F = 0}, \quad H^1(R, G) \cong R/F(R^{N + 1}), \quad H^{\ge 2}(R, G) \cong 0.
\ee
When $G$ is also wound, these references (with \cite{CGP15}*{Lemma B.1.7 (1)$\Leftrightarrow$(2)}, if one prefers) give more: we may choose an $F$ as above for which the polynomial of principal parts
$$
\tst F_{\mathrm{prin}} := \sum_{i = 0}^N f_{i n_i} t_i^{p^{n_i}}
$$
has no nontrivial zeros in $k^{N + 1}$, this is critical for working with wound groups in practice. In addition, once a unipotent group $G$ is given by the vanishing locus of a $p$-polynomial $F$, it is sometimes useful to project onto a proper subset of the coordinates $t_i$: this realizes $G$ as an extension of a power of $\mathbb{G}_a$ by the (possibly nonsmooth) unipotent group given by the vanishing locus of the $p$-polynomial obtained from $F$ by setting these $t_i$ to $0$. This also shows that woundness tends to be destroyed by passing to quotients.
\epp

Wound groups contain no $\mathbb{G}_{a,\, k}$'s as $k$-subgroups, but they may contain nontrivial unirational $k$-subgroups or even be unirational themselves. The latter are somewhat manageable as follows.

\begin{lemma}[Rosengarten, a positive answer to {\cite{Ach19}*{Question 2.21}}]\label{lem:Ros}
Every unirational, wound, unipotent $k$-group $G$ that is \emph{minimal}, in the sense that each of its unirational $k$-subgroups is either trivial or $G$ itself, is commutative. The nontrivial, unirational, wound, unipotent $k$-groups are all of dimension $\ge p-1$, and if $k$ is separably closed of characteristic $p > 0$, then they are precisely the following $k$-groups with $a \in k \setminus k^p$ and $n \ge 1$:
\[
\tst G_{a,\, n} \colonequals \mathrm{Res}_{k(a^{1/p^{n-1}})/k}\left( \mathrm{Res}_{k(a^{1/p^{n}})/k(a^{1/p^{n-1}})}\left( \mathbb{G}_{m,\, k(a^{1/p^n})} \right)/\mathbb{G}_{m,\, k(a^{1/p^{n-1}})} \right).
\]
\end{lemma}

\begin{proof}
This lemma and its proof are due to Rosengarten. Since $G$ is unirational, it is smooth and connected (see \S\ref{pp:fundamental-filtration}~\ref{m:Guni}). By \S\ref{pp:DG} and the cc$k$p filtration, the derived subgroup $\mathscr{D}(G)$ is smaller than $G$ and unirational. Thus, if $G$ is minimal, then $\mathscr{D}(G) = 1$, so that $G$ is commutative.

The dimension assertion follows from the rest, so we now assume that $k$ is separably closed and that $G$ is a nontrivial, unirational, wound, unipotent $k$-group that is minimal (so also commutative). By {\cite{Ros21c}*{proof of Theorem 1.3 on pp.~442--443}} (or {\cite{Ach19}*{Remark 2.6(ii)}}), the unirationality implies that $G$ is a quotient of $\mathrm{Res}_{A/k}(\mathbb{G}_{m,\, A})$ for some finite $k$-algebra $A$. Since $\mathrm{Res}_{A/k}(\mathbb{G}_{m,\, A})$ is unirational and decomposes into a product according to the factors of $A$, the minimality of $G$ allows us to assume that $A$ is local.  By the proof of \Cref{lem:res-unipotent}, the surjection $\mathrm{Res}_{A/k}(\mathbb{G}_{m,\, A}) \twoheadrightarrow \mathrm{Res}_{A^{\mathrm{red}}/k}(\mathbb{G}_{m,\, A^{\mathrm{red}}})$ has a split unipotent kernel, so, since $G$ is wound, we may pass to $A^{\mathrm{red}}$ to assume that $A$ is a finite, purely inseparable field extension $k' / k$. Since $k'$ is a sum of monogenic subextensions $k \subset k_i \subset k'$, we check on Lie algebras that $\mathrm{Res}_{k'/k}(\mathbb{G}_{m,\, k'})$ is a quotient of $\prod_i \mathrm{Res}_{k_i/k}(\mathbb{G}_{m,\, k_i})$. Since $G$ is minimal and nontrivial, we may therefore replace $k'$ by some $k_i$ to reduce to the case when $k' = k(a^{1/p^n})$ for some $a \in k \setminus k^p$ and $n \ge 1$. By decreasing $n$ if needed and again using the minimality of $G$, we may assume that the map $\mathrm{Res}_{k(a^{1/p^{n-1}})/k}(\mathbb{G}_{m,\, k(a^{1/p^{n-1}})}) \to G$ vanishes, so that
\be\label{eqn:RR-surj}
\mathrm{Res}_{k(a^{1/p^n})/k}\left(\mathbb{G}_{m,\, k(a^{1/p^n})}\right)/\mathrm{Res}_{k(a^{1/p^{n-1}})/k}\left(\mathbb{G}_{m,\, k(a^{1/p^{n-1}})}\right) \twoheadrightarrow G.
\ee
The source of this surjection is $G_{a,\, n}$, see \S\ref{pp:conv}. Moreover, $G_{a,\, n}$ is wound and given by the vanishing of some $p$-polynomial $F$ of degree $p$: indeed, this holds for $G_{a,\, 1}$ (see {\cite{Oes84}*{chapitre VI, proposition~5.3}}) and is stable under restrictions of scalars along purely inseparable field extensions (decomposed into successive extensions of degree $p$). Thus, by {\cite{Ros25}*{proof of Proposition 7.7}}, every quotient of $G_{a,\, n}$ by a nontrivial $k$-subgroup is split unipotent. As $G$ is wound, \eqref{eqn:RR-surj} must then be an isomorphism.

The identification of the source of \eqref{eqn:RR-surj} with $G_{a,\, n}$ shows that this wound, unipotent $k$-group is unirational, so it remains to show that it is minimal; equivalently, by the argument above, that it admits no $G_{b,\, m}$ with $b \in k \setminus k^p$ and $m \ge 1$ as a proper $k$-subgroup. By dimension considerations, this is not possible if $n = 1$, so, arguing by induction for all $k$ at once, we assume that $n > 1$ is minimal for which some $G_{a,\, n}$ is not minimal, and hence properly contains some $G_{b,\, m}$. The unipotent $k$-group $G_{b,\, m}$ splits over $k(b^{1/p^m})$, so $G_{a,\, n}$ cannot remain wound over $k(b^{1/p^m})$. The extensions $k(a^{1/p^n})$ and $k(b^{1/p^m})$ are then not disjoint over $k$, that is, $k(a^{1/p}) = k(b^{1/p})$. On the other hand, by the definition of the Weil restriction, we have a nonzero $k(a^{1/p})$-homomorphism $(G_{b,\, m})_{k(a^{1/p})} \to G_{a^{1/p},\, n-1}$. \Cref{lem:res-unipotent} ensures that $(G_{b,\, m})_{k(b^{1/p})}$ is an extension over $k(b^{1/p})$ of $G_{b^{1/p},\, m-1}$ (interpreted to be zero if $m=1$) by a split unipotent group. Since $G_{a^{1/p},\, n-1}$ is wound, we obtain a nonzero $k(b^{1/p})$-homomorphism $G_{b^{1/p},\, m-1} \to G_{a^{1/p},\, n-1}$, which, by induction, must be surjective. Then, however, the dimension of $G_{b^{1/p},\, m-1}$ is at least that of $G_{a^{1/p},\, n-1}$, so the dimension of $G_{b,\, m}$ is at least that of $G_{a,\, n}$, and hence that $G_{b,\, m}$ cannot be a proper $k$-subgroup of $G_{a,\, n}$.
\end{proof}

\bd\label{def:strongly-wound}
A unipotent $k$-group $G$ is \emph{strongly wound}
\begin{itemize}
    \item if $G$ has no nontrivial unirational $k$-subgroups, that is, if $G^{\mathrm{uni}} = 1$; equivalently,
    \item if $G$ has no nontrivial, \emph{commutative}, unirational $k$-subgroups (see \Cref{lem:Ros}); equivalently,
    \item if $G$ receives no nontrivial $k$-homomorphism from a unirational $k$-group.
\end{itemize}
\ed

\beg
By \Cref{lem:Ros}, a wound, unipotent $k$-group of dimension $< p-1$ is strongly~wound.
\eeg

By \S\ref{pp:fundamental-filtration}~\ref{m:Guni}, this condition is insensitive to base change to a separable field extension. It is also inherited by the subquotients of the cc$k$p filtration as follows.

\begin{lemma}\label{lem:sw-cckp}
Let $G$ be a connected, smooth, strongly wound, unipotent $k$-group. Every subquotient of the cc$k$p filtration of $G$ is strongly wound \up{and connected, smooth, unipotent}.
\end{lemma}

\begin{proof}
Let $G_1 \lhd G$ be the largest connected, smooth, central, $p$-torsion $k$-subgroup of $G$. If $(G/G_1)^{\mathrm{uni}} = 1$, then we may replace $G$ by $G/G_1$ and iterate, so, for the sake of contradiction, we choose a nontrivial, unirational $k$-subgroup $H \le G/G_1$ and let $\widetilde{H} \le G$ be its preimage. We have a central extension of $k$-groups
\[
1 \to G_1 \to \widetilde{H} \to H \to 1
\]
and, to obtain a contradiction with the definition of $G_1$, we will show that $\widetilde{H}$ is central in $G$ and $p$-torsion. For this, we imitate the proof of \cite{Con15b}*{Proposition 3.2}.

We lose no generality by base changing to $k^s$ to assume that $k$ is separably closed (see \S\ref{pp:fundamental-filtration}~\ref{m:cckp}). To see that $\widetilde{H}$ is central in $G$, for a $g \in G(k)$ we consider the $k$-morphism $\widetilde{H} \to G$ given by $h \mapsto ghg^{-1}h^{-1}$. Since $G_1$ is central in $G$, this morphism factors through a $k$-morphism $H \to G$. However, $G$ is strongly wound, so every $k$-morphism from a unirational $k$-variety (such as $H$) to $G$ that sends some $k$-point to the identity is constant. Consequently, the commutators $ghg^{-1}h^{-1}$ all vanish and, since $g$ is arbitrary and $k$ is separably closed, we get that $\widetilde{H}$ is central in $G$. Similarly, since $G_1$ is $p$-torsion and central in $G$, the $p$-power morphism $\widetilde{H} \to G$ factors through a $k$-morphism $H \to G$ that must a posteriori be constant, to the effect that $\widetilde{H}$ is $p$-torsion.
\end{proof}

For treating nonsmooth, commutative, unipotent groups, the following lemma is useful. 

\begin{lemma}[\cite{BLR90}*{Section 10.2, Lemma 13}]\label{lem:embed-wound}
A commutative, connected, \up{resp., wound}~unipotent $k$-group is a subgroup of a commutative, connected, smooth, \up{resp., wound} unipotent $k$-group. \qed
\end{lemma}

\section{Pseudo-properness and extension results for sections}

In \S\ref{sec:pseudo-proper}, we introduce the notion of pseudo-properness, which strengthens that of pseudo-completeness due to Tits, and we classify pseudo-proper $k$-groups. In \S\ref{sec:extension-generic}, we establish general extension results for sections of torsors under pseudo-proper groups, and thus reduce our main result \Cref{thm:main-GS} to connected, smooth, affine groups. We use these extension results in \S\ref{sec:GP-pseudo-proper} to show that $G/P$ is pseudo-proper for every pseudo-parabolic $k$-subgroup $P \le G$ of a smooth, affine $k$-group $G$.

\csub[Pseudo-properness and pseudo-completeness] \label{sec:pseudo-proper}

A $k$-variety $X$ is proper if and only if $X(R) \xrightarrow{\sim} X(K)$ for every discrete valuation ring $R$ that is a $k$-algebra with fraction field $K$. By only requiring this valuative criterion to hold for certain subclasses of $R$, we obtain the following variants of the notion of properness over imperfect $k$.

\bd\label{def:pseudo-proper}
A $k$-algebraic space $X$ is
\benumr
    \item \label{m:PP-i} \emph{pseudo-proper} if it is separated, of finite type, and for every discrete valuation ring $R$ that is a geometrically regular $k$-algebra with fraction field $K$, we have
    \begin{equation}\label{eqn:pseudo-proper}
        X(R) \xrightarrow{\sim} X(K);
    \end{equation}
    
    \item \label{m:PP-ii} \emph{pseudo-complete} if it is separated, of finite type, and~\eqref{eqn:pseudo-proper} holds for every discrete valuation ring $R$ that is a $k$-algebra whose residue field is separable over~$k$.
\end{enumerate}
\ed

The notion of pseudo-completeness has already been introduced by Tits in \cite{Tit13}*{cours 1992--1993, section 2.3}, see also \cite{CGP15}*{Definition C.1.1}.

By \cite{Mat89}*{Section 28, Lemma 1 on p. 216, Theorem 28.7}, a discrete valuation ring $R$ that is a $k$-algebra whose residue field is separable over $k$ is geometrically regular over $k$, so
\be\ba\label{eqn:pseudo-fpc}
\xymatrix{
& \text{proper} \ar@{=>}[rd] & & \\
\text{finite} \ar@{=>}[rd] \ar@{=>}[ru]  & & \text{pseudo-proper} \ar@{=>}[r] & \text{pseudo-complete}\\
& \text{pseudo-finite} \ar@{=>}[ru] & \\
}
\ea\ee
where the implication $\text{pseudo-finite} \Longrightarrow \text{pseudo-proper}$ follows by noting that pseudo-properness may be tested on $X^{\mathrm{gred}}$, see \Cref{rem:geometric-regularity} below (and compare also with the more general \Cref{thm:main-GS} below). If $k$ is perfect, then these implications, except for the two that point top right, are all equalities, see \cite{SP}*{Lemma~\href{https://stacks.math.columbia.edu/tag/0H1X}{0H1X}} and \Cref{rem:pseudo-finite-perfect}. Over an imperfect $k$, these implications are all strict thanks to \Cref{eg:pseudo-etale}, the following example, and \Cref{eg:pseudo-complete} below.

\beg[pseudo-proper $\not\Rightarrow$ proper]\label{eg:pseudo-proper}
For a finite extension $k^\prime/k$ and a proper (or merely pseudo-proper) $k^\prime$-algebraic space $X^\prime$, the restriction of scalars $X := \mathrm{Res}_{k^\prime/k}X^\prime$ (see~\S\ref{pp:conv}) is pseudo-proper (\emph{loc.~cit.}). In fact, pseudo-properness ascends along proper morphisms, so an algebraic space that is proper over this~$X$, for instance, a closed subspace of~$X$, is also pseudo-proper. For nonseparable $k'/k$, such~$X$ are often not proper; see~\cite{CGP15}*{Example A.5.6}.

This example shows why it is natural to allow algebraic spaces in~\Cref{def:pseudo-proper}: namely,~$X$ need not be a scheme if~$X'$ is not quasi-projective. For instance,~$X$ is not a scheme if $k^\prime/k$ is separable quadratic and $X^\prime$ has two $k^\prime$-points that are not contained in any affine open subscheme of $X^\prime_{\overline{k}^\prime}$, and Hironaka gave examples of proper such~$X^\prime$ with $k^\prime = \mathbb{C}$ that are not projective.
\eeg

\begin{remark}\label{rem:geometric-regularity}
The map~\eqref{eqn:pseudo-proper} is injective because~$X$ is separated (see~\cite{EGAI}*{corollaire~9.5.6}), so we are only imposing its surjectivity. In particular, since the geometric regularity of~$R$ forces~$K$ to be separable over~$k$, by~\eqref{eqn:reduce-to-gred}, a separated $k$-algebraic space~$X$ of finite type is pseudo-proper (resp., pseudo-complete) if and only if so is $X^{\mathrm{gred}}$. Moreover, since~$R$ is $1$-dimensional, its geometric regularity in~\ref{m:PP-i} is equivalent to geometric normality: concretely, it means that $R \otimes_k k'$ is required to be a semilocal Dedekind ring for every finite field extension $k'/k$. By~\cite{EGAIV2}*{proposition~5.13.7} and the Popescu theorem (see~\S\ref{pp:conv}), in~\ref{m:PP-i} the strict Henselization of~$R$ is still geometrically regular over~$k$, so, by~\cite{SP}*{Lemma~\href{https://stacks.math.columbia.edu/tag/0ARH}{0ARH}}, we may restrict to strictly Henselian~$R$ in~\ref{m:PP-i} (and also in~\ref{m:PP-ii}).
\end{remark}

\begin{lemma}\label{lem:pseudo-base-change}
Let $X$ be a $k$-algebraic space and let $k^\prime/k$ be a field extension.
\benum
    \item\label{m:PBC-a} If $X_{k^\prime}$ is pseudo-proper \up{resp., pseudo-complete} over $k^\prime$, then so is $X$ over $k$.
    \item\label{m:PBC-b} If $k^\prime/k$ is separable and $X$ is pseudo-proper \up{resp., pseudo-complete} over $k$, then so is $X_{k^\prime}$ over $k^\prime$.
\end{enumerate}
\end{lemma}

\begin{proof}
By fpqc descent \cite{SP}*{Lemmas~\href{https://stacks.math.columbia.edu/tag/0421}{0421} and \href{https://stacks.math.columbia.edu/tag/041U}{041U}}, being separated or of finite type may be checked after base change to~$k^\prime$. Thus, for~\ref{m:PBC-b}, it suffices to note that since~$k^\prime/k$ is separable, a ring (resp., a field) that is geometrically regular (resp., separable) as a~$k^\prime$-algebra is also geometrically regular (resp., separable) as a~$k$-algebra.

In the remaining~\ref{m:PBC-a}, we begin with the case when~$k^\prime/k$ is separable. This assumption and \cite{CGP15}*{argument on top of page~583} ensure that for every discrete valuation ring~$R$ that is a geometrically regular~$k$-algebra, there is a local injection~$R \to R^\prime$ over~$k \to k^\prime$ such that~$R^\prime$ is a localization of~$R \otimes_k k^\prime$ at some prime ideal and also a discrete valuation ring that is geometrically regular as a~$k^\prime$-algebra. It then remains to apply \cite{SP}*{Lemma~\href{https://stacks.math.columbia.edu/tag/0ARH}{0ARH}} to get the separable case of~\ref{m:PBC-a}. An analogous argument works for an arbitrary~$k^\prime/k$ granted that we focus on the pseudo-completeness variant (although the latter also results by adapting the argument that follows). In general, due to the settled~\ref{m:PBC-b} and the separable case of~\ref{m:PBC-a}, it remains to settle the pseudo-proper variant of~\ref{m:PBC-a} in the case when~$k^\prime/k$ is purely inseparable. In this case, we fix a discrete valuation ring~$R$ that is a geometrically regular~$k$-algebra with fraction field~$K$ and a~$K$-point of~$X$ that we wish to extend to an~$R$-point. By the pseudo-properness of~$X_{k^\prime}$ over~$k^\prime$ and a limit argument, there is a finite, purely inseparable field extension~$\ell/k$ such that the induced~$(K \otimes_k \ell)$-point of~$X_\ell$ extends to an~$(R \otimes_k \ell)$-point, equivalently, such that the induced~$K$-point of~$\mathrm{Res}_{\ell/k}(X_\ell)$ extends to an~$R$-point. To get the desired~$R$-point of~$X$, it remains to note that, since~$X$ is separated, the map~$X \hookrightarrow \mathrm{Res}_{\ell/k}(X_\ell)$ of the adjunction is a closed immersion, see \cite{BLR90}*{bottom of p.~197} and \cite{SP}*{Lemma~\href{https://stacks.math.columbia.edu/tag/03KP}{03KP}}.
\end{proof}

\begin{proposition}\label{prop:PP-groups}
Let $G$ be a $k$-group scheme of finite type.
\benum
    \item\label{m:PPG-0} $G$ is pseudo-finite if and only if $(G^{\mathrm{gred}})^0 = 1$.
    \item\label{m:PPG-i} $G$ is pseudo-complete if and only if $G^{\mathrm{sm},\, \mathrm{lin}}$ is wound unipotent.
    \item\label{m:PPG-ii} $G$ is pseudo-proper if and only if $G^{\mathrm{sm},\, \mathrm{lin}}$ is strongly wound unipotent, equivalently, if and only if $G$ has no nontrivial unirational $k$-subgroup \up{that is\ucolon $G^{\mathrm{uni}} = 1$}.
\end{enumerate}

In particular, a pseudo-abelian variety is pseudo-proper, and it is proper if and only if it is an abelian variety. In each of these statements, the corresponding pseudo-finiteness, pseudo-completeness, pseudo-properness, or properness conclusion also holds for every \'etale locally trivial $G$-torsor over $k$.
\end{proposition}

\begin{proof}
By \Cref{lem:pseudo-base-change}, \Cref{rem:Xgred-etale}, \S\ref{pp:fundamental-filtration}, and \Cref{rem:geometric-regularity}, the conditions in question are insensitive to base change to $k^s$ and only depend on $G^{\mathrm{gred}}$, so we may assume that $k$ is separably closed and $G$ is smooth. Each $G$-torsor is then trivial and each connected component of $G$ is $k$-isomorphic to $G^0$, so it is enough to treat $G$ under the assumption that it is connected.

Since $k$ is separably closed, $G$ is pseudo-finite if and only if $G(k)$ is finite. This happens if and only if the closure of $G(k)$ in $G$ is $0$-dimensional, that is, if and only if $(G^{\mathrm{gred}})^0 = 1$.

For the rest, we first argue the claim about pseudo-abelian varieties. By \S\ref{pp:pseudo-abelian}, it suffices to show that each pseudo-abelian variety $G$ is pseudo-proper. There are many ways to see this: for instance, by the Popescu theorem (see \S\ref{pp:conv}) and the Bosch--L\"utkebohmert--Raynaud criterion recalled in \S\ref{pp:fundamental-filtration}~\ref{m:Guni}, it suffices to check that there is no nontrivial, connected, unirational (and hence smooth) $k$-subgroup $U \le G$. The abelian variety $G^{\mathrm{av}}$ has no such $k$-subgroup, so $U$ lies in $G^{\mathrm{lin}}$, and so is affine. The definition of a pseudo-abelian variety then implies that $U$ is trivial. For another argument for the pseudo-properness of $G$, see \Cref{thm:pseudo-extend}~\ref{m:PE-ii} below and its proof.

For the remaining~\ref{m:PPG-i} and~\ref{m:PPG-ii}, since pseudo-completeness or pseudo-properness may be tested after restricting to strictly Henselian discrete valuation rings $R$ (subject to the respective condition on $R$, see \Cref{rem:geometric-regularity}), and thanks to the already settled case of pseudo-abelian varieties, we may replace $G$ by $G^{\mathrm{sm},\, \mathrm{lin}}$ to assume that our connected, smooth $G$ is affine. Since neither $\mathbb{G}_a$ nor $\mathbb{G}_m$ are pseudo-complete, they cannot occur as $k$-subgroups of any pseudo-complete $G$. In particular, every pseudo-complete, connected, smooth, affine $k$-group $G$ satisfies $G^{\mathrm{tor}} = 1$, and so is even wound unipotent by \S\ref{pp:fundamental-filtration}~\ref{m:Guni}. Conversely, the pseudo-completeness of wound unipotent groups is a result of Tits, see \cite{CGP15}*{Lemma C.1.8}, that follows by combining the $\mathrm{cc}k\mathrm{p}$ filtration with \eqref{eqn:unipotent-coho}. This settles~\ref{m:PPG-i}.

Pseudo-properness of strongly wound unipotent groups in~\ref{m:PPG-ii} follows from the $\mathrm{cc}k\mathrm{p}$ filtration, \Cref{lem:sw-cckp}, and the same Bosch--L\"utkebohmert--Raynaud criterion from \S\ref{pp:fundamental-filtration}~\ref{m:Guni}. Conversely, if a smooth, wound unipotent $k$-group $G$ is pseudo-proper, then the same criterion implies that $G$ has no nontrivial, commutative, unirational $k$-subgroups, so that $G$ is strongly wound by \Cref{def:strongly-wound}. Thanks to~\ref{m:PPG-i}, this gives the first “if and only if” in~\ref{m:PPG-ii}. For the second “if and only if” in~\ref{m:PPG-ii}, it then remains to review \S\ref{pp:fundamental-filtration}~\ref{m:Gtor}.
\end{proof}

\beg[pseudo-complete $\not\Rightarrow$ pseudo-proper]\label{eg:pseudo-complete}
For every finite, purely inseparable field extension $k'/k$, by \cite{Oes84}*{chapitre VI, section 5.1, lemme} (alternatively, see \Cref{lem:Ros} above),  
\[
U := \mathrm{Res}_{k'/k}(\mathbb{G}_m)/\mathbb{G}_m
\]
is a wound unipotent, unirational (in particular, not strongly wound) $k$-group, so, by \Cref{prop:PP-groups}~\ref{m:PPG-ii}, it is pseudo-complete but not pseudo-proper unless $k' = k$. In general, to appreciate the difference between pseudo-properness and pseudo-completeness of a separated $k$-algebraic space $X$ of finite type, consider a smooth affine $k$-curve $C$ and a closed point $c \in C$. If $X$ is pseudo-proper, then every $k$-morphism $C - c \to X$ extends uniquely to a $k$-morphism $C \to X$; in contrast, if $X$ is only pseudo-complete over $k$, then such an extension is guaranteed only when $k(c)/k$ is separable.
\eeg

\begin{remark}\label{rem:res}
Pseudo-properness is intriguing from the point of view of resolving singularities. Namely, for a smooth, integral, pseudo-proper $k$-scheme $X$ and any regular, proper, integral compactification $X \subset \overline{X}$ with $\overline{X} \setminus X$ a divisor, by \eqref{eqn:pseudo-proper}, we have $\overline{X}^{\mathrm{sm}} = X$, that is, the $k$-smooth locus of $\overline{X}$ is not larger than $X$. In other words, $\overline{X}$ is regular but not smooth anywhere along the boundary $\overline{X} \setminus X$.  Such an $\overline{X}$ ought to exist by the resolution of singularities conjecture but remains elusive even in the cases of \Cref{prop:PP-groups}, for instance, for pseudo-abelian varieties or for connected, smooth, strongly wound unipotent groups. The pseudo-properness of these two types of groups and \cite{BLR90}*{Section 10.3, Theorem 1} show that they never possess a smooth compactification over $k$ unless they are abelian varieties to begin with.
\end{remark}

\csub[Extending generic sections and Grothendieck--Serre for pseudo-abelian varieties]
\label{sec:extension-generic}
An integral scheme $S$ with fraction field $K$ and any finite $S$-scheme $X$ satisfy
\be\label{eqn:finite-extend}
X(S) \xrightarrow{\sim} X(K),
\ee
see, for instance, \cite{modular-description}*{Lemma 3.1.9}. This is useful in many contexts, for instance, in moduli theory, where $X$ may parametrize automorphisms of some object and we may wish to extend automorphisms over $K$ to those of an entire family over $S$. We extend this to pseudo-finite $X$ (see \Cref{def:pseudo-finite}) as follows.

\begin{proposition}\label{prop:pseudo-finite-extend}
Let $S$ be a normal integral $k$-scheme whose function field $K$ is separable over $k$ \up{an automatic assumption if $S$ is geometrically reduced}. For a pseudo-\up{locally quasi-finite} $k$-algebraic space $X$ \up{resp.,~a $k$-group scheme $G$ locally of finite type}, we have
\[
X(S) \xrightarrow{\sim} X(K) \quad \text{{\upshape(}resp.,}\quad (G/(G^{\mathrm{gred}})^0)(S) \xrightarrow{\sim} (G/(G^{\mathrm{gred}})^0)(K)).
\]
\end{proposition}

\begin{proof}
By \eqref{eqn:reduce-to-gred} and \Cref{rem:Xgred-etale}, we may replace $X$ by $X^{\mathrm{gred}}$ to reduce to the case when $X$ is a separated, \'etale $k$-scheme. The separatedness ensures that the map in question is injective, see \cite{EGAI}*{corollaire 9.5.6}. Thus, for its surjectivity, we may replace $X$ by an affine open containing the $K$-point of interest to reduce to when $X$ is finite \'etale. The finite case, however, is a special instance of \eqref{eqn:finite-extend}. The claim about $G$ follows from that about $X$ and \Cref{eg:pseudo-finite-G}.
\end{proof}

An integral, regular scheme $S$ with fraction field $K$ and a torsor $E$ under an abelian $S$-scheme $A$ satisfy
\be\label{eqn:abelian-scheme-extend}
E(S) \xrightarrow{\sim} E(K), \quad\quad \text{so that, consequently,} \quad\quad H^1(S, A) \hookrightarrow H^1(K, A),
\ee
see, for instance, \cite{BLR90}*{Section 8.4, Corollary 6} or \cite{torsors-regular}*{Remark 3.1.7} (essentially, it suffices to view $A$ as the dual of its dual $A^\vee$ and use the regularity to extend line bundles), alternatively \cite{Bha12a}*{Proposition 4.2 and Remark 4.5} (this last argument is based on finding projective lines in positive-dimensional fibers of modifications of regular schemes, see also \cite{KS15}*{Theorem~B.0.1}). This Grothendieck–Serre type conclusion for abelian schemes is useful in many contexts.

We extend this to pseudo-abelian varieties, in particular, we establish the pseudo-abelian variety case of our main Grothendieck–Serre \Cref{thm:main-GS} in part \ref{m:PE-b}~\ref{m:PE-ii} of the following \Cref{thm:pseudo-extend}. Its part \ref{m:PE-a}~\ref{m:PE-ii} also extends \Cref{prop:pseudo-finite-extend} to torsors over $S$ under pseudo-finite $k$-groups. The key statement is \eqref{eqn:Gabber-lemma} in \ref{m:PE-a} and was announced by Gabber in \cite{Gab12}*{bottom of p.~2371} under an additional unipotence assumption.

\begin{theorem}\label{thm:pseudo-extend}
Let $G$ be a $k$-group scheme locally of finite type, let $S$ be an integral $k$-scheme, and let $E$ be a $G$-torsor over $S$.
\benum
    \item \label{m:PE-0} We have
    \begin{equation}\label{eqn:to-G0}
\qq    (E/G^0)(S) \xrightarrow{\sim} (E/G^0)(K).
    \end{equation}
    
    \item \label{m:PE-a} If $S$ is geometrically normal over $k$ and every $\overline{k}$-torus of $G_{\overline{k}}$ lies in $(G^{\mathrm{gred}})_{\overline{k}}$, then
    \begin{equation}\label{eqn:Gabber-lemma}
\qq    (E/(G^{\mathrm{gred}})^0)(S) \xrightarrow{\sim} (E/(G^{\mathrm{gred}})^0)(K).
    \end{equation}
    
    \item \label{m:PE-b} If $S$ is geometrically normal over $k$ and regular and every $\overline{k}$-torus of $G_{\overline{k}}$ lies in $(G^{\mathrm{gred}})_{\overline{k}}$, then
    \begin{equation}\label{eqn:E-Gsmlin}
\qq    (E/G^{\mathrm{sm},\, \mathrm{lin}})(S) \xrightarrow{\sim} (E/G^{\mathrm{sm},\, \mathrm{lin}})(K).
    \end{equation}
\end{enumerate}

In particular, in the setting of \ref{m:PE-0} \up{resp.,~\ref{m:PE-a}, resp.,~\ref{m:PE-b}},
\benumr
    \item \label{m:PE-i} each generically trivial $G$-torsor over $S$ reduces to a generically trivial $G^0$-torsor \up{resp.,~$(G^{\mathrm{gred}})^0$-torsor, resp.,~$G^{\mathrm{sm},\, \mathrm{lin}}$-torsor} over $S$\uscolon
    \item \label{m:PE-ii} if $G^0$ is trivial \up{resp.,~is pseudo-finite, resp.,~is pseudo-proper}, then
    \begin{equation}\label{eqn:pseudo-AV-extend}
    E(S) \xrightarrow{\sim} E(K) \quad \text{and} \quad \mathrm{Ker}(H^1(S, G) \rightarrow H^1(K, G)) = \{*\}.
    \end{equation}
\end{enumerate}
\end{theorem}

\begin{proof}
Preimages in $E$ of scheme-valued points of $E/G^0$ correspond to reductions of $E$ to a $G^0$-torsor, and similarly for $(G^{\mathrm{gred}})^0$ and $G^{\mathrm{sm},\, \mathrm{lin}}$, so \eqref{eqn:to-G0}, \eqref{eqn:Gabber-lemma}, and \eqref{eqn:E-Gsmlin} imply \ref{m:PE-i}. They also imply~\ref{m:PE-ii} in the cases \ref{m:PE-0}, \ref{m:PE-a}, and reduce it to connected, smooth, affine $G$ in the case \ref{m:PE-b} because if $G^0$ is trivial (resp.,~pseudo-finite; resp.,~pseudo-proper), then $G^0$ is trivial (resp.,~$(G^{\mathrm{gred}})^0$ is trivial; resp.,~$G^{\mathrm{sm},\, \mathrm{lin}}$ is pseudo-proper). To conclude \ref{m:PE-ii} in the case \ref{m:PE-b}, we may work \'etale locally on $S$, so, by the smoothness of $G$, may assume that $E$ is trivial. The pseudo-properness of $G$ then extends every $K$-point of $G$ uniquely to a $U$-point for some open $U \subset S$ with complement of codimension $\geq 2$ (depending on the $K$-point). The affineness of $G$, the codimension condition, and, for instance, \cite{flat-purity}*{Lemma 7.2.7(b)} (recalled in \eqref{eqn:S-depth-2} below), then extend a $U$-point uniquely to an $S$-point.

\benum
    \item By working \'etale locally on $S$ to trivialize the $G/G^0$-torsor $E/G^0$ (see \S\ref{pp:fundamental-filtration}~\ref{m:G0}), we find from \eqref{eqn:finite-extend} that
    \[
\qq    (E/G^0)(S) \xrightarrow{\sim} (E/G^0)(K).
    \]
    
    \item Since $(G^{\mathrm{gred}})^0$ is closed in $G$, the quotient $E/(G^{\mathrm{gred}})^0$ is a separated $S$-algebraic space (see~\S\ref{pp:conv}), so \cite{EGAI}*{corollaire 9.5.6} and its proof ensure that \eqref{eqn:Gabber-lemma} is injective. For its surjectivity, by working locally, we may assume that $S$ is affine. In addition, since the preimage in $E$ of any $S$-point of $E/G^0$ is a $G^0$-torsor, by \ref{m:PE-0}, we may assume that $G$ is connected, so of finite type. If then $G$ is already smooth, so that $G = G^{\mathrm{gred}}$, then $E/(G^{\mathrm{gred}})^0$ is $S$-finite and \eqref{eqn:Gabber-lemma} is a special case of \eqref{eqn:finite-extend}. Thus, it suffices to reduce to smooth $G$, and for this it is enough to argue that we lose no generality by replacing $G$ by successively deeper subgroups $G^{(i)}$ of the filtration \eqref{eqn:gred-filtration}. Moreover, since $(G^{(i)})^{\mathrm{gred}} = G^{\mathrm{gred}}$, each $G^{(i)}$ inherits our assumption on the $\overline{k}$-tori. Thus, for a finite, purely inseparable field extension $k'/k$ and $k'$-group quotient $Q \colonequals G_{k'}/(G_{k'})^{\mathrm{red}}$, it suffices to show how to replace $G$ by the stabilizer $G_q$ of the unique $q \in (\mathrm{Res}_{k'/k}(Q))(k)$. For this, since $G^{\mathrm{gred}} \leq G_q$ (see \S\ref{pp:gred-filtration}), our assumption on the $\overline{k}$-tori and \Cref{lem:GGMB} ensure the closedness of the orbit $G/G_q \cong G \cdot q \subset \mathrm{Res}_{k'/k}(Q)$. Consider the affine $S$-scheme defined by the contracted product $E \times^G (\mathrm{Res}_{k'/k}(Q))_S$. The adjunction counit $(\mathrm{Res}_{k'/k}(Q))_{k'} \rightarrow Q$ is $G_{k'}$-equivariant (see \S\ref{pp:gred-filtration}), so it gives an $S$-morphism
    \begin{equation}\label{eqn:EQ}
\qq    E \times^G (\mathrm{Res}_{k'/k}(Q))_S \xrightarrow{\sim} \mathrm{Res}_{S_{k'}/S}(E_{S_{k'}} \times^{G_{k'}} Q_{S_{k'}}),
    \end{equation}
    which we check to be an isomorphism by working fppf locally on $S$ to trivialize $E$. Since $Q$ is a finite $k'$-scheme, its $S_{k'}$-form $E_{S_{k'}} \times^{G_{k'}} Q_{S_{k'}}$ is a finite $S_{k'}$-scheme, and $S_{k'}$ is normal and integral by the geometric normality of $S$. Thus, \eqref{eqn:EQ} and \eqref{eqn:finite-extend} give
    \[
\qq    (E \times^G (\mathrm{Res}_{k'/k}(Q))_S)(S) \xrightarrow{\sim} (E \times^G (\mathrm{Res}_{k'/k}(Q))_S)(K).
    \]
    Since $G \cdot q \subset \mathrm{Res}_{k'/k}(Q)$ is closed, so is $E \times^G (G \cdot q)_S \subset E \times^G (\mathrm{Res}_{k'/k}(Q))_S$. Moreover, fppf locally on $S$, we have $E/G_q \xrightarrow{\sim} E \times^G (G \cdot q)$. Thus, also
    \[
\qq    (E/G_q)(S) \xrightarrow{\sim} (E/G_q)(K).
    \]
    Now, a $K$-point of $E/(G^{\mathrm{gred}})^0$ lifts to an $S$-point after replacing $E$ by the corresponding $G_q$-subtorsor, reducing the proof of \eqref{eqn:Gabber-lemma} to the case of $G_q$, as desired.
    
    \item Similarly to the proof of \ref{m:PE-a}, \eqref{eqn:Gabber-lemma} reduces \eqref{eqn:E-Gsmlin} to the case when $G$ is connected and smooth. Then $G/G^{\mathrm{lin}}$ is an abelian variety (see \S\ref{pp:fundamental-filtration}~\ref{m:Glin}), so \eqref{eqn:abelian-scheme-extend} applies, and the argument of \ref{m:PE-0} reduces to when $G$ is connected and affine (but possibly not smooth). To conclude \eqref{eqn:E-Gsmlin}, we review \S\ref{pp:fundamental-filtration}~\ref{m:Gsmlin} and apply \eqref{eqn:Gabber-lemma} once again. \qedhere
\end{enumerate}
\end{proof}

\begin{remark} \label{rem:PAV-extend-hypotheses} 
In the case when $G = \mathrm{Res}_{k^\prime/k}(G^\prime)$ for a finite, purely inseparable extension $k^\prime/k$ and an abelian variety $G^\prime$ over $k^\prime$, it is remarkable that \Cref{thm:pseudo-extend}~\ref{m:PE-b}~\ref{m:PE-ii} holds without assuming that $S$ is geometrically regular. This is the advantage of approaching \ref{m:PE-b}~\ref{m:PE-ii} indirectly, via \ref{m:PE-a}.
\end{remark}

\begin{remark}
The condition on the $\overline{k}$-tori of $G_{\overline{k}}$ is needed in \Cref{thm:pseudo-extend}: without it \eqref{eqn:Gabber-lemma} fails already for some pseudo-finite $G$, see \cite{FG21}*{Example 7.2 (b) and Remark 7.3 (b)}.
\end{remark}

\csub[Pseudo-properness of $G/P$]
\label{sec:GP-pseudo-proper}

We wish to improve a result of Borel--Tits \cite{BT78}*{Proposition 1} that reappeared in \cite{Tit13}*{cours~1992--1993, section~2.5} and \cite{CGP15}*{Proposition C.1.6}: for a pseudo-parabolic subgroup $P$ of a smooth, affine, connected $k$-group $G$, the quotient $G/P$ is not only pseudo-complete as these references showed, but is even pseudo-proper as we now argue. The arguments there were built on the pseudo-completeness of wound unipotent groups and do not generalize because the wound groups that are relevant there are not pseudo-proper (see \Cref{prop:PP-groups}~\ref{m:PPG-ii}). In contrast, we build our argument on the extension result \Cref{thm:pseudo-extend} and the comparison map $i_G$ of \S\ref{pp:pred-iG}.

\begin{theorem} \label{thm:pseudo-flag}
For a connected, smooth, affine $k$-group $G$ and a pseudo-parabolic $k$-subgroup $P \subset G$, the quotient $G/P$ is pseudo-proper and quasi-projective over $k$.
\end{theorem}

\beg \label{eg:pseudo-flag}
If $G = \mathrm{Res}_{k^\prime/k}(G^\prime)$ for a finite field extension $k^\prime/k$ and a reductive $k^\prime$-group $G^\prime$, then $P = \mathrm{Res}_{k^\prime/k}(P^\prime)$ for some parabolic $k^\prime$-subgroup $P^\prime \le G^\prime$ (see \S\ref{pp:pseudo-parabolic}). In this case, \cite{CGP15}*{Corollary~A.5.4~(3)} ensures that $G/P \cong \mathrm{Res}_{k^\prime/k}(G'/P')$, which is pseudo-proper by \Cref{eg:pseudo-proper} and \cite{SGA3IIInew}*{expos\'e XXVI, proposition 1.2}.
\eeg

\begin{proof}[Proof of Theorem \uref{thm:pseudo-flag}]
The quasi-projectivity of $G/P$ follows from the results reviewed in \S\ref{pp:k-groups-lft}. For its pseudo-properness, by \S\ref{pp:pseudo-parabolic}, we may assume that $G$ is pseudo-reductive, of minimal type. We let $i_G$ be the comparison map as in \S\ref{pp:pred-iG}, so that $i_G(G)/i_G(P)$ is pseudo-proper by \Cref{lem:ultra-to-res} and \Cref{eg:pseudo-flag}.  As for $G/P$ itself, we then consider a strictly Henselian discrete valuation ring $R$ that is a geometrically regular $k$-algebra with fraction field $K$ and an $x \in (G/P)(K)$ to be extended to an $R$-point. The target of the map
\be \label{eqn:pseudo-flag-minimal}
G/P  \rightarrow i_G(G)/i_G(P)
\ee
 is pseudo-proper, so there is a unique $y \in (i_G(G)/i_G(P))(R)$ that extends the image of $x$. Since $R$ is strictly Henselian and $i_G(P)$ is smooth, $y$ lifts to an $R$-point of $i_G(G)$, so, by \Cref{lem:minimal-to-ultra}, the $y$-fiber of the map \eqref{eqn:pseudo-flag-minimal} is a generically trivial torsor over $R$ under a $k$-group that is a direct factor of $\mathrm{Ker}(i_G)$, with $x$ as a generic trivialization. Since $\mathrm{Ker}(i_G)$ is pseudo-finite (see \S\ref{pp:pred-iG}), so are its direct factors over $k$. Thus, \Cref{thm:pseudo-extend}~\ref{m:PE-a}~\ref{m:PE-ii} ensures that $x$ extends to an $R$-point of the torsor in question. In particular, $x$ extends to an $R$-point of $G/P$, so that $G/P$ is indeed pseudo-proper by \Cref{rem:geometric-regularity}. 
\end{proof}

\beg
To illustrate how the pseudo-properness of $G/P$ is useful in practice, consider a smooth, integral $k$-curve $C$ with the function field $K$ and a $G$-torsor $E$ over $C$. Sections of $E/P$ correspond to reductions of $E$ to a $P$-torsor, and the pseudo-properness of $G/P$, applied over an \'etale cover of $C$ trivializing $E$, implies that
$$
(E/P)(C) \xrightarrow{\sim} (E/P)(K).
$$
In other words, every generic $P$-reduction of $E$ extends uniquely to a $P$-reduction of $E$ over all of $C$. If we only knew that $G/P$ was pseudo-complete, then we would only know the same extendability for those $P$-reductions of $E$ that are defined over some dense open $C' \subset C$ such that the residue field of every point of $C \setminus C^\prime$ is separable over $k$.
\eeg

\Cref{prop:PP-groups}~\ref{m:PPG-i} and \Cref{thm:pseudo-flag} combine into the following single statement. 

\begin{corollary} \label{cor:GP-pseudo-proper}
For a $k$-group scheme $G$ locally of finite type, an \'etale locally trivial $G$-torsor $E$ over $k$, and a pseudo-parabolic $k$-subgroup $P \le G^{\mathrm{sm},\, \mathrm{lin}}$, the connected components of the quotient $E/P$ are pseudo-proper and quasi-projective.
\end{corollary}

\bpf By considering a finite field extension over which a given connected component of $E$ acquires a rational point, we see that the components of $E/P$ are quasi-compact. Thus, by the results reviewed in \S\ref{pp:k-groups-lft}, they are also quasi-projective. For their remaining pseudo-properness, by \Cref{lem:pseudo-base-change}~\ref{m:PBC-a}, we may base change to $k^s$ to reduce to when $k$ is separably closed, so that $E$ is trivial and $E/P \simeq G/P$.

Since $P$ is smooth (see \S\ref{pp:pseudo-parabolic}), so is the map $G \twoheadrightarrow G/P$. Thus, by \S\ref{pp:Xgred}, the preimage of $(G/P)^{\mathrm{gred}}$ is precisely $G^{\mathrm{gred}}$, that is, $(G/P)^{\mathrm{gred}} \cong G^{\mathrm{gred}}/P$. By \Cref{rem:geometric-regularity}, we may therefore replace $G$ by $G^{\mathrm{gred}}$ to reduce to when $G$ is smooth. Once $G$ is smooth, since $k$ is separably closed, the target of
$$
G/P \twoheadrightarrow G/G^0
$$
is a disjoint union of copies of $\mathrm{Spec}(k)$ and the preimage of each one of them is a copy of $G^0/P$. This reduces us to considering $G^0/P$, that is, we may assume that $G$ is smooth and connected. We then consider the smooth map
$$
G/P \twoheadrightarrow G/G^{\mathrm{sm},\,\mathrm{lin}}
$$
whose target, by \S\ref{pp:fundamental-filtration}~\ref{m:Gsmlin}, is a pseudo-abelian variety. Since pseudo-properness may be tested using strictly Henselian, geometrically regular discrete valuation rings $R$ over $k$ (see \Cref{rem:geometric-regularity}) and, by strict Henselianity, the $R$-fibers of this map are isomorphic to $G^{\mathrm{sm,\, lin}}/P$, the sought pseudo-properness of $G/P$ follows by combining the pseudo-abelian variety case settled in \Cref{prop:PP-groups}~\ref{m:PPG-i} and the connected, smooth, affine case settled in \Cref{thm:pseudo-flag}. 
\epf

\section{Purity and extension results for torsors}
\label{chap:purity} 

Our argument for \Cref{thm:main-GS} is built on purity results for $G$-torsors. The classical cases of these results concern finite schemes, abelian varieties, and reductive groups, and we now generalize to pseudo-finite schemes, pseudo-abelian varieties, and quasi-reductive groups. These generalizations seem subtle, however: the statements require new, sometimes delicate hypotheses and the proofs reach classical cases in somewhat intricate ways. To carry them out, we critically use the comparison maps $i_G$ of \S\S\ref{sec:iG}--\ref{sec:pseudo-reductive}, in particular, we use the new affineness result \Cref{prop:Brauer-affine}.


\csub[Purity for torsors under pseudo-finite groups] \label{sec:purity-pfinite}

A well-known purity result of Moret-Bailly says that for a regular scheme $S$, a closed subset $Z \subset S$ of codimension $\ge 2$, and a finite, flat $S$-group $G$, we have an equivalence of categories
\be \label{eqn:finite-flat-purity}
\{G\text{-torsors over }S\} \xrightarrow{\sim} \{G\text{-torsors over }S \setminus Z\},
\ee
see \cite{MB85b}*{lemme 2} or \cite{flat-purity}*{Theorem 7.1.3}. We extend this to commutative \emph{pseudo-finite} $k$-groups $G$ in \Cref{thm:purity-pseudo-finite} below. Since pseudo-finite groups appear as kernels of the comparison maps $i_G$ presented in \S\ref{pp:pav-iG} and \S\ref{pp:pred-iG}, this purity for pseudo-finite groups is a building block towards similar results for pseudo-abelian varieties and pseudo-reductive groups. The full faithfulness in \eqref{eqn:finite-flat-purity} is a special case of the Hartogs extension principle for sections: for a scheme $S$ that is of depth $\ge 2$ along a closed subset $Z \subset S$, and for an affine $S$-scheme $X$, we have
\be \label{eqn:S-depth-2}
X(S) \xrightarrow{\sim} X(S\setminus Z),
\ee
see, for instance, \cite{flat-purity}*{Lemma 7.2.7 (b)}, or perhaps also \cite{torsors-regular}*{Section 1.3.9} for further review. We now upgrade this Hartogs principle to a similar extendability result for fppf local triviality of torsors in \Cref{lem:Res-fppf-trivial}.

\begin{lemma}\label{lem:B-Res}
For a finite, locally free scheme map $S^\prime \to S$ and a flat, locally finitely presented $S^\prime$-group algebraic space $G$, the counit map $(\mathrm{Res}_{S^\prime/S}(G))_{S^\prime} \to G$ gives a monomorphism of stacks
\begin{equation}\label{eqn:stack-mono}
\mathbf{B}(\mathrm{Res}_{S^\prime/S}(G)) \hookrightarrow \mathrm{Res}_{S^\prime/S}(\mathbf{B}G)
\end{equation}
whose essential image, for a variable $S$-scheme $T$, consists of those $G$-torsors over $T^\prime \ce S^\prime \times_S T$ that trivialize fppf locally over $T$ \up{and not merely over $T^\prime$}. In particular, this monomorphism is an equivalence whenever $G$ is smooth.
\end{lemma}

\begin{proof}
Except for the last assertion, the claims are all special cases of \cite{Gir71}*{chapitre V, proposition~3.1.3} (and of its proof). When $G$ is smooth, to show that every $G$-torsor over $T^\prime$ trivializes fppf locally on $T$, a limit argument allows us to assume that $T$ is strictly Henselian local. Then our $T$-finite $T^\prime$ is a union of strictly Henselian local schemes, so that, since $G$-torsors inherit smoothness (see \S\ref{pp:conv}), they are all trivial over $T^\prime$.
\end{proof}

To analyze \eqref{eqn:stack-mono} beyond smooth groups, we will use the following lemma.

\begin{lemma}\label{lem:Res-fppf-trivial}
For a finite, purely inseparable field extension $k^\prime/k$, a $k^\prime$-group $G$ that is either finite or unipotent, and a $k$-scheme $S$ that is of depth $\geq 2$ along a closed subset $Z \subset S$, a $G$-torsor $E$ over $S^\prime \ce S_{k^\prime}$ trivializes fppf locally on $S$ if and only if its restriction to $(S \setminus Z)_{k^\prime}$ trivializes fppf locally on $S \setminus Z$, in other words, the following square is Cartesian\ucolon
\[
\xymatrix{
\left(\mathbf{B}(\mathrm{Res}_{S^\prime/S}(G_{S^\prime}))\right)(S) \ar@{^(->}[r] \ar@{^(->}[d] & \left(\mathrm{Res}_{S^\prime/S}(\mathbf{B}G_{S^\prime})\right)(S) \ar@{^(->}[d] \\
\left(\mathbf{B}(\mathrm{Res}_{S^\prime/S}(G_{S^\prime}))\right)(S \setminus Z)  \ar@{^(->}[r] & \left(\mathrm{Res}_{S^\prime/S}(\mathbf{B}G_{S^\prime})\right)(S \setminus Z).
}
\]
\end{lemma}

\begin{proof}
The equivalent reformulation in terms of the Cartesian square, as well as the full faithfulness of the horizontal arrows of this square, follows from \Cref{lem:B-Res}. The affineness of $G$, the depth hypothesis, and \eqref{eqn:S-depth-2} give the full faithfulness of the vertical arrows. Moreover, the ``only if'' assertion is clear, so we focus on the ``if,'' for which we may work locally on $S$, and hence, by a limit argument, assume that $S$ is local.

We choose a $k^\prime$-group embedding $G \hookrightarrow \widetilde{G}$ with $\widetilde{G} = \mathrm{GL}_{n,\, k^\prime}$ if $G$ is finite (resp., with $\widetilde{G}$ some $k^\prime$-group of upper unitriangular matrices when $G$ is unipotent) and note that in both cases the homogeneous space $\widetilde{G}/G$ is affine, see \cite{torsors-regular}*{Section 1.3.8}, \S\ref{pp:conv}, and \Cref{lem:unip-quotient}. By \Cref{lem:res-unipotent}, the kernels of the counit maps
\[
c\colon (\mathrm{Res}_{k^\prime/k}(G))_{k^\prime} \twoheadrightarrow G \quad \text{and} \quad \widetilde{c}\colon(\mathrm{Res}_{k^\prime/k}(\widetilde{G}))_{k^\prime} \twoheadrightarrow \widetilde{G}
\]
are unipotent. \Cref{lem:unip-quotient} then shows that $\mathrm{Ker}(\widetilde{c})/\mathrm{Ker}(c)$ is affine. However, then the quotient  $\mathrm{Res}_{k^\prime/k}(\widetilde{G})/\mathrm{Res}_{k^\prime/k}(G)$ is also affine: after base change to $k^\prime$, it suffices to note that the map
\[
\left(\mathrm{Res}_{k^\prime/k}(\widetilde{G})/\mathrm{Res}_{k^\prime/k}(G)\right)_{k^\prime} \cong (\mathrm{Res}_{k^\prime/k}(\widetilde{G}))_{k^\prime}/(\mathrm{Res}_{k^\prime/k}(G))_{k^\prime} \twoheadrightarrow \widetilde{G}/G
\]
is affine because fppf locally on the target its source is isomorphic to the base change of $\mathrm{Ker}(\widetilde{c})/\mathrm{Ker}(c)$. We conclude that both the source and the target of the map
\[
i\colon \mathrm{Res}_{k^\prime/k}(\widetilde{G})/\mathrm{Res}_{k^\prime/k}(G) \hookrightarrow \mathrm{Res}_{k^\prime/k}(\widetilde{G}/G)
\]
are affine, in particular, by \eqref{eqn:S-depth-2}, an $S$-point of $\mathrm{Res}_{k^\prime/k}(\widetilde{G}/G)$ lies in $\mathrm{Res}_{k^\prime/k}(\widetilde{G})/\mathrm{Res}_{k^\prime/k}(G)$ if and only if so does the induced $(S \setminus Z)$-point.

With this in mind, we consider the commutative square
\[
\xymatrix{
\left(\mathrm{Res}_{k^\prime/k}(\widetilde{G})/\mathrm{Res}_{k^\prime/k}(G)\right)(S \setminus Z) \ar[r] \ar@{^(->}[d]^-{i(S\setminus Z)} & H^1(S \setminus Z,\, \mathrm{Res}_{k^\prime/k}(G)) \ar@{^(->}[d] \\
\left(\mathrm{Res}_{k^\prime/k}(\widetilde{G}/G)\right)(S \setminus Z)  \ar[r] & H^1((S \setminus Z)_{k^\prime}, G)
}
\]
whose horizontal arrows are the connecting maps of long exact cohomology sequences and vertical maps are injective by the above and by \Cref{lem:B-Res}. Since $S$ is local, $S^\prime$ is semilocal, so all $\widetilde{G}$-torsors over $S^\prime$ are trivial. In particular, $E$ induces a trivial $\widetilde{G}$-torsor, so, since $E|_{(S\setminus Z)_{k^\prime}}$ trivializes fppf locally on $S \setminus Z$, \Cref{lem:B-Res} and \cite{Gir71}*{chapitre III, proposition 3.2.2} ensure that the class of $E|_{(S\setminus Z)_{k^\prime}}$ comes from some $\alpha \in \left(\mathrm{Res}_{k^\prime/k}(\widetilde{G})/\mathrm{Res}_{k^\prime/k}(G)\right)(S\setminus Z)$  that is unique up to the left action of $(\mathrm{Res}_{k^\prime/k}(\widetilde{G}))(S\setminus Z)$. Similarly, the analogous square over $S$ shows that the class of $E$ in $H^1(S^\prime, G)$ comes from some $\beta \in \left(\mathrm{Res}_{k^\prime/k}(\widetilde{G}/G)\right)(S)$ that is unique up to the left action of $(\mathrm{Res}_{k^\prime/k}(\widetilde{G}))(S)$. However, $\mathrm{Res}_{k^\prime/k}(\widetilde{G})$ is affine, so $(\mathrm{Res}_{k^\prime/k}(\widetilde{G}))(S) = (\mathrm{Res}_{k^\prime/k}(\widetilde{G}))(S \setminus Z)$ by \eqref{eqn:S-depth-2}. In particular, we may arrange that $\beta|_{S \setminus Z} = \alpha$. At this point, the conclusion of the previous paragraph implies that $\beta$ is actually an $S$-point of $\mathrm{Res}_{k^\prime/k}(\widetilde{G})/\mathrm{Res}_{k^\prime/k}(G)$. Then, however, $\beta$ lifts to an $S$-point of $\mathrm{Res}_{k^\prime/k}(\widetilde{G})$ fppf locally on $S$, to the effect that $E$ trivializes fppf locally on $S$, as desired.
\end{proof}

We now extend the purity \eqref{eqn:finite-flat-purity} to torsors under pseudo-finite groups.

\begin{theorem}\label{thm:purity-pseudo-finite}
Let $S$ be a geometrically regular $k$-scheme, let $Z \subset S$ be a closed subset of codimension $\geq 2$, and let $G$ be a $k$-group scheme locally of finite type with $G^0$ pseudo-finite and commutative. For every gerbe $\mathscr{B}$ over $S$ isomorphic to $\mathbf{B}G$ \'etale locally on $S$, we have 
\[
\mathscr{B}(S) \xrightarrow{\sim} \mathscr{B}(S\setminus Z),
\]
in particular, for every $S$-group $\mathscr{G}$ isomorphic to $G$ \'etale locally on $S$, we have
\[
H^1(S, \mathscr{G}) \xrightarrow{\sim} H^1(S\setminus Z, \mathscr{G}) \quad \text{and, if } \mathscr{G} \text{ is commutative, also} \quad H^2(S, \mathscr{G}) \hookrightarrow H^2(S\setminus Z, \mathscr{G}).
\]
\label{eqn:purity-pseudo-finite}
\end{theorem}

\begin{proof}
It suffices to settle the claim about $\mathscr{B}$: the conclusions about $H^1$ (resp., $H^2$) follow from the rest by choosing $\mathscr{B} = \mathbf{B}\mathscr{G}$ (resp., by letting $\mathscr{B}$ be a $\mathscr{G}$-gerbe over $S$ that trivializes over $S\setminus Z$ and using \cite{SP}*{Lemma~\href{https://stacks.math.columbia.edu/tag/0DLS}{0DLS}} to see that this $\mathscr{B}$ is isomorphic to $\mathbf{B}G$ \'etale locally on $S$). By descent, the claim about $\mathscr{B}$ is \'etale local on $S$, so we may assume that $\mathscr{B} = \mathbf{B}G$ and, by also working \'etale locally and combining Noetherian induction with spreading out, that $S$ is strictly Henselian local and $Z$ is its closed point.\footnote{This specific assumption on $Z$ will not be used in this proof, but it is convenient to do all the straightforward preliminary reductions at once here and then simply also refer to them in the later proofs of \Cref{thm:purity-pproper,thm:purity-pcomplete}.} By \Cref{cor:pfinite-affine} (and \cite{SP}*{Lemma~\href{https://stacks.math.columbia.edu/tag/02L5}{02L5}}), the connected components of $G$ are affine, so, by \eqref{eqn:S-depth-2} and fpqc descent, the map $\mathbf{B}G(S) \to \mathbf{B}G(S\setminus Z)$ is fully faithful. Thus, we only need to show that every $G$-torsor over $S\setminus Z$ extends to a $G$-torsor over $S$ (necessarily uniquely up to a unique isomorphism).

For this, we first treat the case when $G^0 = 1$, that is, when $G$ is twisted constant (see \S\ref{pp:fundamental-filtration}~\ref{m:G0}). Due to its strict Henselianity, $S$ is automatically a $k^s$-scheme, so this case follows from the purity for the \'etale fundamental group \cite{SGA2new}*{expos\'e X, th\'eor\`eme 3.4} (with \cite{SGA3II}*{expos\'e X, corollaire~5.14} to ensure that every connected component of every torsor under a constant group over $S\setminus Z$ is finite \'etale).

In general, by applying the settled twisted constant case to $G/G^0$, and since $S$ is strictly Henselian, we find that every $(G/G^0)$-torsor over $S\setminus Z$ is trivial. In particular, every $G$-torsor over $S\setminus Z$ reduces to a $G^0$-torsor over $S\setminus Z$, so that we may assume that $G$ is connected, and so also pseudo-finite and commutative. With these assumptions, we will even prove that $H^2_Z(S,G) \cong 0$. For this, since such vanishing is stable under extensions, \S\ref{pp:gred-filtration} and the pseudo-finiteness of $G$ allow us to assume that $G \leq \mathrm{Res}_{k^\prime/k}Q$ for a finite, purely inseparable field extension $k^\prime/k$ and a finite, commutative $k^\prime$-group~$Q$. After base change to $k^\prime$ and using \Cref{lem:res-unipotent}, \S\ref{pp:conv}, and the finiteness of $Q$, the $k$-group $(\mathrm{Res}_{k^\prime/k}Q)/G$ is affine. Thus, since $S$ is of depth $\geq 2$ along $Z$, the Hartogs extension for sections \eqref{eqn:S-depth-2} gives $H^1_Z(S, (\mathrm{Res}_{k^\prime/k}Q)/G) \cong 0$. The exact sequence
\[
H^1_Z(S, (\mathrm{Res}_{k^\prime/k}Q)/G) \to H^2_Z(S, G) \to H^2_Z(S, \mathrm{Res}_{k^\prime/k}Q)
\]
then reduces the desired vanishing of $H^2_Z(S, G)$ to the case when $G = \mathrm{Res}_{k^\prime/k}Q$. Moreover, the vanishing of this $H^2_Z$ follows from the equivalence of categories $\mathscr{B}(S) \xrightarrow{\sim} \mathscr{B}(S\setminus Z)$ for gerbes $\mathscr{B}$ over $S$ that become isomorphic to $\mathbf{B}G$ \'etale locally on $S$. Thus, by repeating the same reductions as at the beginning of the proof, we are left with showing that every $G$-torsor over $S\setminus Z$ extends (necessarily uniquely up to a unique isomorphism) to a $G$-torsor over $S$. Equivalently, by \Cref{lem:B-Res}, since $G \cong \mathrm{Res}_{k^\prime/k}(Q)$, we need to show that every $Q$-torsor over $(S\setminus Z)_{k^\prime}$ that trivializes fppf locally on $S\setminus Z$ extends to a $Q$-torsor over $S_{k^\prime}$ that trivializes fppf locally on $S$. Since $Q$ is finite, $S_{k^\prime}$ is regular, and $Z_{k^\prime} \subset S_{k^\prime}$ is of codimension $\geq 2$, the purity result \eqref{eqn:finite-flat-purity} supplies the extension to a $Q$-torsor over $S_{k^\prime}$. It then remains to note that, by \Cref{lem:Res-fppf-trivial}, this torsor automatically trivializes fppf locally on $S$ because the same holds over $S\setminus Z$.
\end{proof}

As we now show, purity for torsors under pseudo-finite groups of \Cref{thm:purity-pseudo-finite} implies that the map on $H^2$ in \eqref{eqn:purity-pseudo-finite} need not be surjective even when $G$ is finite. This also shows that the assumptions of \cite{flat-purity}*{Theorem 1.1.1} are sharp.

\begin{corollary}
If $k$ is imperfect of characteristic $p$, then for every $2$-dimensional, geometrically regular, local $k$-algebra $(R, \mathfrak{m})$, its punctured spectrum $U_R \ce \mathrm{Spec}(R) \setminus \mathfrak{m}$, and a $k$-group $G$ that is either $\alpha_p$ or a nontrivial commutative extension of $\mathbb{Z}/p\mathbb{Z}$ by $\mu_p$, there is a $G$-gerbe over $U_R$ that does not extend to a $G$-gerbe over $R$, in particular, $H^3_{\mathfrak{m}}(R, G) \not\cong 0$.
\end{corollary}

\begin{proof}
By \cite{SP}*{Lemma~\href{https://stacks.math.columbia.edu/tag/0AVZ}{0AVZ}}, since $R$ is $2$-dimensional, we have $H^2_{\mathfrak{m}}(R) \not\cong 0$, equivalently, there is a nontrivial $\mathbb{G}_a$-torsor $E$ over $U_R$. On the other hand, \cite{Tot13}*{Lemmas 6.3 and 7.1} supply a commutative, pseudo-finite extension $\widetilde{G}$ of $\mathbb{G}_a$ by $G$ over $k$.  By \Cref{thm:purity-pseudo-finite}, every $\widetilde{G}$-torsor over $U_R$ extends to a $\widetilde{G}$-torsor over $R$. Thus, since $E$ does not extend to a $\mathbb{G}_a$-torsor over $R$ (else it would be trivial), it does not lift to a $\widetilde{G}$-torsor. In particular, via the long exact cohomology sequence associated to the extension, $E$ gives rise to a nontrivial $G$-gerbe over $U_R$.  By the cohomology sequence and \Cref{thm:purity-pseudo-finite} again, this $G$-gerbe does not extend to a $G$-gerbe over $R$.
\end{proof}

\csub[Purity for torsors under pseudo-proper groups]\label{sec:purity-pproper}

Purity for torsors under abelian schemes says that for a regular scheme $S$, a closed subset $Z \subset S$ of codimension $\ge 2$, and an abelian scheme $G$ over $S$, we have an equivalence of categories
\be\label{eqn:ab-sch-purity}
\{G\text{-torsors over }S\} \xrightarrow{\sim} \{G\text{-torsors over }S \setminus Z\}.
\ee
 Indeed, \eqref{eqn:abelian-scheme-extend} supplies the full faithfulness even with $S \setminus Z$ replaced by the union of the generic points of $S$ and so also shows that every $G$-torsor over $S \setminus Z$ has finite order, at which point the extendability to $S$ of a relevant $G[n]$-torsor over $S \setminus Z$ for some $n > 0$ follows from \eqref{eqn:finite-flat-purity}. 

We now extend \eqref{eqn:ab-sch-purity} to analogous purity for torsors under pseudo-abelian varieties, more generally, for torsors under smooth, pseudo-proper $k$-groups.

\begin{theorem}\label{thm:purity-pproper}
Let $S$ be a geometrically regular $k$-scheme, let $Z \subset S$ be a closed subset of codimension $\geq 2$, and let $G$ be a $k$-group scheme locally of finite type with $G^0$ pseudo-proper \up{see Proposition~\uref{prop:PP-groups}~\uref{m:PPG-ii}} and either smooth or commutative. For every gerbe $\mathscr{B}$ over $S$ isomorphic to $\mathbf{B}G$ \'etale locally on $S$, we have 
\[
\mathscr{B}(S) \xrightarrow{\sim} \mathscr{B}(S\setminus Z),
\]
in particular, for every $S$-group $\mathscr{G}$ isomorphic to $G$ \'etale locally on $S$, we have
\be\label{eqn:purity-pproper}
H^1(S, \mathscr{G}) \xrightarrow{\sim} H^1(S\setminus Z, \mathscr{G}) \quad \text{and, if } \mathscr{G} \text{ is commutative, also} \quad H^2(S, \mathscr{G}) \hookrightarrow H^2(S\setminus Z, \mathscr{G}).
\ee
\end{theorem}

\begin{proof}
As in the proof of \Cref{thm:purity-pseudo-finite}, we may assume that $S$ is strictly Henselian local and that $Z$ is its closed point. As there, we need to show that every $G$-torsor (resp., every isomorphism of $G$-torsors) over $S\setminus Z$ extends uniquely to a $G$-torsor (resp., an isomorphism of $G$-torsors) over $S$. Since $S$ is strictly Henselian, it is automatically a $k^s$-scheme, so every $(G/G^0)$-torsor over $S$ is a disjoint union of copies of $S$ (see \S\ref{pp:fundamental-filtration}~\ref{m:G0}) and every $G$-torsor over $S$ is then, as an $S$-scheme, a disjoint union of $G^0$-torsors. By applying \Cref{thm:purity-pseudo-finite} to the \'etale $k$-group $G/G^0$, we may therefore pass to $G^0$ to assume that $G$ is connected. (To extend $G$-torsor isomorphisms, we note that, by schematic density of $S \setminus Z$ in $S$, it suffices to extend them as scheme isomorphisms.)

The case when our connected $G$ is commutative reduces to when it is smooth: indeed, then torsor isomorphisms extend by \Cref{thm:pseudo-extend}~\ref{m:PE-ii}, the smooth $k$-subgroup $G^{\mathrm{gred}}$ inherits pseudo-properness, and, by \Cref{eg:pseudo-finite-G}, the quotient $G/G^{\mathrm{gred}}$ is pseudo-finite, so \Cref{thm:purity-pseudo-finite} and \eqref{eqn:purity-pproper} applied to $G^{\mathrm{gred}}$ (including the aspect about $H^2$) imply the claim for $G$. Once our connected $G$ is smooth, since $S$ is strictly Henselian, every $G$-torsor over $S$ is trivial. In particular, the desired conclusion is stable under extensions and we may assume that $G$ is either a pseudo-abelian variety or a connected, smooth, strongly wound unipotent group (see \S\ref{pp:fundamental-filtration}~\ref{m:Gsmlin} and \Cref{prop:PP-groups}~\ref{m:PPG-ii}). Moreover, by \Cref{thm:pseudo-extend}~\ref{m:PE-ii}, it even suffices to show that every $G$-torsor over $S\setminus Z$ is generically trivial.

In the pseudo-abelian variety case, we consider the comparison map of \S\ref{pp:pav-iG}:
\[
i_{G}\colon G \to \mathrm{Res}_{k^\prime/k}(\overline{G}),
\]
and we aim directly for the a priori stronger conclusion that $H^2_Z(S, G) \cong 0$. For this, since $\mathrm{Ker}(i_G)$ is commutative and pseudo-finite, \Cref{thm:purity-pseudo-finite} allows us to replace $G$ by $i_G(G)$; in other words, we may assume that $i_G$ is injective. At this point, since $\mathrm{Coker}(i_G)$ is affine, \eqref{eqn:S-depth-2} gives $H^1_Z(S, \mathrm{Coker}(i_G)) \cong 0$, to the effect that we may assume that $G = \mathrm{Res}_{k^\prime/k}(\overline{G})$.  By the same reasoning as in the proof of \Cref{thm:purity-pseudo-finite} and by \Cref{lem:B-Res}, it now suffices to show that every $\overline{G}$-torsor over $(S\setminus Z)_{k^\prime}$ extends to a $\overline{G}$-torsor over $S_{k^\prime}$. However, $\overline{G}$ is an abelian variety, $S_{k^\prime}$ is regular, and $Z_{k^\prime} \subset S_{k^\prime}$ is of codimension $\geq 2$, so this extendability follows from \eqref{eqn:ab-sch-purity}.

In the remaining unipotent case, by \Cref{lem:sw-cckp}, we may assume that our connected, smooth, strongly wound, unipotent $k$-group $G$ is commutative and $p$-torsion. Moreover, by the Popescu theorem (see \S\ref{pp:conv}), we may assume that $S$ is the strict Henselization of a smooth $k$-scheme at a point. By excision \cite{torsors-regular}*{Proposition 4.2.1}, the question of extending $G$-torsors over $Z$ only depends on the base change to the completion of $S$ along $Z$. Thus, by passing to this completion, we may assume that $S$ is not only strictly Henselian but also complete.  At this point, the desired generic triviality of $G$-torsors over $S\setminus Z$ becomes a special case of \Cref{lem:SH-trivial}~\ref{m:SHT-ii}.
\end{proof}

\begin{lemma}\label{lem:SH-trivial}
Let $p$ be the characteristic exponent of $k$, let $S$ be a locally Noetherian $k$-scheme, let $Z \subset S$ be a closed subset of codimension $\geq 2$, and let $G$ be a unipotent $k$-group scheme.
\benum
    \item\label{m:SHT-i} If $S$ is of depth $\geq 3$ along $Z$ \up{so that $Z \subset S$ is of codimension $\geq 3$}, then every $G$-torsor over $S \setminus Z$ extends uniquely to a $G$-torsor over $S$, in fact, for every $S$-gerbe $\mathscr{B}$ isomorphic to $\mathbf{B}G$ \'etale locally on $S$, we have $\mathscr{B}(S) \xrightarrow{\sim} \mathscr{B}(S\setminus Z)$.
    \item\label{m:SHT-ii} If $S$ is geometrically regular, strictly Henselian, local, $Z$ is its closed point, $S$ is complete along $Z$, and $G$ is connected, smooth, commutative, and $p$-torsion, then every $G$-torsor over $S\setminus Z$ is generically trivial.
    \item\label{m:SHT-iii} If $k$ is separably closed, $S \cong \mathrm{Spec}(\widehat{\mathscr{O}}_{\mathbb{A}^2_k,\, z})$ for a closed point $z = Z$ of $\mathbb{A}^2_k \cong \mathrm{Spec}(k[s,t])$ cut out set-theoretically by $s^{p^n} - a$ and $t^{p^m} - b$ for some $a, b \in k$, $n,m \geq 0$, and $G$ is connected, smooth, commutative, and $p$-torsion, then every $G$-torsor over $S\setminus z$ trivializes both over $S[\frac{1}{s^{p^n}-a}]$ and over $S[\frac{1}{t^{p^m}-b}]$.
\end{enumerate}
\end{lemma}

\begin{proof}
Both the assumptions and the conclusion of \ref{m:SHT-iii} are strictly stronger than those of \ref{m:SHT-ii}.
\benum
    \item As in the proof of \Cref{thm:purity-pseudo-finite}, it suffices to show that $G$-torsors extend, and we may assume that $S$ is local. Since $G$ is unipotent, it is a $k$-subgroup of some upper unitriangular matrix $k$-group $U$, and $U/G$ is affine; its sections correspond to $G$-reductions of the trivial $U$-torsor, see \S\ref{pp:conv}. This affineness and \eqref{eqn:S-depth-2} give $(U/G)(S) \xrightarrow{\sim} (U/G)(S\setminus Z)$. Thus, since every $U$-torsor over our local $S$ is trivial, we may assume that $G = U$. We then need to show that every $G$-torsor over $S\setminus Z$ is trivial and, by passing to subquotients of $G$, we reduce to $G = \mathbb{G}_a$. It then suffices to note that $H^1(S\setminus Z, \mathscr{O}_{S\setminus Z}) \cong H^2_Z(S, \mathscr{O}_S) \cong 0$ because $S$ is Noetherian local of depth $\geq 3$.

    \item By the Popescu theorem (see \S\ref{pp:conv}) and a limit argument, we may assume that $S$ is the completion of the strict Henselization of a smooth, irreducible $k$-scheme $X$ of dimension $d \geq 0$ at the generic point $y$ of some irreducible closed $Y \subset X$ of codimension $\geq 2$. Moreover, by \ref{m:SHT-i}, we may assume that $S$ is $2$-dimensional, so that $Y$ is of codimension $2$.    By the presentation theorem \cite{CTHK97}*{Theorem 3.1.1}, at the cost of shrinking $X$ around $y$, we may find a smooth $k$-morphism $\psi\colon X \to \mathbb{A}^{d-1}_k$ of relative dimension $1$ that makes $Y$ finite over $\mathbb{A}^{d-1}_k$. Since $Y \subset X$ is of codimension $2$, its image $\psi(Y) \subset \mathbb{A}^{d-1}_k$ is a divisor with the generic point $\psi(y)$.     By the presentation theorem again, now applied with $\psi(Y) \subset \mathbb{A}^{d-1}_k$ in place of $Y \subset X$, after shrinking $X$ around $y$ we may find a smooth $k$-morphism $\psi'\colon X \to \mathbb{A}^{d-2}_k$ of relative dimension $2$ that makes $Y$ finite over $\mathbb{A}^{d-2}_k$. Now $\psi'(y)$ is the generic point of $\mathbb{A}^{d-2}_k$, so we may replace $k$ by $k(t_1, \dotsc, t_{d-2})$ (see also \S\ref{pp:fundamental-filtration}~\ref{m:RuskG}) to reduce to the case when $X$ is of dimension $2$ and $Y$ is its closed point $y$.    Moreover, since the strict Henselization of $X$ at $y$ does not change if we consider $X_{k^s}$ instead, we may also assume that $k$ is separably closed.

    We apply the presentation theorem \cite{CTHK97}*{Theorem 3.1.1} one final time to find, after shrinking $X$ around $y$, an \'etale $k$-morphism $\varphi\colon X \to \mathbb{A}^2_k$ with $k_{\varphi(y)} \xrightarrow{\sim} k_y$, so also with $\widehat{\mathscr{O}}_{\mathbb{A}^2_k,\, \varphi(y)} \cong \widehat{\mathscr{O}}_{X,\, y}$.     By its universal property, the completion of the strict Henselization of a Noetherian local ring is insensitive to first passing to an initial completion, so our $S$ is also the completion of the strict Henselization of $\mathbb{A}^2_k$ at $\varphi(y)$. Moreover, since $k$ is separably closed, the strict Henselization is superfluous.     All in all, we may assume that $S \cong \mathrm{Spec}(\widehat{\mathscr{O}}_{\mathbb{A}^2_k,\, z})$ for some closed point $z \in \mathbb{A}^2_k$. The extension $k_z/k$ is purely inseparable, so, letting $s$ and $t$ be the standard coordinates on $\mathbb{A}^2_k$, there are some $a, b \in k$ and $n, m \geq 0$ for which the Artinian local closed subscheme $\mathrm{Spec}(k[s,t]/(s^{p^n} - a, t^{p^m} - b)) \subset \mathbb{A}^2_k$ is an infinitesimal thickening of $z$.   We are now in the setting of \ref{m:SHT-iii}, so it suffices to settle the latter.
    
    \m
    By \S\ref{pp:wound-structure}, our $k$-group $G$ is given by the vanishing locus of some $p$-polynomial
\[
\tst \qq F = \sum_{i=0}^N \sum_{j=0}^{n_i} f_{ij} u_i^{p^j} \in k[u_0, \dotsc, u_N]
\]
with $f_{in_i} \neq 0$ whenever $n_i \geq 0$. By symmetry, it suffices to argue that every $G$-torsor over $S \setminus z$ trivializes over $S[\frac{1}{t^{p^m} - b}]$.

We set $B \ce \widehat{\mathscr{O}}_{\mathbb{A}^2_k,\, z}$, so that concretely $B$ is the $(s^{p^n} - a, t^{p^m} - b)$-adic completion of $k[s,t]$. Moreover, let $A$ be the $(s^{p^n} - a)$-adic completion of $k[s]$, so that $B$ is the $(t^{p^m} - b)$-adic completion of $A[t]$. Concretely, $A$ is the ring of formal series
\be\label{eqn:a-series}
 \textstyle \qqqq \sum_{i \geq 0} (a_{i,0} + a_{i,1}s + \cdots + a_{i,p^n-1}s^{p^n-1})(s^{p^n} - a)^i \quad \text{with} \quad a_{i,i'} \in k,
\ee
where multiplication uses iterative division by $s^{p^n} - a$ until the residue polynomial has degree $\leq p^n-1$. Similarly, $B$ is the ring of formal series
\[
\textstyle \q \sum_{j \geq 0} (b_{j,0} + b_{j,1}t + \cdots + b_{j,p^m-1}t^{p^m-1})(t^{p^m} - b)^j \quad \text{with} \quad b_{j,j'} \in A.
\]
We consider the fraction fields $K_A$ and $K_B$ of $A$ and $B$, respectively, as well as the $(t^{p^m} - b)$-adic completion $C$ of $K_A[t]$. Similarly, $C$ is the ring of formal power series
\[
\textstyle \qq \sum_{j \geq 0} (c_{j,0} + c_{j,1}t + \cdots + c_{j,p^m-1}t^{p^m-1})(t^{p^m} - b)^j \quad \text{with} \quad c_{j,j'} \in K_A.
\]
By excision \cite{torsors-regular}*{Proposition 4.2.1}, a $G$-torsor over $S\setminus z$ amounts to a $G$-torsor over $B[\frac{1}{t^{p^m} - b}]$ whose restriction to $C[\frac{1}{t^{p^m} - b}]$ extends to a $G$-torsor over $C$. Via the explicit description of $G$-torsors in terms of $F$ (see \eqref{eqn:unipotent-coho}), this amounts to an $\alpha \in B[\frac{1}{t^{p^m} - b}]$ and a $\beta \in C$ for which we have $\alpha - \beta = F(\gamma)$ for some $\gamma \in (C[\frac{1}{t^{p^m} - b}])^{N+1}$. The coordinates of $\gamma$ are formal series $\sum_{j \in \mathbb{Z}} (c_{j,0} + c_{j,1}t + \cdots + c_{j,p^m-1}t^{p^m-1})(t^{p^m} - b)^j$ with $c_{j,j'} \in K_A$ and only finitely many $c_{j,j'}$ nonzero for $j<0$. Collecting the terms with $j\geq 0$ gives an element $\gamma' \in C^{N+1}$. By replacing $\beta$ with $\beta + F(\gamma')$ (see \eqref{eqn:unipotent-coho}), we reduce to the case when $c_{j,j'} = 0$ for all $j \geq 0$. In turn, each $c_{j,j'}$ is a formal series $\sum_{i \in \mathbb{Z}} (a_{i,0} + a_{i,1}s + \cdots + a_{i,p^n-1}s^{p^n-1})(s^{p^n} - a)^i$ with $a_{i,i'} \in k$ and only finitely many $a_{i,i'}$ nonzero for $i<0$. Collecting the terms with $i\geq 0$ for all the nonzero $c_{j,j'}$ with $j<0$ gives an element $\gamma'' \in (A[\frac{1}{t^{p^m}-b}])^{N+1}$. Replacing $\alpha$ by $\alpha - F(\gamma'')$ reduces us further to the case when $c_{j,j'} = 0$ whenever $j\geq 0$ and also $a_{i,i'} = 0$ for $i\geq 0$. At this point, each coordinate of $\gamma$ is a finite sum of terms $a_{i,i'} s^{i'} t^{j'} (s^{p^n} - a)^i (t^{p^m} - b)^j$ with $i, j < 0$ and $i'<p^n$, $j'<p^m$. Consequently, $F(\gamma)$ is a finite sum of terms of this form because $F$ is a $p$-polynomial and, for any $r\geq 0$, $i' p^r < p^n p^r$ and $j' p^r < p^m p^r$, ensuring that division by $(s^{p^n}-a)$ and $(t^{p^m}-b)$ only involves powers with negative exponents. On the other hand, $\beta$ involves only nonnegative powers of $(t^{p^m} - b)$, while $\alpha$ only involves coefficients $b_{i,i'}$ whose formal expansions only use nonnegative powers of $(s^{p^n} - a)$. Thus, $\alpha - \beta = F(\gamma)$ is possible only if $F(\gamma) = 0$, that is, $\alpha = \beta$. Then $\alpha \in B$ by \eqref{eqn:S-depth-2}. In other words, the restriction of our $G$-torsor to $B[\frac{1}{t^{p^m} - b}]$ extends to $B$ itself. Since $B$ is strictly Henselian and $G$ is smooth, it follows that every $G$-torsor over $S\setminus z$ trivializes over $S[\frac{1}{t^{p^m} - b}]$, as desired.
 \qedhere
\eenum
\end{proof}

\csub[Purity for torsors under pseudo-complete groups]
\label{sec:purity-pcomplete}

Purity for torsors under wound unipotent groups that we establish in \Cref{thm:purity-pcomplete} is a genuinely new phenomenon that has no analogue over perfect fields. Indeed, for $\mathbb{G}_a$ we have
\be\label{eqn:A2}
\{\mathbb G_a\text{-torsors over }\mathbb A^2_k\} \xrightarrow{\not\sim} \{\mathbb G_a\text{-torsors over }\mathbb A^2_k \setminus \{(0, 0) \}\}
\ee
 because $\mathbb A^2_k \setminus \{(0, 0) \}$ is not affine. In contrast, this equivalence holds with $\mathbb{G}_a$ replaced by any smooth, wound, unipotent $k$-group, or even any smooth, pseudo-complete $k$-group, due to the following extension of \Cref{thm:purity-pproper} to pseudo-complete groups.

\begin{theorem}\label{thm:purity-pcomplete}
Let $S$ be a geometrically regular $k$-scheme, let $Z \subset S$ be a closed subset of codimension $\geq 2$, and let $G$ be a $k$-group scheme locally of finite type with $G^0$ pseudo-complete \up{for instance, wound unipotent, see Proposition~\uref{prop:PP-groups}~\uref{m:PPG-i}} and either smooth or commutative. If either
\benumr
    \item\label{m:PPC-i} every $z \in Z$ of codimension $2$ in $S$ lies in a geometrically regular $k$-subscheme $S_z \subset S$ of codimension $> 0$ \up{when $k_z/k$ is separable, we may take $S_z = z$}; or
    \item\label{m:PPC-ii} $G^0 \simeq \mathrm{Res}_{\ell/k}(C/T)$ for a finite field extension $\ell/k$, a commutative pseudo-reductive $\ell$-group $C$, and a maximal $\ell$-torus $T \subset C$ \up{all such $G^0$ are wound unipotent}; or
    \item\label{m:PPC-iii} $G^0$ is wound unipotent and $[k : k^p] = p$, where $p$ is the characteristic exponent of $k$; or
    \item\label{m:PPC-iv} $G_{k^s}^0$ is unirational and contains no nonzero proper unirational $k^s$-subgroups;
\end{enumerate}
then, for every gerbe $\mathscr{B}$ over $S$ isomorphic to $\mathbf{B}G$ \'etale locally on $S$, we have
\[
\mathscr{B}(S) \xrightarrow{\sim} \mathscr{B}(S \setminus Z),
\]
in particular, for every $S$-group $\mathscr{G}$ isomorphic to $G$ \'etale locally on $S$, we have
\[
H^1(S, \mathscr{G}) \xrightarrow{\sim} H^1(S\setminus Z, \mathscr{G}) \quad \text{and, if }\mathscr{G}\text{ is commutative, also} \quad H^2(S, \mathscr{G}) \hookrightarrow H^2(S\setminus Z, \mathscr{G}).
\]
\end{theorem}
The cases \ref{m:PPC-ii}--\ref{m:PPC-iv} show that the purity conclusion of Theorem~\uref{thm:purity-pproper} may continue to hold for some $k$-groups $G$ that are not pseudo-proper.

\begin{proof}
As in the proof of \Cref{thm:purity-pseudo-finite}, we may assume that $S$ is strictly Henselian local and $Z$ is its closed point $z$. As there, we need to show that every $G$-torsor over $S \setminus z$ extends uniquely (up to a unique isomorphism) to a $G$-torsor over $S$, and likewise for torsor isomorphisms. As in the proof of \Cref{thm:purity-pproper}, we may assume that $G$ is connected. Moreover, as there, the case when our connected $G$ is commutative reduces to the smooth case. In the latter, we combine the stability of the desired conclusion under extensions with \Cref{thm:purity-pseudo-finite,thm:purity-pproper} (and \Cref{prop:PP-groups}) to replace $G$ by $G^{\mathrm{sm},\,\mathrm{lin}}$ and reduce to when $G$ is connected, wound, and unipotent. At this point $G$ is affine, so the uniqueness aspect follows from \eqref{eqn:S-depth-2}. Since $G$ is even unipotent, for the remaining triviality of $G$-torsors over $S \setminus z$, \Cref{lem:SH-trivial}~\ref{m:SHT-i} allows us to assume that $S$ is $2$-dimensional.

\benumr
    \item 
The triviality of $G$-torsors over $S \setminus z$ is stable under extensions, so we may pass to the subquotients of the filtration of $G$ by its iterated cc$k$p kernels (see \S\ref{pp:fundamental-filtration}~\ref{m:cckp}) to reduce to when $G$ is commutative and $p$-torsion.

We let $R$ be the coordinate ring of $S$ and fix a geometrically regular, closed $S' \subset S$ of codimension $> 0$ that is the $S_z$ in the assumption \ref{m:PPC-i}. The Popescu theorem (see \S\ref{pp:conv}) ensures that $R$ is a filtered direct limit of local rings $R_i$ of smooth $k$-schemes with local transition maps. Since the transition maps are local, for large $i$, a regular sequence that cuts out $S' \subset S$ descends to a regular sequence in $R_i$. This descended sequence cuts out a closed subscheme $S'_i \subset \mathrm{Spec}(R_i)$ of codimension $> 0$ that is geometrically regular (see \cite{SP}*{Lemma~\href{https://stacks.math.columbia.edu/tag/0381}{0381}}): indeed, $S'_{k^\prime} \subset \mathrm{Spec}(R \otimes_k k^\prime)$ is regular for every finite, purely inseparable field extension $k^\prime/k$, which only happens if $(S'_i)_{k^\prime} \subset \mathrm{Spec}(R_i \otimes_k k^\prime)$ is regular. Thus, a limit argument allows us to replace $R$ by some $R_i$ to reduce to the case when $R$ is a local ring of a smooth $k$-scheme, at the cost of no longer assuming that $R$ is strictly local and instead needing to show that every $G$-torsor over $S \setminus z$ extends to a $G$-torsor over $S$.

To simplify $R$ further, we now use our assumption on $S'$. By spreading out, $S' \subset S$ is a localization of a closed immersion $\mathcal{S}' \subset \mathcal{S}$ of irreducible smooth $k$-schemes, and we let $s\colon \mathscr{O}_{\mathcal{S},\, z} \twoheadrightarrow \mathscr{O}_{\mathcal{S}',\, z}$ be the resulting surjection of local rings (with $R = \mathscr{O}_{\mathcal{S},\, z}$). By \cite{EGAIV4}*{corollaire~17.11.3} and the Jacobi criterion \cite{BLR90}*{Section 2.2, Proposition 7}, there is a local ring $A$ of some maximal ideal of some $k[T_1, …, T_d]$ for which $\mathscr{O}_{\mathcal{S},\, z}$ is a localization of a smooth $A$-algebra in such a way that $\mathscr{O}_{\mathcal{S}',\, z}$ is an \'etale $A$-algebra. By base change, the $\mathscr{O}_{\mathcal{S}',\, z}$-algebra $\widetilde{R} := \mathscr{O}_{\mathcal{S}',\, z} \otimes_{A} \mathscr{O}_{\mathcal{S},\, z}$ is then a localization of a smooth $\mathscr{O}_{\mathcal{S}',\, z}$-algebra and comes equipped with a ``diagonal'' section $\widetilde{s} \colon \widetilde{R} \twoheadrightarrow \mathscr{O}_{\mathcal{S}',\, z}$. Since $\widetilde{s}$ factors through the quotient $\widetilde{R} \twoheadrightarrow \mathscr{O}_{\mathcal{S}',\, z} \otimes_{A} \mathscr{O}_{\mathcal{S}',\, z}$, whose target is \'etale over $\mathscr{O}_{\mathcal{S}', \, z}$, and so has the target of $\widetilde{s}$ as a direct factor, we conclude from \cite{EGAIV4}*{corollaire 17.9.5} that the \'etale map $\mathscr{O}_{\mathcal{S},\, z} \rightarrow  \widetilde{R}$ induces an isomorphism between the completion of $\mathscr{O}_{\mathcal{S},\, z}$ along the kernel of $s$ and that of $\widetilde{R}$ along the kernel of $\widetilde{s}$. By \cite{BLR90}*{Section 3.1, Proposition 2},  the latter completion is isomorphic to $\mathscr{O}_{\mathcal{S}',\, z}\llb t_1, \dotsc, t_n\rrb $. Thus, the completion of $R$ along the ideal of $S'$ is a formal power series ring $\overline{R}\llb t_1, \dotsc, t_n\rrb $, where $\overline{R}$ is a local ring of a smooth $k$-scheme and $n \le 2$ (because $R$ is~$2$-dimensional).

The completion map $R \rightarrow \overline{R}\llb t_1, \dotsc, t_n\rrb $ is faithfully flat and an isomorphism on residue fields at maximal ideals, so, for extending our $G$-torsor over $S\setminus z$ to a $G$-torsor over $S$, excision \cite{torsors-regular}*{Proposition 4.2.1} reduces us to when $R$ is $\overline{R}\llb t_1, …, t_n\rrb $, at the cost of losing the assumption that $R$ is a localization of a smooth $k$-algebra. We then apply the Popescu theorem $\overline{R}\llb t_1, \dotsc, t_{n - 1}\rrb $ and use a limit and algebraization argument based on \cite{Hitchin-torsors}*{Theorem 2.3.3 (c) (or Corollary 2.3.5 (a))} to reduce further to the case when $n = 1$, in other words, to when $R$ is $\overline{R}\llb t\rrb $ for an essentially smooth $k$-algebra $\overline{R}$ of dimension $1$ (because $R$ is of dimension $2$).

By spreading out, $\overline{R}$ is the local ring at a generic point of an irreducible divisor in a smooth, affine $k$-scheme $\overline{\mathcal{S}}$ of dimension $d > 0$. Moreover, we may assume that $k$ is imperfect (so infinite): otherwise, our smooth, wound, unipotent $k$-group $G$ is trivial, and so are its torsors. By the presentation lemma \cite{CTHK97}*{Theorem 3.1.1}, after shrinking around $\overline{R}$, our $\overline{\mathcal S}$ admits a smooth map of relative dimension $1$ to $\mathbb{A}^{d - 1}_k$ that makes the divisor in question finite over $\mathbb{A}^{d - 1}_k$. We may then replace $k$ by the function field $k(t_1, \dotsc, t_{d - 1})$  of $\mathbb{A}^{d - 1}_k$ to reduce to when $\overline{R}$ is a local ring at a closed point of a smooth, affine $k$-curve (the assumptions on $G$ are preserved because $k(t_1, \dotsc, t_{d - 1})$ is separable over $k$, see \S\ref{pp:fundamental-filtration}). By the uniqueness of the sought $G$-torsor extension and by Galois descent, we have the liberty of base changing to any finite Galois extension of $k$, so by a limit argument and further passage to $t$-adic completion, we may also assume that $k$ is separably closed.

By \cite{split-unramified}*{Lemma 6.3}, the completion of $\overline{R}$ is isomorphic to that of $\mathbb{A}^1_k$ at some closed point. Thus, by excision \cite{torsors-regular}*{Proposition 4.2.1} again, we are reduced to when $R = \overline{R}\llb t\rrb $, where $\overline{R}$ is the completion of $k[s]$ at some maximal ideal, explicitly, since $k$ is separably closed, at the ideal $(s^{p^n} - a) \subset k[s]$ for some $a \in k$. At this point, we are in the setting \Cref{lem:SH-trivial}~\ref{m:SHT-iii}, so we apply it to get that every $G$-torsor $E$ over $S \setminus z$ trivializes both over $\overline{R}\llp  t\rrp $ and over $K\llb t\rrb $ where $K$ is the fraction field of $\overline{R}$. It remains to note that, since $K$ is separable over $k$ and $G$ is pseudo-complete (see \Cref{prop:PP-groups}~\ref{m:PPG-ii}), any trivialization of $E$ over $\overline{R}\llp  t\rrp $ extends uniquely to a trivialization of $E$ over $K\llb t\rrb $, that is, to a trivialization over all of $S \setminus z$.

\m
The woundness of $G^0$ follows from the fact that $\mathrm{Ext}^1_{\ell}(\mathbb{G}_a, T) = 0$, see \cite{DG70}*{chapitre~III, section 6, no.~5, corollaire 5.2}. By \Cref{lem:B-Res}, we may assume that $\ell = k$, so that $G \simeq C/T$. By the purity for the Brauer group \cite{brauer-purity}*{Theorem 5.3} (in fact, already by \cite{Gro68c}*{th\'eor\`eme~6.1~b)}), we have $H^2(S \setminus z, T) \cong 0$, to the effect that every $G$-torsor over $S \setminus z$ lifts to a $C$-torsor over $S \setminus z$. To conclude, we will show that every $C$-torsor over $S \setminus z$ extends to a $C$-torsor over $S$ (equivalently, is trivial). 
    
    We consider the map $i_C\colon C\to \mathrm{Res}_{k^\prime/k}(\overline{C})$ of \S\ref{pp:pred-iG}, where $\overline{C}$ is a $k^\prime$-torus because it is both commutative and reductive. Since $\mathrm{Ker}(i_C)$ is commutative and pseudo-finite, \Cref{thm:purity-pseudo-finite} applies to it, and so allows us to replace $C$ by $i_C(C)$ to assume that $C$ is ultraminimal. Moreover, since $\mathrm{Res}_{k^\prime/k}(\overline{C})/C$ is affine (see \S\ref{pp:pred-iG}), as in the proof of \Cref{lem:SH-trivial}~\ref{m:SHT-i}, we may replace $C$ by $\mathrm{Res}_{k'/k}(\overline{C})$. Then, however, \Cref{lem:B-Res} reduces us to the case when $k^\prime = k$, so to when $C$ is a torus. The torus case holds, for instance, by \cite{brauer-purity}*{Theorem 6.1}.

\m
As in the proof of \ref{m:PPC-i}, we reduce to when our connected, smooth, wound, unipotent $k$-group $G$ is commutative and $p$-torsion. By \cite{Ros25}*{Theorem 1.8}, the assumption $[k : k^p] = p$ then ensures that $G \simeq C /T$ for some commutative pseudo-reductive $k$-group $C$ and a maximal $k$-torus $T \subset C$. Thus, \ref{m:PPC-ii} gives the claim. 

\m
As in the proof of \ref{m:PPC-i}, we may base change to any finite Galois subextension of $k^s/k$. By our assumption on $G_{k^s}$ and \Cref{lem:Ros}, after some such base change $G$ becomes a $\mathrm{Res}_{\ell/k}(C/T)$ as in \ref{m:PPC-ii}. Thus, \ref{m:PPC-ii} gives the claim. \qedhere
\end{enumerate}
\end{proof}

\begin{remark}
\Cref{thm:purity-pcomplete}~\ref{m:PPC-i} and \ref{m:PPC-iii} continue to hold if instead of assuming that $G$ is smooth, we assume that $G^0$ is wound unipotent and commutative. Indeed, we may first use \Cref{thm:purity-pseudo-finite} to replace $G$ by $G^0$ as in the proof above and then use \Cref{lem:embed-wound} to embed our commutative, connected, wound, unipotent $G$ as a $k$-subgroup of a commutative, connected, smooth, wound, unipotent $k$-group $\widetilde{G}$. Since $\widetilde{G}/G$ is affine (see \S\ref{pp:conv}), as in the proof of \Cref{lem:SH-trivial}~\ref{m:SHT-i}, we replace $G$ by $\widetilde{G}$ to reduce to when $G$ is, in addition, smooth, a case settled in \Cref{thm:purity-pcomplete}~\ref{m:PPC-i} and \ref{m:PPC-iii}.
\end{remark}

\begin{remark}
The proof of \Cref{thm:purity-pcomplete}~\ref{m:PPC-i} simplifies significantly if every $z \in Z$ of codimension $2$ in $S$ has $k_z/k$ separable: by the Cohen structure theorem \cite{Mat89}*{Theorems 28.3 and 30.6 (i)}, then the completion of $\mathscr{O}_{S,\, z}$ is $k$-isomorphic to $k_z\llb s, t\rrb $, so passage to power series rings becomes much more direct than in the proof above.
\end{remark}

\begin{remark}
\Cref{thm:purity-pcomplete}~\ref{m:PPC-ii} fails if we assume instead that $G^0 \simeq A/A^{\mathrm{ab}}$ is the largest unipotent quotient of a pseudo-abelian variety $A$ over $k$ as in \eqref{eqn:pseudo-abelian}: by \cite{Tot13}*{Corollary~7.3}, even $\mathbb{G}_{a,\, k}$ is of this form for suitable $k$, and the conclusion of \Cref{thm:purity-pcomplete} fails for it, see \eqref{eqn:A2}.
\end{remark}

\csub[Auslander--Buchsbaum extension for torsors under quasi-reductive groups]

A well-known extension result for torsors under reductive groups that is ultimately based on the Auslander--Buchsbaum formula to treat the key case of vector bundles says that for a regular scheme $S$ of dimension $2$, a closed subset $Z \subset S$ of codimension $2$ (so that $Z$ consists of isolated points of height $2$), and a reductive $S$-group $G$, pullback gives an equivalence of categories 
\be\label{eqn:reductive-extn}
\{G\text{-torsors over }S\} \xrightarrow{\sim} \{G\text{-torsors over }S \setminus Z\},
\ee
 see \cite{CTS79}*{Corollary 6.14}. We generalize this to quasi-reductive groups as follows and simultaneously extend \Cref{thm:purity-pproper,thm:purity-pcomplete} beyond pseudo-proper or pseudo-complete smooth $k$-groups $G$.

\begin{theorem}\label{thm:AB-extn}
Let $S$ be a geometrically regular $k$-scheme of dimension $\le 2$, let $Z \subset S$ be a closed subset of codimension $\ge 2$ \up{so that $Z$ consists of isolated points of height $2$}, and let $G$ be a smooth $k$-group scheme with $G^{\mathrm{sm},\, \mathrm{lin}}$ quasi-reductive. Suppose that either
\benumr
    \item\label{m:ABE-i} $G^{\mathrm{sm},\, \mathrm{lin}}$ is pseudo-reductive\uscolon or
    \item\label{m:ABE-ii} every $z \in Z$ lies in a geometrically regular $k$-subscheme $S_z \subset S$ of codimension $> 0$ \up{when $k_z/k$ is separable, we may take $S_z = z$}.
\end{enumerate}
For every gerbe $\mathscr{B}$ over $S$ isomorphic to $\mathbf{B}G$ \'etale locally on $S$, we have 
\[
\mathscr{B}(S) \xrightarrow{\sim} \mathscr{B}(S\setminus Z),
\]
in particular, for every $S$-group $\mathscr{G}$ isomorphic to $G$ \'etale locally on $S$, we have
\[
H^1(S, \mathscr{G}) \xrightarrow{\sim} H^1(S\setminus Z, \mathscr{G}) \quad \text{and, if }\mathscr{G}\text{ is commutative, also} \quad H^2(S, \mathscr{G}) \hookrightarrow H^2(S\setminus Z, \mathscr{G}).
\]
\end{theorem}

\bpf
By descent, we may work \'etale locally on $S$, so a spreading out argument allows us to assume that $S$ is strictly Henselian, local, of dimension $2$, and that $Z$ is its closed point. Then, since $G$ is $k$-smooth, every $G$-torsor over $S$ is trivial and we need to show that
\benuma
    \item\label{m:AB-1} $G(S) \xrightarrow{\sim} G(S \setminus Z)$; and
    \item\label{m:AB-2} every $G$-torsor over $S \setminus Z$ is trivial.
\end{enumerate}
The separatedness of $G$ gives the injectivity in \ref{m:AB-1}, see \cite{flat-purity}*{Lemma 7.2.7~(a)}, while the surjectivity is stable under extensions of smooth groups because each such group has no nontrivial torsors over $S$. To get the surjectivity for $G$ it then suffices to note that it holds for $G/G^{\mathrm{sm},\, \mathrm{lin}}$ by \Cref{thm:purity-pproper} and for $G^{\mathrm{sm},\, \mathrm{lin}}$ by \eqref{eqn:S-depth-2}. 

The remaining \ref{m:AB-2} is likewise stable under extensions of smooth groups and holds for $G/G^{\mathrm{sm},\, \mathrm{lin}}$ by \Cref{thm:purity-pproper}, so we may replace $G$ by $G^{\mathrm{sm},\, \mathrm{lin}}$ to assume that $G$ is quasi-reductive. Moreover, \Cref{thm:purity-pcomplete}~\ref{m:PPC-i} supplies \ref{m:AB-2} for the smooth, wound unipotent $k$-group $\mathscr{R}_{\mathrm{u},\, k}(G)$, so we may replace $G$ by $G/\mathscr{R}_{\mathrm{u},\, k}(G)$ to assume that $G$ is pseudo-reductive. For pseudo-reductive $G$, we first consider its central, pseudo-finite $k$-subgroup $\mathscr{C}_G$ discussed in \S\ref{pp:pseudo-minimal}. By \Cref{thm:purity-pseudo-finite}, every $\mathscr{C}_G$-torsor over $S \setminus Z$ extends to a  $\mathscr{C}_G$-torsor over $S$. This lets us replace $G$ by $G/\mathscr{C}_G$, in other words, we may assume that our pseudo-reductive $G$ is of minimal type. We then consider the comparison map of \S\ref{pp:pred-iG}:
$$
i_{G}\colon G \rightarrow \mathrm{Res}_{k^\prime/k}(\overline{G}).
$$
Since $G$ is of minimal type, its pseudo-finite $k$-subgroup $\mathrm{Ker}(i_{G})$ is commutative (see \S\ref{pp:pseudo-minimal}), so we may use \Cref{thm:purity-pseudo-finite} again to replace $G$ by $i_G(G)$, and hence reduce to when our pseudo-reductive $G$ is ultraminimal. As in the proof of \Cref{lem:SH-trivial}~\ref{m:SHT-i}, the affineness of $\mathrm{Res}_{k^\prime/k}(\overline{G})/G$ (see \S\ref{pp:pred-iG}) and \S\ref{eqn:S-depth-2} then allow us to replace $G$ by $\mathrm{Res}_{k^\prime/k}(\overline{G})$. Finally, \Cref{lem:B-Res} allows us to replace $\mathrm{Res}_{k^\prime/k}(\overline{G})$ by $\overline{G}$ to reduce to when $G$ is reductive. In the reductive case, \ref{m:AB-2} follows from \eqref{eqn:reductive-extn}.
\epf

\begin{remark}\label{rem:qred-optimal}
To show that the quasi-reductivity assumption is optimal in \Cref{thm:AB-extn}, we suppose that $\mathscr{R}_{\mathrm{us},\, k}(G) \neq 1$ (with $G$ still smooth), we let $S$ be the strict Henselization of the origin of $\mathbb{A}^2_k$ with $Z \subset S$ the closed point, and we will show that $G$ has a nontrivial torsor over $S\setminus Z$, contrary to what the conclusion of \Cref{thm:AB-extn} would predict. By the proof of \Cref{thm:AB-extn}, nontrivial $\mathbb{G}_a$-torsors over $\mathbb{A}^2_k \setminus Z$ restrict to nontrivial $\mathbb{G}_a$-torsors over $S \setminus Z$. Such restrictions induce nontrivial $\mathscr{R}_{\mathrm{us},\, k}(G)$-torsors $E$ over $S \setminus Z$ (see the next sentence or the proof of \Cref{lem:SH-trivial}~\ref{m:SHT-i}). If such an $E$ induced the trivial $G$-torsor, then $E$ would be the preimage of an $(S \setminus Z)$-point of $G/\mathscr{R}_{\mathrm{us},\, k}(G)$. By \Cref{thm:AB-extn}, this $(S\setminus Z)$-point extends to an $S$-point, so $E$ would extend to an $\mathscr{R}_{\mathrm{us},\, k}(G)$-torsor over $S$, a contradiction. The importance of quasi-reductivity, relatedly, the failure of \Cref{thm:AB-extn} for $\mathbb{G}_a$ invalidates several more direct ways to attack \Cref{thm:AB-extn}.
\end{remark}

\begin{remark}\label{rem:playful}
Although there are other ways to see this, \Cref{rem:qred-optimal} gives a playful proof of the fact that the quotient $G/\mathscr{Z}$ of a quasi-reductive $k$-group $G$ by a central $k$-subtorus $\mathscr{Z}$ is still quasi-reductive. Indeed, with $Z \subset S$ as in \Cref{rem:qred-optimal}, it suffices to note that, by \cite{brauer-purity}*{Theorem~1.3} (in fact, already by \cite{Gro68c}*{th\'eor\`eme 6.1 b)}), every $G/\mathscr{Z}$-torsor over $S \setminus Z$ lifts to a $G$-torsor over $S \setminus Z$, and so, by \Cref{thm:AB-extn}, it extends to a $G/\mathscr{Z}$-torsor over $S$, and hence is trivial.
\end{remark}


\section{Classification of $G$-torsors over $\mathbb{P}^1_k$}

The endpoint of the geometric approach to the Grothendieck--Serre conjecture is the study of torsors over the relative $\mathbb{P}^1$. As we explain in \S\ref{sec:torsors-P1}, the extension results of Chapter \ref{chap:purity} lead to a quick classification of $G$-torsors over $\mathbb{P}^1_k$ for quasi-reductive $G$. For this, we begin in \S\ref{sec:P1-unipotent} with the auxiliary and more direct unipotent case.

\csub[Torsors over $\mathbb{P}^n_S$ under unipotent groups]
\label{sec:P1-unipotent}

Torsors over a relative $\mathbb{P}^n_S$ under unipotent $S$-groups descend to $S$ as in \Cref{prop:unipotent-P1} below, whose proof uses the following auxiliary lemma.

\begin{lemma}\label{lem:P1-alg}
For every $n \ge 0$, base change induces an equivalence of categories
\[
\{ \text{finite }k\text{-schemes}\} \xrightarrow{\sim} \{\text{finite, flat }\mathbb{P}^n_k\text{-schemes } F \text{ such that }\mathscr{O}_F \simeq \mathscr{O}_{\mathbb{P}^n_k}^{\oplus d} \text{ as $\mathscr O_{\mathbb P^n_k}$-modules} \},
\]
with an inverse functor given by $F \mapsto s^*(F)$ for any fixed $s \in \mathbb{P}^n_k(k)$.
\end{lemma}

\begin{proof}
Since $\mathrm{Hom}_{\mathscr{O}_{\mathbb{P}^n_k}}(\mathscr{O}_{\mathbb{P}^n_k},\mathscr{O}_{\mathbb{P}^n_k}) \cong k$, the triviality of $\mathscr{O}_F$  as a vector bundle  implies that commutative $\mathscr{O}_{\mathbb{P}^n_k}$-algebra structures on $\mathscr{O}_F$ correspond to commutative $k$-algebra structures on $s^*(\mathscr{O}_F)$. 
\end{proof}

\begin{proposition}\label{prop:unipotent-P1}
Let $S$ be an algebraic space and let $\mathscr{G}$ be a flat, locally finitely presented, quasi-separated $S$-group algebraic space each of whose geometric $S$-fibers $\mathscr{G}_{\overline{K}}$ has a unipotent linear part $(\mathscr{G}_{\overline{K}})^{\mathrm{lin}}$ \up{see \uS\uref{pp:k-groups-lft} and \uS\uref{pp:fundamental-filtration}~\uref{m:Glin}}, and suppose that $\mathscr{G}$ is an extension of an fpqc locally constant $S$-group $\mathscr{C}$ by a finitely presented $S$-group algebraic space. For every $S$-gerbe $\mathscr{B}$ for the fppf topology isomorphic to $\mathbf{B}\mathscr{G}$ fppf locally on $S$, we have $\mathscr{B}(S) \xrightarrow{\sim} \mathscr{B}(\mathbb{P}^n_S)$ for every $n \ge 1$, in particular,
\[
\mathscr{G}(S) \xrightarrow{\sim} \mathscr{G}(\mathbb{P}^n_S) \quad \text{and} \quad H^1(S, \mathscr{G}) \xrightarrow{\sim} H^1(\mathbb{P}^n_S, \mathscr{G}).
\]
\end{proposition}

\begin{proof}
By working fppf locally on $S$, we may assume that $\mathscr{B} = \mathbf{B}\mathscr{G}$. We then first claim that
\begin{equation}\label{eqn:ind-quasi-affine}
\mathscr{G}(S) \xrightarrow{\sim} \mathscr{G}(\mathbb{P}^n_S).
\end{equation}

For this, since $\mathbb{P}^n_S$ has an $S$-point, we may focus on the surjectivity and, by a descent and limit argument, may assume that $S$ is affine and even local. The rigidity lemma \cite{MFK94}*{Proposition 6.1} then reduces us to $S$ being the spectrum of a field $K$, and fpqc descent allows us to assume that $K$ is algebraically closed. To then check that every $\mathbb{P}^n_K$-point of $\mathscr{G}$ comes from a $K$-point, by \eqref{eqn:reduce-to-gred} and translation, we may assume that $\mathscr{G}$ is smooth and connected. However, every $\mathbb{P}^n_K$-point of a smooth, connected $K$-group comes from a $K$-point: for abelian varieties this is \cite{Mil86a}*{Corollary~3.8}, for affine groups this results from $\Gamma(\mathbb{P}^n_K, \mathscr{O}) \cong K$, and the general case follows from \S\ref{pp:fundamental-filtration}~\ref{m:Glin}.

By applying \eqref{eqn:ind-quasi-affine} fpqc locally on $S$, we find that
$$
\underline{\mathrm{Isom}}_{\mathscr G}(E, E^\prime)(S) \xrightarrow{\sim} \underline{\mathrm{Isom}}_{\mathscr G}(E, E^\prime)(\mathbb{P}^n_S)
$$
for all $\mathscr{G}$-torsors $E$ and $E^\prime$ over $S$, in other words, the functor $\mathscr{B}(S) \hookrightarrow \mathscr{B}(\mathbb{P}^n_S)$ is fully faithful. To argue that it is also essentially surjective, by fpqc descent again, it then suffices to show that every $\mathscr{G}$-torsor $E$ over $\mathbb{P}^n_S$ trivializes fpqc locally on $S$. For this, letting $s$ be an $S$-point of $\mathbb{P}^n_S$ and working fpqc locally again, we may assume that $s^*(E)$ is trivial and $\mathscr{C}$ is constant. 

Every $\mathscr C$-torsor over $\mathbb{P}^n_S$ trivializes even \'etale locally on $S$: for this, by a limit argument and the proper base change theorem \cite{SGA4III}*{expos\'e XII, corollaire 5.5 (ii)}, we may assume that $S$ is the spectrum of a separably closed field $k$, then note that, by \cite{SGA3II}*{expos\'e X, corollaire 5.14}, every connected component of a $\mathscr C$-torsor over $\mathbb{P}^n_k$ is finite \'etale over $\mathbb{P}^n_k$, and, finally, recall that this \'etale cover splits by \cite{SGA1new}*{expos\'e XI, proposition 1.1} (with \cite{SP}*{Theorem \href{https://stacks.math.columbia.edu/tag/0BTY}{0BTY}} to pass to $\overline{k}$). By applying this to the $\mathscr C$-torsor induced by $E$, we reduce the sought fpqc local triviality of $E$ to the case when $\mathscr C$ is trivial, so that $\mathscr G$ is finitely presented over $S$, with the triviality of $s^*(E)$ again arranged by a preliminary fpqc base change on $S$.

Once $\mathscr{G}$ is finitely presented over $S$, we may first do a descent and limit argument to reduce to Noetherian $S$ and then, by a further limit argument and fpqc base change that uses \cite{EGAIII1}*{chapitre 0, proposition 10.3.1}, even to $S$ being the spectrum of a complete Noetherian local ring with an algebraically closed residue field. We seek a trivialization of $E$ over $\mathbb{P}^n_S$ whose restriction to $s$ recovers a fixed trivialization of $s^*(E)$, and, by \eqref{eqn:ind-quasi-affine}, there is at most one such. This way of choosing a preferred trivialization uniquely and compatibly with base change in $S$ allows us to apply the continuity formula \cite{BHL17}*{Corollary 1.6} (with \cite{torsors-regular}*{Section 1.2.3} to ensure that $E$ is an algebraic space) to reduce to when $S$ is Artinian local with an algebraically closed residue field $K$.

At this point the vanishing of the cohomology $H^1(\mathbb{P}^n_k, \mathscr{O}) \cong H^2(\mathbb{P}^n_k, \mathscr{O}) \cong 0$ and the deformation theory of $\mathscr G$-torsors, more precisely, \cite{Ill72}*{th\'eor\`eme 2.4.4, page 209}, imply that $E$ is the unique deformation of $E|_{\mathbb{P}^n_k}$ to a $\mathscr G$-torsor over $\mathbb{P}^n_S$ (see the second half of the proof of \cite{totally-isotropic}*{Proposition~3.1~(b)} for more details). In particular, to prove that $E$ is trivial, we may base change to the special fiber and reduce to the case when $S = \mathrm{Spec}(K)$.

Recalling that $\mathscr C$ is trivial and $\mathscr G$ is of finite type, consider the finite $K$-group $F := \mathscr G/(\mathscr G^{0})^{\mathrm{red}}$ and the $F$-torsor $\overline{E}$ over $\mathbb{P}^n_K$ induced by $E$. Since $F/F^0$ is \'etale (see \S\ref{pp:fundamental-filtration}~\ref{m:G0}), the $F/F^0$-torsor $\overline{E}/F^0$ is a finite \'etale cover of $\mathbb{P}^n_K$. Thus, $\overline{E}/F^0$ is trivial by \cite{SGA1new}*{expos\'e XI, proposition 1.1}, so we may assume that $F = F^0$, so that $F$ is infinitesimal. To show that $\overline{E}$ descends to an $F$-torsor over $K$, and hence is trivial, by \Cref{lem:P1-alg}, it suffices to show that $\mathscr{O}_{\overline{E}}$ is a trivial vector bundle over $\mathbb{P}^n_K$. For this, by \cite{OSS80}*{Theorem 3.2.1} (which is stated over the complex numbers but whose proof works over any base field), it suffices to show the same after restricting to any line $L \subset \mathbb{P}^n_K$, so we may assume that $n = 1$. However, the vector bundle $\mathscr{O}_{\overline{E}}$ trivializes over the finite, flat cover $\overline{E} \rightarrow \mathbb P^1_K$, and hence also over $(\overline{E})^{\mathrm{red}}$. Since $F$ is infinitesimal, the function field of the integral $K$-curve $(\overline{E})^{\mathrm{red}}$ is purely inseparable over $K(t)$, so it is contained in some $K(t^{1/p^\ell})$ with $\ell \ge 0$. By passing to normalizations, we conclude that $\mathscr{O}_{\overline{E}}$ trivializes over some finite, flat cover $\mathbb{P}^1_K \rightarrow \mathbb{P}^1_K$. By then decomposing $\mathscr{O}_{\overline{E}}$ into a sum of line bundles $\mathscr{O}(n)$, we get that $\mathscr{O}_{\overline{E}}$ is trivial, as desired. 

Since $\overline{E}$ is trivial, $E$ reduces to a  $(\mathscr G^{0})^{\mathrm{red}}$-torsor over $\mathbb P^n_K$ and we may assume that $\mathscr G$ is smooth and connected. Our assumptions and \S\ref{pp:fundamental-filtration}~\ref{m:Glin} now imply that $\mathscr G$ is an extension of an abelian variety $A$ by a unipotent group. By \eqref{eqn:abelian-scheme-extend},
$$
H^1(\mathbb{P}^n_K, A) \hookrightarrow H^1(k(\mathbb{P}^n_K), A),
$$
so every $A$-torsor over $\mathbb{P}^n_K$ has finite order, in other words, reduces to an $A[m]$-torsor for some $m\geq 0$. Since the $K$-group $A[m]$ is finite, by the argument above about $F$, it has no nontrivial torsors over $\mathbb{P}^n_K$, and hence neither does $A$. Consequently, we may assume that $\mathscr{G}$ is unipotent, and, by passing to subquotients, that it is either finite or $\mathbb{G}_{a,\, K}$ (see \S\ref{pp:conv}). In the finite case, the argument above about $F$ suffices, while $\mathbb{G}_{a,\, K}$ has no nontrivial torsors over $\mathbb P^n_K$ because $H^1(\mathbb{P}^n_K, \mathscr{O}) \cong 0$. All in all, our $\mathscr{G}$-torsor $E$ over $\mathbb{P}^n_K$ is trivial, so it descends to $K$. 
\end{proof}

\begin{remark}\label{rem:unip-higher}
As for higher cohomology, for every commutative, unipotent $k$-group $U$ and every $k$-algebraic space $S$, we have $H^i(S, U) \xrightarrow{\sim} H^i(\mathbb{P}^n_S, U)$ for all $i$. Indeed, to see that $U \xrightarrow{\sim} R\pi_*(U_{\mathbb{P}^n_S})$ where $\pi \colon \mathbb{P}^n_S \rightarrow S$ is the structure map, we may assume that $k$ is algebraically closed, handle $U = \mathbb{G}_{a,\, k}$ by \cite{EGAIII1}*{proposition 2.1.15}, and then handle a general $U$ by d\'evissage, \Cref{prop:unipotent-P1}, and \cite{SGA3II}*{expos\'e XVII, corollaire 1.7}.
\end{remark}

Beyond unipotent groups, \Cref{prop:unipotent-P1} helps reduce the structure group of $G$-torsors over $\mathbb{P}^n_S$ to the (unirational) $k$-subgroup $G^{\mathrm{tor}} \le G$ generated by the $k$-tori of $G$ (see \S\ref{pp:fundamental-filtration}~\ref{m:Gtor}) as follows. 

\begin{corollary}\label{cor:reduce-to-Gt}
Let $G$ be a smooth $k$-group scheme, let $S$ be a $k$-scheme, and let $E$ be a $G$-torsor over $\mathbb{P}^n_S$ with $n \ge 1$. Pullback along an $s \in \mathbb{P}^n_S(S)$ gives an equivalence
\[
\{\text{reductions of } E \text{ to a } G^{\mathrm{tor}}\text{-torsor}\} \xrightarrow{\sim} \{\text{reductions of } s^*(E) \text{ to a } G^{\mathrm{tor}}\text{-torsor}\}.
\]

In particular, if $s^*(E)$ is trivial, then $E$ reduces to a $G^{\mathrm{tor}}$-torsor over $\mathbb{P}^n_S$ whose $s$-pullback is trivial \up{resp.,~then $E$ is even trivial if $G^{\mathrm{tor}} = 1$, that is, if $G^{\mathrm{sm},\, \mathrm{lin}}$ is unipotent}.
\end{corollary}

\begin{proof}
The category of reductions of $E$ (resp.,~of $s^*(E)$) to a $G^{\mathrm{tor}}$-torsor is equivalent to the set $(E/G^{\mathrm{tor}})(\mathbb{P}^n_S)$ (resp.,~$(s^*(E/G^{\mathrm{tor}}))(S)$). Moreover, $E/G^{\mathrm{tor}}$ is a $(G/G^{\mathrm{tor}})$-torsor over $\mathbb{P}^n_S$, so \Cref{prop:unipotent-P1} gives the claim once we argue that $(G/G^{\mathrm{tor}})^{\mathrm{lin}}$ is unipotent. The latter follows from \S\ref{pp:fundamental-filtration}~\ref{m:Glin}, \ref{m:Gsmlin}, and \ref{m:Gtor}. 
\end{proof}

The following consequence of \Cref{cor:reduce-to-Gt} extends \cite{totally-isotropic}*{Theorem 4.2} beyond reductive groups.

\begin{corollary}
For a smooth $k$-group scheme $G$ with $G^{\mathrm{sm},\, \mathrm{lin}}$ unipotent and a $k$-algebra $A$, no nontrivial $G$-torsor $E$ over $\mathbb{A}^1_A$ trivializes over the punctured formal neighborhood $A\llp t^{-1}\rrp$ of infinity.
\end{corollary}

\begin{proof}
For convenience, we set $s \colonequals t^{-1}$. Since $G$ is smooth, by \cite{Hitchin-torsors}*{Corollary 2.1.22 (b) (with Example 2.1.18)}, our $E$ trivializes already over the punctured Henselization $A\{s\}[\frac{1}{s}]$ at infinity. In particular, by patching \cite{MB96}*{Theorem 5.5}, it extends to a $G$-torsor over $\mathbb{P}^1_A$ that is trivial at infinity. \Cref{cor:reduce-to-Gt} then implies that $E$ is trivial. 
\end{proof}

\csub[Torsors over $\mathbb{P}^1_k$ under quasi-reductive groups]
\label{sec:torsors-P1}

We apply the Auslander--Buchsbaum extension \Cref{thm:AB-extn} for torsors under quasi-reductive groups to obtain a classification of torsors over $\mathbb{P}^1_k$ in \Cref{thm:torsors-P1}. For this, we view the projective space as open $[\mathbb{A}^n_k/\mathbb{G}_m]$ as follows and then classify torsors over this stack in \Cref{lem:maps-full}.

\beg\label{eg:Ank-CC}
Letting $\mathbb{G}_m$ act on $\mathbb{A}^{n + 1}_k$ by scaling the coordinates, the open complement of $\mathbf{B}\mathbb{G}_m \cong [\{0\}/\mathbb{G}_m]$ in the quotient $[\mathbb{A}^{n + 1}_k/\mathbb{G}_m]$ is the projective space $\mathbb{P}^n_k \cong [(\mathbb{A}^{n + 1}_k \setminus \{0\})/\mathbb{G}_m]$. The structure map $\mathbb{P}^n_k \rightarrow \mathbf{B}\mathbb{G}_m$ classifies the line bundle $\mathscr{O}(-1)$, that is, as we now argue, the square
\[
\xymatrix{
\mathbb{A}_k^{n+1} \setminus \{0\} \ar[r] \ar[d] & \mathrm{Spec}(k) \ar[d] \\
\mathbb{P}^n_k \ar[r]^-{\varphi_{[\mathcal{O}(-1)]}} & \mathbf{B}\mathbb{G}_m
}
\]
is a $\mathbb{G}_m$-equivariant and Cartesian. To see this, first recall from \cite{SP}*{Lemma \href{https://stacks.math.columbia.edu/tag/01NE}{01NE}} (or \cite{EGAII}*{th\'eor\`eme 4.2.4}) that the fiber product of the outer bottom part of the square represents the functor that sends a $k$-scheme $S$ to the set of isomorphism classes of line bundles $\mathscr{L}$ on $S$ equipped both with $s_0,\dots,s_n\in H^0(S,\mathscr{L})$ that have no common zeros and with a $\sigma \colon \mathscr{O}_S\xrightarrow{\sim}\mathscr L^\vee$: explicitly, $\mathscr{L}$ is the pullback of $\mathscr{O}(1)$ along the map $S\to \mathbb{P}^n_k$ given by $x \mapsto [s_0(x):\dots : s_n(x)]$. On the other hand, $\mathbb{A}^{n+1}_k\setminus\{0\}$ represents the functor that sends $S$ to the set of $(n + 1)$-tuples $t_0, \dotsc, t_n \in H^0(S, \mathscr{O}_S)$ with no common zeros. The two functors are identified via $t_i = \sigma^\vee \circ s_i$, and this is $\mathbb{G}_m$-equivariant because $\mathbb{G}_m$ acts on both sides by scaling via a character of weight $-1$ (so $t_i \mapsto \lambda^{-1}t_i$, etc.).
\eeg

Torsors over $[\mathbb{A}^n_k/\mathbb{G}_m]$ under smooth, affine $k$-groups may be classified as follows.

\begin{lemma}\label{lem:maps-full}
In Example~\uref{eg:Ank-CC}, for a smooth $k$-group $G$, we have
\begin{equation}\label{eqn:M}
H^1(\mathbf{B}\mathbb{G}_m, G) \xrightarrow{\sim} H^1([\mathbb{A}^n_k/\mathbb{G}_m], G)
\end{equation}
via pullback along the structure map $[\mathbb{A}^n_k/\mathbb{G}_m] \rightarrow \mathbf{B}\mathbb{G}_m$, and if $G$ is also affine, then restricting along the origin $\mathbf{B}\mathbb{G}_m \hookrightarrow [\mathbb{A}^n_k/\mathbb{G}_m]$ gives a full, essentially surjective functor
\begin{equation}\label{eqn:Gpb}
\{G\text{-torsors over }[\mathbb{A}^n_k/\mathbb{G}_m]\} \to \{G\text{-torsors over }\mathbf{B}\mathbb{G}_m\}.
\end{equation}
\end{lemma}

\bpf
The composition $\mathbf B\mathbb{G}_m\hookrightarrow [\mathbb{A}^n_k/\mathbb{G}_m]\rightarrow \mathbf B \mathbb G_m$ is the identity, so \eqref{eqn:M} is injective, and \eqref{eqn:Gpb} is essentially surjective. To show that \eqref{eqn:M} is also surjective, we need to argue that every $G$-torsor $E$ over $[\mathbb A^n_k/\mathbb G_m]$ descends to a $G$-torsor over $\mathbf B\mathbb G_m$. For this, we may replace $G$ by an inner form to assume that the restriction of $E$ to some $k$-point of $\mathbb{P}^{n - 1}_k \cong [(\mathbb A^n_k \setminus \{0\})/\mathbb G_m] \subset \mathbf [\mathbb A^n_k/\mathbb G_m]$ is trivial, see \cite{Gir71}*{chapitre III, remarque 2.6.3}. \Cref{cor:reduce-to-Gt} then ensures that $E|_{\mathbb{P}^{n - 1}_k}$ reduces a $G^{\mathrm{tor}}$-torsor over $\mathbb{P}^{n - 1}_k$, in particular, to a $G^{\mathrm{sm},\, \mathrm{lin}}$-torsor over $\mathbb{P}^{n - 1}_k$. At this point, by \eqref{eqn:E-Gsmlin}, our $E$ itself reduces to a $G^{\mathrm{sm},\, \mathrm{lin}}$-torsor over $[\mathbb A^n_k/\mathbb G_m]$, so we may assume that $G$ is affine. For affine $G$, however, the surjectivity in \eqref{eqn:M} follows from the fullness of \eqref{eqn:Gpb}. 

The claim about \eqref{eqn:Gpb} is a special case of \cite{Wed24}*{Corollary 2.9}. Alternatively, the remaining fullness of the functor \eqref{eqn:Gpb} results from the coherent completeness of $[\mathbb{A}^{n + 1}_k/\mathbb{G}_m]$ along its stacky origin $\mathbf B\mathbb{G}_m$ (see \cite{AHR20}*{Theorem 1.3 or the proof of Lemma 5.18}) and a deformation-theoretic argument that is analogous to that of \cite{AHHL21}*{Lemma 3.3}.\footnote{\label{foot:Alper}In \emph{op.~cit.},~it is useful to first replace $G$ by $G \times_k \mathbb G_m$ and modify $f$ and $f^\prime$ by adding the structure map to the resulting new factor $\mathbf B \mathbb G_m$ in their target. This artificial maneuver does not affect the deformability of $2$-isomorphisms relevant for the argument, but it makes the maps $f$ and $f^\prime$ representable, so that the deformation-theoretic input of \cite{Ols06a} cited in the argument applies (\emph{op.~cit.}~was written for representable maps, see \cite{Ols06a}*{Remark 1.8}), compare also with \cite{Alp10}*{Proposition 5.1 and its proof}.}
\epf

The following final input to \Cref{thm:torsors-P1} is widely known for reductive $k$-groups.

\begin{lemma}\label{lem:SGconj}
For a $k$-group scheme $G$ locally of finite type, the maximal split $k$-tori $S \le G$ are pairwise $G(k)$-conjugate and, for any such $S$, we have
\[
\mathrm{Hom}_{k\text{-\upshape{gp}}}(\mathbb{G}_m, S)/N_G(S)(k) \xrightarrow{\sim} \mathrm{Hom}_{k\text{-\upshape{gp}}}(\mathbb{G}_m, G)/G(k)
\]
\end{lemma}

\bpf
The $G(k)$-conjugacy of maximal split $k$-tori is \cite{CGP15}*{Proposition C.4.5 (1)}. Thus, since every $k$-homomorphism $\lambda \colon \mathbb{G}_m \rightarrow G$ factors through some maximal split $k$-torus of $G$, the displayed map is surjective. For its injectivity, fix a maximal split $k$-torus $S \le G$ and consider $k$-homomorphisms $\lambda, \lambda^\prime\colon \mathbb{G}_m \rightarrow S$ with $\lambda(-) = g\lambda^\prime(-)g^{-1}$ for some $g \in G(k)$. Both $S$ and $gSg^{-1}$ are maximal $k$-split tori of $G$ through which $\lambda$ factors, so both lie in the largest connected, smooth, affine $k$-subgroup $G^{\mathrm{sm},\, \mathrm{lin}} \le G$, and then even in the centralizer $Z_{G^{\mathrm{sm},\, \mathrm{lin}}}(\lambda)$ of $\lambda$. However, \cite{SGA3II}*{expos\'e XI, corollaire 5.3} ensures that $Z_{G^{\mathrm{sm},\, \mathrm{lin}}}(\lambda)$ is a closed $k$-subgroup of $G$. The conjugacy of maximal split $k$-tori applies to it and gives an $h \in (Z_{G^{\mathrm{sm},\, \mathrm{lin}}}(\lambda))(k)$ with $S = hgSg^{-1}h^{-1}$. Now $hg \in N_G(S)(k)$ and $hg\lambda^\prime(-)(hg)^{-1} = h\lambda(-)h^{-1} = \lambda(-)$, so that $\lambda$ and $\lambda^\prime$ are $N_G(S)(k)$-conjugate. 
\epf

\begin{theorem}\label{thm:torsors-P1}
Let $G$ be a $k$-group scheme locally of finite type such that every $\overline{k}$-torus of $G_{\overline{k}}$ lies in $(G^{\mathrm{gred}})_{\overline{k}}$ \up{see Remark \uref{rem:condition-star}} and let $E$ be a $G$-torsor over $\mathbb{P}^1_k$. No nontrivial $G$-torsor over $k[t]$ \up{resp.,~over $k[t]_{1 + tk[t]}$\uscolon resp.,~over $k\{t\}$\uscolon resp.,~over $k\llb t \rrb$} trivializes after inverting $t$, and the following conditions are equivalent:
\benumr
    \item\label{m:TP-i} $E|_{\mathbb{A}^1_k}$ is trivial\uscolon
    \item\label{m:TP-ii} $E$ is Zariski locally trivial\uscolon
    \item\label{m:TP-iii} $E$ is generically trivial\uscolon
\end{enumerate}
and if $G$ is smooth, then they are also equivalent to
\benumr\setcounter{enumi}{3}
    \item\label{m:TP-iv} $s^*(E)$ is trivial for some $s \in \mathbb{P}^1_k(k)$\uscolon
    \item\label{m:TP-v} $s^*(E)$ is trivial for every $s \in \mathbb{P}^1_k(k)$\uscolon
    \item\label{m:TP-vi} $E|_{k((t))}$ is trivial, where $t$ is the standard coordinate of $\mathbb{A}^1_k$.
\end{enumerate}

If $G$ is smooth, then, for every $s \in \mathbb{P}^1_k(k)$, the sequence of pointed sets
\begin{equation}\label{eqn:H1-P1}
\{*\} \to H^1_{\mathrm{Zar}}(\mathbb{P}^1_k, G) \to H^1(\mathbb{P}^1_k, G)\xrightarrow{s^*} H^1(k, G)\to \{*\}
\end{equation}
is exact. If $G$ is smooth with $G^{\mathrm{sm},\, \mathrm{lin}}$ quasi-reductive \up{see \uS\uref{pp:fundamental-filtration}~\uref{m:Gsmlin}}, then
\begin{equation}\label{eqn:H1et}
H^1(\mathbb{P}^1_k, G) \cong H^1(\mathbf{B}\mathbb{G}_m, G),
\end{equation}
letting $S \le G$ be a maximal $k$-split torus, also
\begin{equation}\label{eqn:H1Zar}
H^1_{\Zar}(\mathbb{P}^1_k,G) \cong \mathrm{Hom}_{k\text{-\upshape{gp}}}(\mathbb{G}_m, G)/G(k) \overset{\ref{lem:SGconj}}{\cong} \mathrm{Hom}_{k\text{-\upshape{gp}}}(\mathbb{G}_m, S)/N_G(S)(k),
\end{equation}
and, in addition, the conditions \ref{m:TP-i}--\ref{m:TP-vi} are also equivalent to
\benumr\setcounter{enumi}{6}
    \item\label{m:TP-vii} $E$ is the inflation of the $\mathbb{G}_m$-torsor $\mathscr{O}(1)$ along some $\lambda\colon \mathbb{G}_{m,\, k} \rightarrow G$.
\end{enumerate}
\end{theorem}

\bpf
By patching \cite{MB96}*{Theorem 5.5}, a $G$-torsor $\mathcal{E}$ over $k[t]$ (resp.,~over $k[t]_{1 + tk[t]}$; resp.,~over $k\{t\}$; resp.,~over $k\llb t\rrb$) that trivializes after inverting $t$ extends to a $G$-torsor over $\mathbb{P}^1_k$ that trivializes over $\mathbb{P}^1_k \setminus \{t = 0\}$ (\emph{loc.~cit.}~applies since $k\llb t \rrb$ is a filtered direct limit of flat, finitely presented $k[t]$-algebras). Thus, the equivalence of \ref{m:TP-i}--\ref{m:TP-iii} implies that $\mathcal{E}$ is trivial. If at least one of \ref{m:TP-i}--\ref{m:TP-iii} holds, then, by \Cref{thm:pseudo-extend}~\ref{m:PE-i}, our $E$ reduces to a generically trivial $G^{\mathrm{sm},\, \mathrm{lin}}$-torsor. In particular, we may assume throughout that $G$ is smooth. 

The patching argument also shows that \ref{m:TP-vi} implies \ref{m:TP-iv}, so since \ref{m:TP-vi} follows from \ref{m:TP-i}, we may discard \ref{m:TP-vi} altogether. Similarly, we may discard \ref{m:TP-v} because it implies \ref{m:TP-iv} and follows from \ref{m:TP-ii}. Moreover, the equivalence of \ref{m:TP-ii} and \ref{m:TP-v} gives \eqref{eqn:H1-P1}.  

If either \ref{m:TP-i}--\ref{m:TP-iii} holds for our smooth $G$, then we saw that, by \Cref{thm:pseudo-extend}~\ref{m:PE-b}~\ref{m:PE-i}, we may assume that $G$ is connected, smooth, and affine. If \ref{m:TP-iv} holds, then \Cref{cor:reduce-to-Gt} allows us to reduce $E$ to a $G^{\mathrm{tor}}$-torsor $E^{\mathrm{tor}}$ with $s^*(E^{\mathrm{tor}})$ trivial, so we may again assume that $G$ is connected, smooth, and affine. Moreover, since $\mathbb{G}_a$ has no nontrivial torsors over affine schemes, \ref{m:TP-i}--\ref{m:TP-iv} are all insensitive to replacing $G$ by $G/\mathscr{R}_{\mathrm{us},\, k}(G)$. All in all, to argue that each of \ref{m:TP-i}--\ref{m:TP-iv} implies the others, we may assume that $G$ is quasi-reductive.

In the remaining case when $G$ is smooth with $G^{\mathrm{sm},\, \mathrm{lin}}$ quasi-reductive, we view $\mathbb{P}^1_k\cong [(\mathbb{A}^2_k\setminus\{0\})/\mathbb{G}_m]$ as the open complement of $\mathbf B\mathbb{G}_m \cong [\{0\}/\mathbb{G}_m]$ in $[\mathbb{A}^2_k/\mathbb{G}_m]$, see \Cref{eg:Ank-CC}. By the Auslander--Buchsbaum extension \Cref{thm:AB-extn}~\ref{m:ABE-ii} for $G$-torsors (applied after pullback along $\mathbb{A}^2_k \rightarrow [\mathbb{A}^2_k/\mathbb{G}_m]$) and \eqref{eqn:M}, we have
$$
H^1(\mathbb{P}^1_k, G) \cong H^1([\mathbb A^2_k/\mathbb G_m], G) \cong H^1(\mathbf B\mathbb G_m, G).
$$
Here the Zariski locally trivial $G$-torsors on $\mathbb{P}^1_k$ correspond to those $G$-torsors over $\mathbf B \mathbb G_m$ that trivialize over the cover $\mathrm{Spec}(k) \rightarrow \mathbf B \mathbb G_m$: indeed, Zariski local triviality implies triviality at every $k$-point of $\mathbb P^1_k$, whereas the $G$-torsors over $\mathbf B\mathbb G_m$ that trivialize over $\mathrm{Spec}(k)$ are the inflations of the tautological $\mathbb G_m$-torsor along some $k$-homomorphism $\lambda \colon \mathbb G_m \rightarrow G$, where $G$ is the automorphism group of a trivial $G$-torsor $E_0$ over $\mathrm{Spec}(k)$ and $\lambda$ is unique up to changing a trivialization of $E_0$, concretely, up to $G(k)$-conjugation. The claimed \eqref{eqn:H1Zar} follows, and we also get that either of \ref{m:TP-i}, \ref{m:TP-ii}, \ref{m:TP-iv}, or \ref{m:TP-vii} implies all of \ref{m:TP-i}--\ref{m:TP-vii}. By spreading out, so does \ref{m:TP-iii} if $k$ is infinite. In the remaining case when $k$ is finite and \ref{m:TP-iii} holds, we already saw that to argue \ref{m:TP-i}--\ref{m:TP-vii}, we may assume that $G$ is quasi-reductive. Then Lang's theorem gives \ref{m:TP-iv} (see \cite{Ser02}*{Chapter III, Section 2.3, Theorem~1${}^\prime$}), and hence, by what we have already argued, \ref{m:TP-i}--\ref{m:TP-vii} all hold.
\epf

\begin{remark}\label{rem:classification-fail}
We now show that \eqref{eqn:H1Zar} fails without the quasi-reductivity assumption. Torsors under $G \colonequals \mathbb{G}_a \rtimes \mathbb{G}_m$, where $\mathbb{G}_m$ acts on $\mathbb{G}_a$ via a character of weight $1$, are all Zariski locally trivial. By \cite{Gir71}*{chapitre III, remarque 2.6.3, proposition 3.3.1 (i)}, the set of isomorphism classes of those $G$-torsors over $\mathbb{P}^1_k$ whose induced $\mathbb{G}_m$-torsor corresponds to $\mathscr{O}(n)$ is identified with $H^1(\mathbb{P}^1_k,\mathscr{O}(n))$. For $n \le 2$, the latter is a nonzero $k$-vector space (see \cite{EGAIII1}*{proposition 2.1.15}), so not a singleton.
\end{remark}

\section{The Birkhoff and Cartan decompositions for quasi-reductive groups}

The classification \Cref{thm:torsors-P1} for torsors over $\mathbb{P}^1_k$ under a quasi-reductive $k$-group and its proof yield the Birkhoff and the Cartan decompositions in \Cref{thm:Birkhoff,thm:Cartan}. For the Iwasawa decomposition, whose proof relies on a case of \Cref{thm:main-GS}, see \Cref{thm:Iwasawa}.

\csub[The Birkhoff decomposition]

The final part of the following decomposition result is similar to \cite{Gr-presheaf}*{Theorem 3.6}.

\begin{theorem}[Birkhoff Decomposition]\label{thm:Birkhoff}
For a $k$-group scheme $G$ locally of finite type with $G^{\mathrm{sm},\, \mathrm{lin}}$ quasi-reductive and a maximal split $k$-torus $S \le G$, we have
\be\ba\label{eqn:Birkhoff}
\textstyle G(k[t^{\pm1}])&= \tst\coprod_{\lambda \in \mathrm{Hom}_{k\text{-\upshape{gp}}}(\mathbb{G}_m,\, S)/N_G(S)(k)}G(k[t^{-1}])t^\lambda G(k[t]),\\
\textstyle G(k(t))&= \tst\coprod_{\lambda\in \mathrm{Hom}_{k\text{-\upshape{gp}}}(\mathbb{G}_m,\, S)/N_G(S)(k)}G(k[t^{-1}])t^\lambda G(k[t]_{(t)}),\\
\textstyle G(k\{t\}[\frac{1}{t}])&=\tst \coprod_{\lambda\in \mathrm{Hom}_{k\text{-\upshape{gp}}}(\mathbb{G}_m,\, S)/N_G(S)(k)}G(k[t^{-1}])t^\lambda G(k\{t\}),\\
\textstyle G(k\llp t\rrp)&=\tst \coprod_{\lambda\in \mathrm{Hom}_{k\text{-\upshape{gp}}}(\mathbb{G}_m,\, S)/N_G(S)(k)}G(k[t^{-1}])t^\lambda G(k\llb t\rrb),
\ea\ee
and
\be \label{eqn:alg-Gr}
\textstyle G(k[t^{\pm1}])/G(k[t]) \xrightarrow{\sim} G(k(t))/G(k[t]_{(t)}) \xrightarrow{\sim} G(k\{t\}[\frac{1}{t}])/G(k\{t\}) \xrightarrow{\sim} G(k\llp t\rrp)/G(k\llb t\rrb).
\ee
\end{theorem}

\bpf
By \eqref{eqn:reduce-to-gred}, the decompositions in question are insensitive to replacing $G$ by $G^{\mathrm{gred}}$, so we may assume that $G$ is smooth. Let $A$ denote either $k[t]$, or $k[t]_{(t)}$, or $k\{t\}$, or $k\llb  t\rrb  $, depending on the respective decomposition in question, so that in all cases $A$ is a filtered direct limit of flat, finitely presented $k[t]$-algebras (by the Popescu theorem \cite{SP}*{Theorem~\href{https://stacks.math.columbia.edu/tag/07GC}{07GC}} when $A = k\llb t\rrb$). By patching \cite{MB96}*{Theorem 5.5}, the double coset space $G(k[t^{-1}])\backslash G(A[\frac 1t])/G(A)$ is identified with the set of isomorphism classes of those $G$-torsors over $\mathbb{P}^1_k$ that trivialize both over $\mathbb{P}^1_k \setminus \{ t = 0\}$ and over $\mathrm{Spec}(A)$. By \Cref{thm:torsors-P1}, such $G$-torsors are all induced from the $\mathbb{G}_m$-torsor $\mathscr{O}(1)$ via some $k$-homomorphism $\lambda \colon \mathbb G_m \rightarrow S$, with two $G$-torsors being isomorphic if and only if the corresponding cocharacters are $N_G(S)(k)$-conjugate. Moreover, when $G = \mathbb{G}_m$, the class in the double coset space $A^\times \backslash A[\frac 1t]^{\times}/ k[t^{-1}]^\times$ that corresponds to $\mathscr{O}(1)$ is given by the element $t$. Thus, the functoriality gives the desired decompositions \eqref{eqn:Birkhoff}.

The maps in \eqref{eqn:alg-Gr} are all surjective by \eqref{eqn:Birkhoff}. Their injectivity follows from \cite{MB96}*{Theorem~5.5}, which implies, for instance, that a $k[t^{\pm1}]$-point of $G$ that extends to a $k\llb t\rrb$-point when restricted to $k\llp t\rrp$ already extends to a $k[t]$-point.
\epf

\begin{remark}
When one approaches the decompositions \eqref{eqn:Birkhoff} purely group-theoretically, the critical part to argue is that the union of double cosets on the right side of a desired equality form a subgroup of the left side. This is quite concrete but tends to be rather delicate even in the reductive case, compare with \cite{Rag94}*{proof of Theorem 3.4}.
\end{remark}

\csub[The Cartan decomposition]

To argue the Cartan decomposition in \Cref{thm:Cartan}, we adapt ideas from \S\ref{sec:torsors-P1} that gave the Birkhoff decomposition. More precisely, we combine the Auslander--Buchsbaum extension \Cref{thm:AB-extn} for torsors under quasi-reductive groups with the approach to the Cartan decomposition introduced in \cite{AHHL21} in the reductive case. The relevant analogue of \Cref{lem:maps-full} is the following lemma.

\begin{lemma}\label{lem:classify-S}
Let $G$ be a smooth $k$-group scheme, let $\mathcal{O}$ be a Henselian discrete valuation ring that is a $k$-algebra whose residue field $\kappa$ is separable over $k$, let $\pi \in \mathcal{O}$ be a uniformizer, and set $S \colonequals [\mathrm{Spec}(\mathcal{O}[s,s^{\prime}]/(ss^{\prime}-\pi))/\mathbb{G}_m]$ where $\mathbb{G}_m$ acts over $\mathcal{O}$ by scaling $s$ \up{resp.,~$s^{\prime}$} via the character of weight $1$ \up{resp.,~$-1$}. Restriction to $s = s^{\prime} = 0$ gives
\[
H^1(S, G) \cong H^1(\mathbf{B}\mathbb{G}_{m,\, \kappa}, G).
\]
\label{eqn:G-tors-S}
Moreover, for a $G$-torsor $E$ over $S$, the following are equivalent\ucolon
\benumr
    \item\label{m:CS-i} $E$ trivializes over the source of the map $\mathrm{Spec}(\kappa) \rightarrow \mathbf{B}\mathbb{G}_{m,\, \kappa} \subset S$\uscolon
    \item\label{m:CS-ii} $E$ trivializes over $S[\frac{1}{s}]$\uscolon
    \item\label{m:CS-iii} $E$ trivializes over $S[\frac{1}{s^{\prime}}]$.
\end{enumerate}
\end{lemma}

\bpf
Our $\mathcal O$ is Henselian and $\kappa$ is separable over $k$, so $\mathcal O$ is an algebra over $\kappa$. Thus, we may replace $k$ by $\kappa$ to reduce to when $\kappa = k$. Then \eqref{eqn:G-tors-S} amounts to saying that $G$-torsors over $S$ descend uniquely to $G$-torsors over $\mathbf B\mathbb G_{m,\, k}$. Thus, since both $S[\frac 1s]$ and $S[\frac 1{s^\prime}]$ are isomorphic to $\mathrm{Spec}(\mathcal O)$, by restricting to the residue field of the latter and using the smoothness of $G$ and the Henselianity of $\mathcal O$ (see \cite{Hitchin-torsors}*{Theorem 2.1.6 (a)}), we find that \eqref{eqn:G-tors-S} gives the equivalence of \ref{m:CS-i}--\ref{m:CS-iii}.

Overall, it remains to descend every $G$-torsor $E$ over $S$ to $\mathbf B\mathbb{G}_{m,\, k}$. For this, by twisting \cite{Gir71}*{chapitre~III, remarque 2.6.3}, we may replace $G$ by an inner form (see \S\ref{pp:k-groups-lft}) to force $E|_{s=s^{\prime}=0}$ to trivialize over the source of the map $\mathrm{Spec}(k) \rightarrow \mathbf B\mathbb G_{m,\, k}$. Since the closed $\{s^{\prime}=0\}\subset S$ is $[\mathbb{A}^1_k/\mathbb{G}_m]$, we then conclude from \eqref{eqn:M} that $E|_{s^{\prime}=0}$ is generically trivial. The closed point of $S[\frac 1s] \cong \mathrm{Spec}(\mathcal{O})$ is the generic point of $\{s^\prime = 0\} \subset S$, so, since $\mathcal{O}$ is Henselian and $G$ is smooth, we conclude that $E|_{S[\frac 1s]}$ is trivial (see \cite{Hitchin-torsors}*{Theorem 2.1.6 (a)}), in particular, that $E$ is generically trivial. At this point, since $S$ is geometrically regular over $k$, \Cref{thm:pseudo-extend}~\ref{m:PE-i} ensures that $E$ reduces to a generically trivial $G^{\mathrm{sm},\,\mathrm{lin}}$-torsor over $S$. Thus, we may assume that $G$ is affine, in which case, since $\mathcal O$ is Henselian, \eqref{eqn:G-tors-S} is a special case of \cite{Wed24}*{Corollary 2.9}.
\epf

\begin{theorem}[Cartan Decomposition]\label{thm:Cartan}
For a $k$-group scheme $G$ locally of finite type with $G^{\mathrm{sm},\, \mathrm{lin}}$ quasi-reductive, a $k$-algebra $\mathcal{O}$ that is a Henselian discrete valuation ring whose residue field is separable over $k$, a uniformizer $\pi \in \mathcal{O}$, and $K \colonequals \mathrm{Frac}(\mathcal{O})$,
\be\label{eqn:Cartan}
\textstyle G(K) = \coprod_{\lambda \in \mathrm{Hom}_{\mathcal{O}\text{-\upshape{gp}}}(\mathbb{G}_{m,\, \mathcal{O}},\, G)/G(\mathcal{O})} G(\mathcal{O})\pi^\lambda G(\mathcal{O}),
\ee
where $\pi^\lambda \colonequals \lambda(\pi) \in G(\mathcal{O})$; moreover, for a maximal split $k$-torus $S \le G$,
\[\ba
\textstyle G(k\{t\}[\frac{1}{t}]) &\tst= \coprod_{\lambda\in \mathrm{Hom}_{k\text{-\upshape{gp}}}(\mathbb{G}_m,\, S)/N_G(S)(k)} G(k\{t\})t^\lambda G(k\{t\}),\\
\textstyle G(k\llp t\rrp ) &\tst= \coprod_{\lambda\in \mathrm{Hom}_{k\text{-\upshape{gp}}}(\mathbb{G}_m,\, S)/N_G(S)(k)} G(k\llb  t\rrb  )t^\lambda G(k\llb  t\rrb  ).
\ea \]
\end{theorem}

\bpf
By \eqref{eqn:reduce-to-gred} and the separability of $K$ over $k$, which results from that of the residue field of $\mathcal O$ (see the sentence containing \eqref{eqn:pseudo-fpc}), the desired decompositions  are insensitive to replacing $G$ by $G^{\mathrm{gred}}$, so we may assume that $G$ is smooth. Moreover, $\mathcal O$ is Henselian and its residue field $\kappa$ is separable over $k$, so $\mathcal O$ is an algebra over $\kappa$. Thus, we may replace $k$ by $\kappa$ to reduce to when $\kappa = k$ (see also \S\ref{pp:fundamental-filtration}). By functoriality, this implies, in particular, that the pullback map
$$
\mathrm{Hom}_{\mathcal{O}\text{-gp}}(\mathbb{G}_{m,\, \mathcal{O}},\, G)/G(\mathcal O) \rightarrow\mathrm{Hom}_{k\text{-gp}}(\mathbb{G}_{m,\, k},\, G)/G(k)
$$
is surjective. It is also injective: for this, since $G(\mathcal O) \twoheadrightarrow G(k)$ by smoothness (see \cite{EGAIV4}*{th\'eor\`eme~18.5.17}), it suffices to note that, for any two $\mathcal O$-homomorphisms $\lambda, \lambda^\prime \colon \mathbb{G}_{m,\, \mathcal{O}} \rightarrow G$ that agree over the residue field $k$ and, as one checks over $K$, that must factor through $G^{\mathrm{sm},\, \mathrm{lin}}$, by \cite{SGA3II}*{expos\'e XI, corollaires 5.2, 5.4}, the subfunctor $\mathrm{Transp}_{G^{\mathrm{sm},\, \mathrm{lin}}}(\lambda, \lambda^\prime) \subset G^{\mathrm{sm},\, \mathrm{lin}}$ parametrizing sections that conjugate $\lambda$ to $\lambda^\prime$ is a residually trivial $Z_{G^{\mathrm{sm},\, \mathrm{lin}}}(\lambda)$-torsor over $\mathcal{O}$, which must then be trivial because $Z_{G^{\mathrm{sm},\, \mathrm{lin}}}(\lambda)$ is smooth (see \emph{loc.~cit.}~or \S\ref{pp:Gm-actions}). The bijectivity we just argued and \Cref{lem:SGconj} imply that it suffices to settle \eqref{eqn:Cartan} and with the index set replaced by~$\mathrm{Hom}_{k\text{-gp}}(\mathbb{G}_{m,\, k},\, G)/G(k)$.

To argue \eqref{eqn:Cartan}, we imitate the method of Alper--Heinloth--Halpern-Leistner from the reductive case \cite{AHHL21}. We consider the stack $S\ce [\mathrm{Spec}(\mathcal{O}[s,s^{\prime}]/(ss^{\prime}-\pi))/\mathbb G_m]$ of \Cref{lem:classify-S}. Its open $S \setminus \{ s = s^\prime = 0\}$ is the glueing of two copies of $\mathrm{Spec}(\mathcal O)$ along $\mathrm{Spec}(K)$. Thus, $G(\mathcal O)\backslash G(K)/G(\mathcal O)$ is identified with the set of isomorphism classes of those $G$-torsors over $S \setminus \{ s = s^\prime = 0\}$ that are trivial on both copies of $\mathrm{Spec}(\mathcal O)$. Since $S$ is geometrically regular over $k$ and $\mathcal O[s, s^\prime]/(ss^\prime - \pi)$ is $2$-dimensional, the quasi-reductivity assumption and \Cref{thm:AB-extn} ensure that $G$-torsors over $S \setminus \{ s = s^\prime = 0\}$ extend uniquely to those over $S$. Thus, by \Cref{lem:classify-S}, pullback of $G$-torsors along the structure map $S \rightarrow \mathbf B\mathbb{G}_{m,\, k}$ identifies the set $G(\mathcal O)\backslash G(K)/G(\mathcal O)$ with the set of isomorphism classes of those $G$-torsors over $\mathbf B\mathbb{G}_{m,\, k}$ that trivialize over $\mathrm{Spec}(k)$. By the proof of \Cref{thm:torsors-P1}, this latter set is identified with $\mathrm{Hom}_{k\text{-gp}}(\mathbb{G}_{m,\, k}, G)/G(k)$. Since over $K$ both $s$ and $s^\prime$ are invertible with $s s^\prime = \pi$, overall a $\lambda \in \mathrm{Hom}_{k\text{-gp}}(\mathbb{G}_{m,\, k}, G)/G(k)$ corresponds to the double coset of $\pi^\lambda$.
\epf

The Cartan decomposition gives the following integrality property of rational points of anisotropic quasi-reductive groups that generalizes pseudo-completeness of wound unipotent groups (see \Cref{prop:PP-groups}~\ref{m:PPG-i}). In the reductive case, this is a theorem of Bruhat and Tits usually proved group-theoretically, for instance, by using buildings, see \cite{Pra82}, \cite{Guo22}*{Proposition 6 (2)}, and \cite{FG21}*{Corollary 3.8}, while our argument is algebro-geometric and in essence originates in \cite{AHHL21}.

\begin{corollary}\label{cor:anisotropic}
For a $k$-group scheme $G$ locally of finite type with $G^{\mathrm{sm}, \mathrm{lin}}$ quasi-reductive, a $k$-algebra $\mathcal{O}$ that is a discrete valuation ring whose residue field $\kappa$ is separable over $k$ and such that $G_{\kappa}$ contains no nontrivial split $\kappa$-torus, and $K \colonequals \mathrm{Frac}(\mathcal{O})$, we have
\[
G(\mathcal{O}) = G(K).
\]
\end{corollary}

\begin{proof}
Certainly, $G(\mathcal{O}) \subset G(K)$. Conversely, for checking that every $K$-point of $G$ extends to an $\mathcal{O}$-point, by considering generators and relations for coordinate algebras of elements of an affine open cover of $G$, we may replace $\mathcal{O}$ by its completion. Once $\mathcal{O}$ is complete, the Cartan decomposition \eqref{eqn:Cartan} (and the first part of the proof of \Cref{thm:Cartan}) implies that $G(\mathcal{O}) = G(K)$. 
\end{proof}

We now use \Cref{thm:Cartan} to quickly reprove some of the main results of \cite{CGP15}*{Appendix C.3}.

\begin{corollary}\label{cor:Ga-Gm}
If a $k$-group scheme $G$ locally of finite type with $G^{\mathrm{sm},\, \mathrm{lin}}$ quasi-reductive has $\mathbb{G}_{a,\, k}$ as a $k$-subgroup, then it also has $\mathbb{G}_{m,\, k}$ as a $k$-subgroup.
\end{corollary}

\begin{proof}
If $\mathbb{G}_{a,\, k} \le G$, then $G(k\llb t\rrb) \subsetneq G(k\llp t\rrp)$. By \Cref{cor:anisotropic} (applied with $\mathcal{O} = k\llb t\rrb$), this means that $G$ has $\mathbb{G}_{m,\, k}$ as a $k$-subgroup. 
\end{proof}

The following consequence of the Cartan decomposition generalizes \cite{BT71}*{corollaire 3.7}.

\begin{corollary}\label{cor:split-parabolic}
Every split unipotent $k$-subgroup $U$ of a $k$-group scheme $G$ locally of finite type lies in $\mathscr{R}_{\mathrm{us},\, k}(P)$ for some pseudo-parabolic $k$-subgroup $P \leq G^{\mathrm{sm},\, \mathrm{lin}}$. The maximal split unipotent $k$-subgroups of $G$ are precisely the unipotent radicals of the minimal pseudo-parabolic $k$-subgroups of $G^{\mathrm{sm},\, \mathrm{lin}}$, and the latter are pairwise $G^{\mathrm{sm},\, \mathrm{lin}}(k)$-conjugate.
\end{corollary}

\bpf
The pairwise $G^{\mathrm{sm},\, \mathrm{lin}}(k)$-conjugacy of the minimal pseudo-parabolic $k$-subgroups, and so also of their (split) unipotent radicals, mentioned in the statement for convenience of later reference, is a result of Borel--Tits \cite{CGP15}*{Theorem C.2.5}. Moreover, by \cite{CGP15}*{Proposition 3.5.14}, an inclusion of pseudo-parabolic $k$-subgroups induces the opposite inclusion on their split unipotent radicals, so it suffices to settle the claim about $U$.

For the latter, since $U \le G^{\mathrm{sm},\, \mathrm{lin}}$, we lose no generality by assuming that $G$ is connected, smooth, and affine. Since $U$ is split, by \cite{CGP15}*{Lemma C.2.2}, it lies in some $G(k)$-conjugate of each pseudo-parabolic $k$-subgroup $P \le G$. Thus, by \S\ref{pp:pseudo-parabolic}, we may iteratively replace $G$ by such conjugates of $P$ to reduce to when $G$ has no proper pseudo-parabolic $k$-subgroups. Moreover, since $U$ is split, it is enough to show that $U \le \mathscr{R}_{\mathrm{u},\, k}(P)$, so we may replace $G$ by its largest pseudo-reductive quotient $G^{\mathrm{pred}}$ to also assume that $G$ is pseudo-reductive (see \S\ref{pp:fundamental-filtration}~\ref{m:RuskG} and \S\ref{pp:pseudo-parabolic}). In this pseudo-reductive case with no proper pseudo-parabolic $k$-subgroups, we claim that $U$ is trivial (so that $P = G$ suffices). Indeed, otherwise \Cref{cor:Ga-Gm} would supply some $\mathbb{G}_{m,\, k} \le G$, which, since $G$ has no nontrivial pseudo-parabolic $k$-subgroups, would be central. By \cite{CGP15}*{Proposition 2.2.12 (3)}, we could then replace $G$ by its quotient $\overline{G}$ by this central $\mathbb{G}_{m,\, k}$: by \cite{CGP15}*{Proposition 1.2.4, Lemma 9.4.1}, this $\overline{G}$ is still pseudo-reductive because $\mathrm{Ext}^1_k(\mathbb G_{a,\, k}, \mathbb G_{m,\, k}) = 0$ thanks to \cite{DG70}*{chapitre III, section 6, no.~5, corollaire 5.2}. By iterating this reduction finitely many times, we would arrive at the case when $G$ has no $\mathbb{G}_{m,\, k}$ as a $k$-subgroup, so that $U = 1$ by \Cref{cor:Ga-Gm}, as desired.
\epf

\section{Torsors over a relative $\mathbb{P}^1_A$ under constant groups}

In \S\ref{sec:swise}, we upgrade the analysis of torsors over $\mathbb P^1_k$ carried out in \S\ref{sec:torsors-P1} to the analysis of torsors over $\mathbb P^1_A$ for any semilocal $k$-algebra $A$. This is a critical step of the geometric approach to the Grothendieck--Serre question and it rests on structural results about the Whitehead group of a quasi-reductive group that we establish in \S\ref{sec:WkG}.

\csub[Unramifiedness of the Whitehead group in the quasi-reductive setting]\label{sec:WkG}

In \Cref{prop:W-unr}, we prove an unramifiedness property of the subgroup $G(k)^+$ generated, roughly speaking, by the ``elementary matrices'' of a quasi-semisimple $k$-group $G$. In the classical case of semisimple groups, the corresponding results about $G(k)^+$ are primarily due to Borel and Tits \cite{BT73}, with an excellent overview by Gille \cite{Gil09}.

\bpp[Notation]
Throughout this section, we fix a connected, smooth, affine $k$-group $G$.
\epp

\bpp[The subgroup $G(k)^+ \le G(k)$] \label{pp:Gkplus} 
As in \cite{BT73}*{section 6.1}, we let $G(k)^+$ denote the (normal) subgroup of $G(k)$ generated by the $U(k)$ for split unipotent $k$-subgroups $U \le G$. By \S\ref{pp:fundamental-filtration}~\ref{m:RuskG}, this $G(k)^+$ is functorial in $G$  and $k$.

By \Cref{cor:split-parabolic}, our $G(k)^+$ is generated by the $G(k)$-conjugates of $\mathscr{R}_{\mathrm{us},\, k}(P)(k)$ for any single minimal pseudo-parabolic $k$-subgroup $P \le G$. If $k$ is infinite, then we even have
$$
G(k)^+ = \langle\mathscr{R}_{\mathrm{us},\,k}(P_\lambda)(k), \mathscr{R}_{\mathrm{us},\,k}(P_{-\lambda})(k)\rangle
$$
for any minimal pseudo-parabolic $k$-subgroup $P_\lambda \le G$ (with notation as in \S\ref{pp:pseudo-parabolic}): for this, it is enough to show that $G(k)$ normalizes the right hand side, or that 
$$
G(k) \overset{?}{=}  U_G(\lambda)(k)U_G(-\lambda)(k)U_G(\lambda)(k)Z_G(\lambda)(k).
$$
For the latter, since $P_G(\lambda) \cong U_G(\lambda) \rtimes Z_G(\lambda)$ and $U_G(-\lambda)\times P_G(\lambda)$ is open in $G$ (see \S\ref{pp:Gm-actions}), it suffices to cover $G$ by the $U_G(\lambda)(k)$-translates of $U_G(-\lambda)P_G(\lambda)$. By multiplying the inverse of the open with the open itself, we get $U_G(\lambda)U_G(-\lambda)P_G(\lambda) = G$. Since $k$ is infinite, the split unipotent $U_G(\lambda)$ has a dense set of $k$-points (see \S\ref{pp:Gm-actions}), so, for each $g \in G(\overline{k})$, \cite{CGP15}*{top of p. 587} gives an $u \in U_G(\lambda)(k)$ whose inverse lies in the open $U_G(-\lambda)P_G(\lambda) g^{-1} \le G_{\overline{k}}$. Thus, $g \in u (U_G(-\lambda)P_G(\lambda))$.
\epp

\bpp[The Whitehead group $W(k, G)$]\label{pp:WkG}
The \emph{Whitehead group} of~$G$ is the quotient
\[
W(k, G) \colonequals G(k)/G(k)^+,
\]
compare with \cite{Gil09}*{Section~1} in the reductive case. It is functorial in~$G$ and in~$k$ (see~\S\ref{pp:Gkplus}) and only depends on the largest quasi-reductive quotient $G^{\mathrm{qred}}$:
\[
W(k, G) \xrightarrow{\sim} W(k, G^{\mathrm{qred}}),
\]
indeed, split unipotent groups have no nontrivial torsors over~$k$, so $G(k) \twoheadrightarrow G^{\mathrm{qred}}(k)$ and, by~\S\ref{pp:Gkplus} and~\S\ref{pp:pseudo-parabolic}, the preimage of $G^{\mathrm{qred}}(k)^+$ is precisely $G(k)^+$. Moreover, as we now argue, for every maximal $k$-split torus $S \le G$, we have
\begin{equation}\label{eqn:plus-ZGS}
G(k) = G(k)^+ Z_G(S)(k), \quad \text{equivalently,} \quad Z_G(S)(k) \twoheadrightarrow W(k, G).
\end{equation}
Indeed, \cite{CGP15}*{Remark~C.2.33} (applied with $\mathscr{G} = G^{\mathrm{qred}}(k)^+$) gives this with~$G$ replaced by $G^{\mathrm{qred}}$, so to deduce it for~$G$ it suffices to recall from \cite{CGP15}*{Lemma~C.2.31} that~$S$ maps onto a maximal $k$-split torus of $G^{\mathrm{qred}}$ and to apply the following sharpening of \cite{CGP15}*{Proposition~A.2.8}.
\epp

\begin{lemma}\label{lem:T-kernel}
For a surjection $\pi\colon G \twoheadrightarrow G^\prime$ of connected, smooth, affine $k$-groups with $\mathrm{Ker}(\pi)$ split unipotent, and a $k$-torus $T \le G$ with image $T^\prime \le G^\prime$, the kernel of the surjection
\[
Z_G(T) \twoheadrightarrow Z_{G^\prime}(T^\prime)
\]
is also split unipotent. In particular,
\[
Z_G(T)(k) \twoheadrightarrow Z_{G^\prime}(T^\prime)(k).
\]
\end{lemma}

\bpf
By \cite{CGP15}*{Propositions A.2.5 and A.2.8}, the map $Z_G(T) \twoheadrightarrow Z_{G^{\prime}}(T^\prime)$ is indeed a surjection of connected, smooth, affine $k$-groups. For the claim about the kernel, we may replace $G$ by the preimage of $Z_{G^{\prime}}(T^\prime)$ and then base change to $k^s$ to assume that $T^\prime$ is central in $G^\prime$ and $k$ is separably closed (see \S\ref{pp:fundamental-filtration}~\ref{m:RuskG}). Then, for any cocharacter $\lambda \colon \mathbb{G}_{m,\, k} \rightarrow T \le G$, by \S\ref{pp:Gm-actions}, we have 
$$
\mathrm{Ker}(\pi) = U_{\mathrm{Ker}(\pi)}(-\lambda) \times Z_{\mathrm{Ker}(\pi)}(\lambda) \times U_{\mathrm{Ker}(\pi)}(\lambda)
$$
as $k$-schemes. In particular, $Z_{\mathrm{Ker}(\pi)}(\lambda)$ inherits split unipotence from $\mathrm{Ker}(\pi)$, see \S\ref{pp:fundamental-filtration}~\ref{m:RuskG}. Thus, thanks to the exact sequence
$$
1 \rightarrow Z_{\mathrm{Ker}(\pi)}(\lambda) \rightarrow Z_G(\lambda) \rightarrow G^\prime \rightarrow 1
$$
(see \S\ref{pp:Gm-actions}) and the evident containment $Z_G(T) \le Z_G(\lambda)$, we may replace $G$ by $Z_G(\lambda)$. Since $k$ is separably closed, if we iterate this for a well-chosen finite set of $\lambda$’s, we reduce to when $Z_G(T) = G$. Then the kernel in question is $\mathrm{Ker}(\pi)$, so is split unipotent.
\epf

The final input needed for the promised unramifiedness property of the Whitehead group is the following extension of the existence of Levi subgroups \cite{CGP15}*{Theorem 3.4.6} to the quasi-reductive case. For the sake of focus on intended use and since our proof essentially combines heavy inputs from \cite{CGP15}, we do not aim for a finer statement that would discuss uniqueness of Levi subgroups.

\begin{lemma}\label{lem:Levi}
For a quasi-reductive $k$-group $G$ that has a split maximal $k$-torus $S \le G$, there exists a split reductive $k$-subgroup $\mathscr{G} \le G$ containing $S$ such that $\mathscr{G}_{\overline{k}} \xrightarrow{\sim} (G_{\overline{k}})^{\mathrm{pred}}$.
\end{lemma}

\bpf
If $k$ is finite, then we may choose $\mathscr{G} \colonequals G$, so we assume that $k$ is infinite. Then \cite{CGP15}*{Theorem C.2.30} supplies a split reductive $k$-subgroup $\mathscr{G}$ containing $S$: the assumptions of \emph{loc.~cit.}~are met because $G$ is quasi-reductive, see \cite{CGP15}*{Proposition B.4.4, Theorem C.2.15, and bottom of~p.~630}. It remains to show that the map $\iota \colon \mathscr{G} \rightarrow G^{\mathrm{pred}}$ is injective: then, by the characterization of our $\mathscr{G}$ given in \cite{CGP15}*{Theorem C.2.30} (in particular, by its uniqueness aspect applied to $\iota(\mathscr{G}) \le G^{\mathrm{pred}}$) and by \cite{CGP15}*{Theorem 3.4.6}, we will get that $\mathscr G \cong \iota(\mathscr G)$ satisfies $\mathscr{G}_{\overline{k}} \xrightarrow{\sim} (G_{\overline{k}})^{\mathrm{pred}}$. 

The injectivity of $\iota$ is not completely general because in characteristic $2$ some reductive groups do have nontrivial normal unipotent subgroups, see \cite{Vas05}*{Theorem 1.2}. Nevertheless, by \cite{CGP15}*{Proposition B.4.4}, our $S$ commutes with the wound unipotent $\mathscr{R}_{\mathrm{u},\, k}(G)$, so also with $\mathrm{Ker}(\iota)$. The $S$-weight decomposition of $\mathrm{Lie}(\mathscr G)$ supplied by \cite{CGP15}*{Theorem C.2.30} then shows that $\mathrm{Lie}(\mathrm{Ker}(\iota)) \subset \mathrm{Lie}(S)$. However, $\mathrm{Ker}(\iota)$ is unipotent, so $\mathrm{Ker}(\iota) \cap S = 0$, and hence $\mathrm{Lie}(\mathrm{Ker}(\iota)) = 0$. We get that $\mathrm{Ker}(\iota)$ is \'etale, so, as a normal subgroup of a connected $\mathscr{G}$, it must be central. Since $\mathscr G$ is reductive, its center is of multiplicative type, so our unipotent $\mathrm{Ker}(\iota)$ is trivial, as desired.
\epf

\begin{proposition}\label{prop:W-unr}
Suppose that $G$ is quasi-semisimple, simply connected \up{see \uS\uref{pp:pred-iG}}, and has a split maximal $k$-torus $S$, let $\mathcal{O}$ be a discrete valuation ring whose residue field is separable over $k$, and set $K \colonequals \mathrm{Frac}(\mathcal O)$. The Whitehead group $W(K, G_K)$ is $\mathcal O$-unramified in the sense that it is a quotient of $G(\mathcal{O})$, in fact, even of $Z_G(S)(\mathcal O)$\ucolon
\[
G(K) = G(K)^+Z_G(S)(\mathcal{O}), \quad \text{equivalently,} \quad Z_G(S)(\mathcal O) \twoheadrightarrow W(K, G_K).
\]
\end{proposition}

\bpf
Granted the inputs established above, the argument is similar to the one given for semisimple groups in \cite{Gil09}*{lemme 4.5 (1)}. Namely, let $S \le G$ be a split maximal $k$-torus, so that, by \eqref{eqn:plus-ZGS},
$$
G(K) = G(K)^+Z_G(S)(K).
$$
By \cite{CGP15}*{Remark C.2.12 (1)}, the $k$-group $Z_G(S)$ is quasi-reductive and, by \Cref{rem:playful}, so is its central quotient $Z_G(S)/S$. The latter contains no $k$-torus, so it must even be wound unipotent, see \S\ref{pp:fundamental-filtration}~\ref{m:Gtor}. In particular, \Cref{prop:PP-groups}~\ref{m:PPG-i} ensures that $(Z_G(S)/S)(\mathcal O) = (Z_G(S)/S)(K)$, so 
$$
Z_G(S)(K) \subset S(K)Z_G(S)(\mathcal O).
$$
It remains to show that $S(K) \subset G(K)^+$, and for this we will use \Cref{lem:Levi}, according to which our split maximal $k$-torus $S$ is a maximal torus of a split reductive $k$-subgroup $\mathscr G \le G$ such that $\mathscr{G}_{\overline{k}} \xrightarrow{\sim} (G_{\overline{k}})^{\mathrm{pred}}$. This isomorphism and our assumptions imply that $\mathscr G$ is semisimple, simply connected. Then $S(K) \subset \mathscr{G}(K)^+$ by \cite{BT73}*{corollaire 6.8}, so that   $\mathscr G(K)^+ \subset G(K)^+$ by \S\ref{pp:Gkplus}.
\epf

\brem
For Henselian $\mathcal{O}$, analogously to semisimple groups treated in \cite{Gil09}*{lemme 4.5~(1)}, there ought to be a more general version of \Cref{prop:W-unr} in which instead of having a split maximal $k$-torus, our quasi-semisimple, simply connected $G$ merely has a sufficiently small pseudo-parabolic $K$-subgroup. As for our argument, a thorny aspect of such a generalization, is that, in general, it seems delicate to check that the split reductive subgroup supplied by \cite{CGP15}*{Theorem~C.2.30} inherits the simple connectedness from $G$ (compare with \cite{BT72}*{corollaire 4.6} in the reductive case). 
\erem

\csub[Sectionwise triviality of torsors over a relative $\mathbb{P}^1_A$]
\label{sec:swise}

The following \Cref{thm:swise} about $A$-sectionwise triviality of $G$-torsors over $\mathbb{P}^1_A$ for smooth $k$-groups $G$ and semilocal $k$-algebras $A$ is a critical final input for our \Cref{thm:main-GS}. It both uses and generalizes \Cref{thm:torsors-P1}, which treated the case when $A$ is a field, and it extends \cite{totally-isotropic}*{Theorem~3.5} (so also the main result of \cite{PS25}), which established a similar conclusion for torsors under reductive groups and used it as an input to establish cases of the Grothendieck--Serre conjecture. 

\bthm\label{thm:swise}
For a smooth $k$-group $G$, a semilocal $k$-algebra $A$, and a $G$-torsor $E$ over $\mathbb{P}^1_A$, if $s^*(E)$ is trivial for a single $s \in \mathbb P^1_A(A)$, then it is trivial for every such $s$.
\ethm

\bpf
By decomposing into components, we may assume that $\mathrm{Spec}(A)$ is connected and, by using the smoothness of $G$ and replacing $A$ by $A^{\mathrm{red}}$, also reduced (see \cite{Hitchin-torsors}*{Theorem 2.1.6 (a)} and \cite{SP}*{Lemma \href{https://stacks.math.columbia.edu/tag/0ALI}{0ALI}}). Since $A$ is semilocal, the $A$-automorphisms of $\mathbb P^1_A$ act transitively on $\mathbb P^1_A(A)$: for every $s \in \mathbb{P}^1_A(A)$, by considering the closed $A$-fibers, we may build an $s^\prime \in \mathbb A^1_A(A)$ disjoint from $s$, then change coordinates to first make $s^\prime$  be $\{t = 0\}$ and then make $s$ be $\{t = \infty\}$. Thus, we may assume that $E|_{\{ t= \infty\}}$ is trivial and need to argue that $E|_{\{ t = 0\}}$ is then also trivial.

\Cref{cor:reduce-to-Gt} allows us to replace $G$ by $G^{\mathrm{tor}}$ to reduce to when $G$ is connected smooth, affine, and generated by tori. In addition, since split unipotent groups have no nontrivial torsors over affine schemes, we may further replace $G$ by $G/\mathscr{R}_{\mathrm{us},\, k}(G)$ to reduce to when $G$ is also quasi-reductive. At this point, \S\ref{pp:DG} ensures that $G/\mathscr D(G)$ is a torus, a quotient of any maximal $k$-torus $T \le G$. By \cite{totally-isotropic}*{Lemma 3.2}, the $(G/\mathscr D(G))$-torsor $\overline{E}$ over $\mathbb P^1_A$ induced by $E$ is the inflation of the $\mathbb G_m$-torsor $\mathscr O(1)$ along some $A$-homomorphism $\mathbb G_{m,\, A} \rightarrow G/\mathscr D(G)$. Thus, for a trivialization $\iota$ of $E|_{A\llp t^{-1}\rrp}$, some $\overline{g} \in (G/\mathscr D(G))(A\llp t^{-1}\rrp )$ is such that the glueing of $\overline{E}|_{\mathbb A^1_A}$ and the trivial $(G/\mathscr D(G))$-torsor over $A\llb t^{-1}\rrb$ along the $\overline{g}$-translate of the trivialization $\overline{\iota}$ of $\overline{E}|_{A\llp t^{-1}\rrp}$ induced by $\iota$ gives the trivial $(G/\mathscr D(G))$-torsor over $\mathbb P^1_A$ (for the relevant glueing technique that includes non-Noetherian $A$, see \cite{Hitchin-torsors}*{Lemma~2.2.11~(b)}). By \cite{Hitchin-torsors}*{Lemma 3.1.6}, we have
$$
(G/\mathscr D(G))(A\llp t^{-1}\rrp) \cong X_*(G/\mathscr D(G))(A) \times (G/\mathscr D(G))(A\llb t^{-1}\rrb)
$$
where a cocharacter $\alpha\colon \mathbb G_{m,\, A} \rightarrow G/\mathscr D (G)$ maps to $\alpha(t^{-1}) \in (G/\mathscr D(G))(A\llp t^{-1}\rrp)$. Elements of $(G/\mathscr D(G))(A\llb t^{-1}\rrb)$ do not affect the isomorphism class of the glued $(G/\mathscr D(G))$-torsor over $\mathbb P^1_A$, so we may assume that $\overline{g} = \alpha(t^{-1})$ for some such $\alpha$. Pulling $E$ back along the $d$-th power map $f_d\colon \mathbb P^1_A \rightarrow \mathbb P^1_A$ given by $[x : y] \mapsto [x^d : y^d]$ for $d > 0$ preserves both our assumption (the triviality of $E|_{\{t = \infty\}}$) and the desired conclusion (the triviality of $E|_{\{ t= 0\}}$), and it replaces $\alpha(t^{-1})$ by $\alpha(t^{-d}) = \alpha^d(t^{-1})$, so we have the liberty of replacing $E$ by its pullback along any $f_d$.  For a sufficiently divisible $d > 0$, however, $\alpha^d$ lifts to a cocharacter $\widetilde\alpha\colon \mathbb G_{m,\, A} \rightarrow T$. Thus, we may first replace $E$ by its pullback along $f_d$ and then by the glueing of $E|_{\mathbb A^1_A}$ and the trivial $G$-torsor over $A\llb t^{-1}\rrb$ along the $\widetilde{\alpha}(t^{-1})$-translate of $\iota$ to reduce to when $\overline{E}$ is a trivial $(G/\mathscr D(G))$-torsor. By the rigidity lemma \cite{MFK94}*{Proposition~6.1}, the trivializations of $\overline{E}$ over $\mathbb P^1_A$ are pulled back bijectively to those of $\overline{E}|_{\{ t = \infty\}}$. In particular, $E$ reduces to a $\mathscr{D}(G)$-torsor over $\mathbb P^1_A$ whose restriction to $\{t = \infty\}$ is trivial, so we may replace $G$~by~$\mathscr{D}(G)$. 

We iterate these reductions---we replace $G$ by $\mathscr D(G)$, then $\mathscr{D}(G)$ by $\mathscr D(G)^{\mathrm{tor}}$, then $\mathscr D(G)^{\mathrm{tor}}$ by $\mathscr{D}(\mathscr D(G)^{\mathrm{tor}})$, and so on---until we are left with the case when $G$ is quasi-semisimple (see \S\ref{pp:DG}). By \cite{CP16}*{Theorem 5.1.3}, there then exist a commutative, affine $k$-group $Z$ that has no nontrivial unipotent $k$-subgroups, a quasi-semisimple, simply connected $k$-group $\widetilde{G}$, and a central extension
$$
1 \rightarrow Z \rightarrow \widetilde{G} \rightarrow G \rightarrow 1
$$
over $k$. The structure theorem \cite{DG70}*{chapitre IV, section 3, th\'eor\`eme 1.1} for commutative, affine $k$-groups ensures that $Z$ is an extension of a unipotent $k$-group $U$ by a $k$-group $M$ of multiplicative type. In addition, since $\widetilde{G}_{\overline{k}}/\mathscr{R}_{\mathrm{u},\, \overline{k}}(\widetilde{G}_{\overline{k}})$ is semisimple, $M$ has no nontrivial $k$-subtori, so it is finite (see \cite{SGA3II}*{expos\'e~XII, proposition 1.12}).

Liftings of $E$ to a $\widetilde{G}$-torsor form a $Z$-gerbe $\mathscr{Z}$ over $\mathbb{P}^1_A$ whose restriction to $\{t = \infty\}$ is trivial. Moreover, \Cref{rem:unip-higher} ensures that $\mathscr{Z}$ reduces to an $M$-gerbe $\mathscr{M}$ over $\mathbb{P}^1_A$. By \cite{totally-isotropic}*{Lemma 3.3} (with \cite{CTS21}*{Theorem 6.1.3}), pulling back along the $d$-th power map $f_d$ makes $\mathscr{M}$ descend to an $M$-gerbe over $A$ whenever $d$ kills $M$. In effect, we may once more replace $E$ by its pullback along $f_d$ to make $\mathscr{Z}$ descend to $A$. By the triviality at $\{t = \infty\}$, this makes $\mathscr{Z}$ (noncanonically) the trivial $Z$-gerbe $\mathbf B Z$, so that the map $\mathscr{Z}(\mathbb P^1_A) \rightarrow \mathscr Z(A)$ given by restricting to $t = \infty$ is essentially surjective. This implies that there is a lifting of $E$ to a $\widetilde{G}$-torsor $\widetilde{E}$ over $\mathbb P^1_A$ such that $\widetilde{E}|_{\{ t = \infty\}}$ is trivial. This allows us to replace $G$ by $\widetilde{G}$ to reduce to the case when $G$ is quasi-semisimple, simply connected.

Once $G$ is quasi-semisimple, simply connected, we let $\ell/k$ be a finite, separable field extension such that $G_\ell$ has a split maximal $\ell$-torus $T \le G_\ell$, we set $d \colonequals [\ell :k ]$, and we choose an $\ell$-cocharacter $\lambda \colon \mathbb G_{m,\, \ell} \rightarrow T \le G_\ell$ such that $P_{\lambda} \le G_\ell$ is a pseudo-Borel $\ell$-subgroup (see \S\ref{pp:pseudo-parabolic}). Similarly to above, we now replace $E$ by its pullback along $f_{d!}$: the point of this maneuver is that, by \Cref{thm:torsors-P1}, then for every residue field $\kappa$ of $A$, the $(G_{\kappa})^{\mathrm{qred}}$-torsor over $\mathbb{P}^1_{\kappa}$ induced by $E$ reduces to the $\mathbb{G}_m$-torsor $\mathscr{O}(d!)$ along some $\kappa$-cocharacter $\mathbb G_{m,\, \kappa} \rightarrow (G_{\kappa})^{\mathrm{qred}}$, so this $(G_{\kappa})^{\mathrm{qred}}$-torsor is trivial away from every divisor of $\mathbb P^1_{\kappa}$ of degree $\le d$. To construct a relevant such divisor, we choose a nonzero primitive element that generates $\ell/k$, consider the resulting embedding $\mathrm{Spec}(\ell) \subset \mathbb G_{m,\, k}$, and let $Y \cong \mathrm{Spec}(A^\prime) \subset \mathbb G_{m,\, A}$ be its $A$-(finite \'etale) base change, so that $A^\prime \cong \ell \otimes_k  A$. The formal completion of $\mathbb P^1_A$ along $Y$ has a formal power series ring $A^\prime\llb \tau\rrb$ as its coordinate ring. Moreover, since $\ell/k$ is separable, we have $Y_{k_{\mathfrak m}} \cong \mathrm{Spec}(A^\prime \otimes_A k_{\mathfrak m}) \cong \mathrm{Spec}(\ell \otimes_k k_{\mathfrak m})$ for every maximal ideal $\mathfrak m \subset A$, with $\ell \otimes_k k_{\mathfrak m}$ a finite product of finite, separable field extensions of $k_{\mathfrak m}$. Therefore, \Cref{prop:W-unr} (with \S\S\ref{pp:Gm-actions}--\ref{pp:pseudo-parabolic} and \S\ref{pp:Gkplus}) ensures that each coset in $(G_{k_{\mathfrak m}})^{\mathrm{qred}}((A^\prime \otimes_A k_{\mathfrak m})\llp\tau\rrp)/ (G_{k_{\mathfrak m}})^{\mathrm{qred}}((A^\prime \otimes_A k_{\mathfrak m})\llb \tau\rrb)$ is represented by a product of elements of $(U_{(G_{k_{\mathfrak m}})^{\mathrm{qred}}}(\pm \lambda))((A^\prime \otimes_A k_{\mathfrak m})\llp\tau\rrp)$. By \S\ref{pp:Gm-actions},
$$
1 \rightarrow U_{\mathscr R_{\mathrm{us},\, k_{\mathfrak m}}(G_{k_{\mathfrak m}})}(\pm \lambda) \rightarrow U_{G_{k_{\mathfrak m}}}(\pm \lambda) \rightarrow U_{(G_{k_{\mathfrak m}})^{\mathrm{qred}}}(\pm \lambda) \rightarrow 1
$$
are short exact sequences of split unipotent $k_{\mathfrak m}$-groups for every maximal $\mathfrak m \subset A$. The maps
\[
(U_G(\pm\lambda))((A^\prime \otimes_A k_{\mathfrak m})\llp \tau\rrp) \twoheadrightarrow (U_{(G_{k_{\mathfrak m}})^{\mathrm{qred}}}(\pm \lambda))((A^\prime \otimes_A k_{\mathfrak m})\llp\tau\rrp)
\]
are therefore surjective. In addition, the map $A^\prime((\tau)) \twoheadrightarrow \prod_{\mathfrak m}(A^\prime \otimes_A k_{\mathfrak m})((\tau))$ is surjective and this surjectivity persists upon applying $(U_G(\pm \lambda))(-)$ because the $\ell$-schemes $U_G(\pm \lambda)$ are both isomorphic to some $\mathbb A^n_\ell$ (see \S\ref{pp:fundamental-filtration}~\ref{m:RuskG}). Thus, overall, the following map is surjective: 
\be\label{eqn:surj}
\textstyle G(A'\llp \tau\rrp) \twoheadrightarrow \prod_{\mathfrak m}(G_{k_{\mathfrak m}})^{\mathrm{qred}}((A^\prime \otimes_A k_{\mathfrak m})\llp \tau\rrp)/ (G_{k_{\mathfrak m}})^{\mathrm{qred}}((A^\prime \otimes_A k_{\mathfrak m})\llb\tau\rrb).
\ee
By construction, for every maximal ideal $\mathfrak m \subset A$, the $(G_{k_{\mathfrak m}})^{\mathrm{qred}}$-torsor over $\mathbb{P}^1_{k_{\mathfrak m}}$ induced by $E$ trivializes over $\mathbb{P}^1_{k_{\mathfrak m}} \setminus Y_{k_{\mathfrak m}}$. Thus, by patching with the trivial torsor over $A^\prime\llb \tau\rrb$ (see \cite{Hitchin-torsors}*{Lemma~2.2.11~(b)}), the surjectivity \eqref{eqn:surj} allows us to modify $E$ along $Y$ without changing $E|_{\mathbb P^1_A \setminus Y}$, thereby reducing us to the case when $E$ induces a trivial $(G_{k_{\mathfrak m}})^{\mathrm{qred}}$-torsor over $\mathbb{P}^1_{k_{\mathfrak m}}$ for every maximal ideal $\mathfrak m \subset A$. Since split unipotent groups have no nontrivial torsors over $\mathbb P^1_{k_{\mathfrak m}}$, this means that $E$ induces a trivial $G_{k_{\mathfrak m}}$-torsor over $\mathbb{P}^1_{k_{\mathfrak m}}$for each $\mathfrak m$. Then, however, the deformation-theoretic \cite{totally-isotropic}*{Proposition~3.1} implies that $E$ is constant, so that $E|_{\{t = 0\}}$  is isomorphic to $E|_{\{t = \infty\}}$, and hence is trivial, as desired.
\epf

The following consequence of \Cref{thm:swise} is inspired by \cite{Gr-presheaf}*{Corollary 2.3 and Remark 2.4}.

\begin{corollary}\label{cor:FG-A}
For a smooth $k$-group $G$ and a semilocal $k$-algebra $A$, no nontrivial $G$-torsor over $A$ trivializes over $A\llp t\rrp$, in other words,
\[
\mathrm{Ker}(H^1(A, G) \rightarrow H^1(A\llp t\rrp, G)) = \{ *\}.
\]
\end{corollary}

\bpf
Let $E$ be a $G$-torsor over $A$ that trivializes over $A\llp t\rrp$. By \cite{Hitchin-torsors}*{Corollary 2.1.22 (b) (with Example 2.1.18)}, this $E$ trivializes already over $A\{t\}[\frac 1t]$. Patching \cite{MB96}*{Theorem 5.5} then gives a $G$-torsor $\mathcal E$ over $\mathbb P^1_A$ that trivializes over $\mathbb P^1_A \setminus \{t = 0\}$ and satisfies $\mathcal E|_{\{t = 0\}} \simeq E$. At this point, \Cref{thm:swise} implies that $E$ is trivial.
\epf

When $A = k$, we recover the following result of Gille \cite{Gil24}*{Theorem 7.1}, which generalized the earlier \cite{FG21}*{Theorem 5.4}.
\begin{corollary}\label{cor:FG-B}
For a $k$-group scheme $G$ locally of finite type, we have
\[
H^1(k, G) \hookrightarrow H^1(k\llp t\rrp , G).
\]
\end{corollary}

\bpf
By twisting \cite{Gir71}*{chapitre III, remarque 2.6.3} (with \cite{SP}*{Lemmas \href{https://stacks.math.columbia.edu/tag/04SK}{04SK} and \href{https://stacks.math.columbia.edu/tag/0421}{0421}} and \S\ref{pp:k-groups-lft} for the representability of the resulting inner form of $G$), we only need to show that every $G$-torsor $E$ that trivializes over $k\llp t\rrp$ is trivial. However, if $E(k\llp t\rrp) \neq \emptyset$, then, by \S\ref{pp:Xgred}, our $E$ reduces to a $G^{\mathrm{gred}}$-torsor over $k$ that trivializes over $k\llp t\rrp$. Since $G^{\mathrm{gred}}$ is $k$-smooth, \Cref{cor:FG-A} then implies that $E$ is trivial, as desired.
\epf

\section{Generically trivial torsors under constant groups are Zariski locally trivial}

In \S\ref{sec:GS-proof}, we settle the Grothendieck--Serre question over an arbitrary base field, and then in \S\ref{sec:counter} we give various examples showing that our hypotheses are sharp. We conclude in \S\ref{sec:Iwasawa} by using our main result to establish the Iwasawa decomposition for arbitrary $k$-group schemes locally of finite type.

\csub[Grothendieck--Serre for torsors under constant groups]
\label{sec:GS-proof}

With the inputs of the previous chapters, we are ready to establish our main result in \Cref{thm:GS-main}. The final stepping stone is the following standard lemma, which encapsulates the geometric arguments that transform the Grothendieck--Serre problem to the study of torsors over the relative $\mathbb P^1$.

\blem\label{lem:GS-to-P1}
For a $k$-group scheme $G$ locally of finite type, a geometrically regular, semilocal $k$-algebra $R$, and a generically trivial $G$-torsor $E$ over $R$, there is a $G$-torsor $\mathcal{E}$ over $\mathbb{P}^1_R$ with $\mathcal{E}|_{\{t = 0\}} \simeq E$ such that $\mathcal{E}|_{\mathbb P^1_R \setminus Z}$ is trivial for some $R$-finite closed $Z\subset\mathbb{A}^1_R$.
\elem

\bpf
The argument is standard but  appears in the literature only for reductive $G$ (see \cite{FP15}, \cite{Pan20a}, \cite{split-unramified}, \cite{totally-isotropic}, among others), so we give it in full. \Cref{thm:pseudo-extend}~\ref{m:PE-i} reduces $E$ to a generically trivial $G^0$-torsor over $R$, so we lose no generality by assuming that $G$ is connected, so of finite type. A limit argument based on the Popescu theorem (see \S\ref{pp:conv}) then lets us assume that $R$ is the semilocal ring of a smooth, affine $k$-scheme $X$ and that $E$ spreads out to a $G$-torsor $\widetilde{E}$ over $X$ that trivializes away from some closed $Y \subset X$ of codimension $> 0$. By decomposing $X$ into connected components, we may assume that $X$ is integral of dimension $d > 0$. By the presentation lemma \cite{CTHK97}*{Theorem 3.1.1} (see also \cite{torsors-complement}*{Lemma 8.1}), at the cost of shrinking $X$ around $\mathrm{Spec}(R)$, there are an affine open $S \subset \mathbb A^{d - 1}_k$ and a smooth morphism $\pi \colon X \rightarrow S$ of relative dimension $1$ for which $Y$ is $S$-finite. By base changing along the map $\mathrm{Spec}(R) \rightarrow S$, we therefore obtain
\benumr
    \item a smooth, affine $R$-scheme $C$ of pure relative dimension $1$ (the base change of $X$);
    \item a section $s \in C(R)$ (obtained from the ``diagonal'' section $\mathrm{Spec}(R) \rightarrow X$);
    \item an $R$-finite closed subscheme $Z \subset C$ (the base change of $Y$);
    \item a $G$-torsor $\mathcal{E}$ over $C$ such that $s^*\mathcal{E} \cong E$ and $\mathcal{E}$ trivializes over $C \setminus Z$ (the base change of $\widetilde{E}$).
\end{enumerate}
By \cite{split-unramified}*{Lemmas 6.1 and 6.3}, we then reduce further to when  there is a quasi-finite, flat $R$-map $C\to \mathbb{A}^1_R$ that maps $Z$ isomorphically onto an $R$-finite closed subscheme $Z^\prime\subset\mathbb{A}^1_R$ for which $Z\cong Z^\prime\times_{\mathbb{A}^1_R}C$. By patching \cite{split-unramified}*{Lemma 7.1}, then $\mathcal E$ descends to a $G$-torsor over $\mathbb A^1_R$ that is trivial away from $Z^\prime$ and whose $s$-pullback is $E$, so we may assume that $C = \mathbb A^1_R$. We then change coordinates to make $s$ be $\{t = 0\}$ and patch $\mathcal E$ with the trivial torsor over $\mathbb P^1_R \setminus Z$ to conclude.
\epf

\begin{theorem}\label{thm:GS-main}
Let $G$ be a $k$-group scheme locally of finite type, let $R$ be a geometrically regular, semilocal $k$-algebra, and let $E$ be a generically trivial $G$-torsor over $R$. If either
\benumr
    \item\label{m:GSM-i} every $\overline{k}$-torus of $G_{\overline{k}}$ lies in $(G^{\mathrm{gred}})_{\overline{k}}$ \up{see Remark~\uref{rem:condition-star}}\uscolon or
    \item\label{m:GSM-ii} $E$ is \'etale locally trivial\uscolon
\end{enumerate}
then $E$ is trivial.
\end{theorem}

\bpf
Set $K \colonequals \mathrm{Frac}(R)$ and let $\overline{E} \subseteq E$ be the scheme-theoretic image of $(E_K)^{\mathrm{gred}} \subseteq E_K$. Since $K/k$ is separable, \S\ref{pp:Xgred} ensures that the $K$-subgroup $(G^{\mathrm{gred}})_K \le G_K$ agrees with $(G_K)^{\mathrm{gred}} \le G_K$, so that its action on $E$ preserves $(E_K)^{\mathrm{gred}}$. The formation of the scheme-theoretic image of a quasi-compact morphism commutes with flat base change (see \cite{SP}*{Lemma \href{https://stacks.math.columbia.edu/tag/089E}{089E}}), so $(E_K)^{\mathrm{gred}} \times_K (G^{\mathrm{gred}})_K$ is scheme-theoretically dense in $\overline{E} \times_R (G^{\mathrm{gred}})_R$. In particular, the action of $(G^{\mathrm{gred}})_R$ on $E$ preserves $\overline{E}$. If $E$ trivializes \'etale locally on $R$, then this commutation with flat base change also allows us to check \'etale locally on $R$ that $\overline{E}$ is even a $G^{\mathrm{gred}}$-torsor over $R$. Therefore, since every $K$-point of $E$ lies in $(E_K)^{\mathrm{gred}}$, we have reduced \ref{m:GSM-ii} to \ref{m:GSM-i}.

As for \ref{m:GSM-i}, \Cref{thm:pseudo-extend}~\ref{m:PE-i} reduces us to when $G = G^{\mathrm{sm},\, \mathrm{lin}}$, that is, when $G$ is connected, smooth, and affine. Moreover, \Cref{lem:GS-to-P1} supplies a $G$-torsor $\mathcal{E}$ over $\mathbb{P}^1_R$ with $\mathcal{E}|_{\{t = 0\}} \simeq E$ such that $\mathcal{E}|_{\mathbb P^1_R \setminus Z}$ is trivial for some $R$-finite closed $Z\subset\mathbb{A}^1_R$. Since $\mathcal{E}|_{\{ t = \infty\}}$ is trivial, \Cref{thm:swise} then implies that $E$ is also trivial, as desired.
\epf

\csub[Counterexamples to more general versions of the Grothendieck--Serre question]
\label{sec:counter}

We illustrate the sharpness of the assumptions of \Cref{thm:GS-main}: we give counterexamples to some overly optimistic generalizations of the Grothendieck--Serre question. In \Cref{eg:fiberwise-unipotent,eg:congruence}, we show that, unlike in the reductive case, we cannot allow $G$ to be defined merely over $R$ even if we assume that the group is smooth, affine, with connected fibers, in \Cref{eg:semilocal} we argue that a straightforward reduction to the local case is unlikely to exist, and in \Cref{eg:condition-star} we recall that we cannot drop the condition on the $\overline{k}$-tori of $G_{\overline{k}}$ lying in $(G^{\mathrm{gred}})_{\overline{k}}$. The connectedness of $R$-fibers is important in these examples because it shows that the failure is not merely bootstrapped in some way from the failure of \eqref{eqn:finite-extend} to hold for quasi-finite $S$-schemes that are not finite. It seems that prior to these examples it was not known whether the Grothendieck--Serre type triviality of generically trivial torsors holds for general smooth, affine groups with connected fibers over regular local rings, although Colliot-Th\'el\`ene--Sansuc \cite{CTS87}*{Example 5.9} gave a smooth, affine counterexample over $\mathbb{R}[x]_{(x)}$ whose generic fiber is a torus and special fiber is $\mathbb{G}_{a,\, \mathbb{R}}\times_{\mathbb{R}}\mu_{2,\, \mathbb{R}}$. 

\beg \label{eg:fiberwise-unipotent}
Suppose that $k$ is an imperfect, finitely generated field of characteristic $p > 0$ and let $c\in k\llb  t\rrb  ^\times$ be such that its class modulo $t$ is not a $p$-th power in $k$. As Gabber pointed out to us, $\{ x=x^p+cy^p+tz\} \le \mathbb{G}_{a,\, k\llb  t\rrb  }^3$ is an example of a smooth, affine $k\llb  t\rrb  $-group $G$ with generic fiber isomorphic to $\mathbb{G}_{a,\, k((t))}^2$ (solve for $z$) and special fiber isomorphic to the product of $\mathbb{G}_{a,\, k}$ and the $1$-dimensional connected, smooth, wound unipotent $k$-group $\{x = x^p + cy^p\} \le \mathbb{G}_{a,\, k}^2$. Consequently, every $G$-torsor over $k((t))$ is trivial, yet \cite{Ros24}*{Theorem 1.6} (with \cite{Hitchin-torsors}*{Theorem 2.1.6 (b)}) ensures that $G$ has infinitely many nontrivial torsors over $k\llb  t\rrb  $, so the triviality of generically trivial torsors fails for $G$ over $k\llb  t\rrb  $. The $k\llb  t\rrb  $-group $G$ has a fiberwise constant reductive rank (equal to $0$), but its split unipotent rank is not fiberwise constant. However, by instead considering the product of $G$ and a smooth, affine $k\llb  t\rrb  $-group $G^\prime$ whose generic fiber is connected, $1$-dimensional, wound unipotent and special fiber is $\mathbb{G}_{a,\,k}$ (for instance, $G^\prime$ could be $\{x = x^p + (1+ t)y^p\} \le \mathbb{G}_{a,\, k\llb  t\rrb  }^2$), we obtain the same failure of triviality of generically trivial torsors under a smooth, affine $k\llb  t\rrb  $-group whose fibers are connected, unipotent, and now have constant split unipotent ranks.
\eeg

\beg \label{eg:congruence}
Gabber suggested the following further counterexample for nonconstant groups in characteristic $0$ that is based on dilatations and congruence subgroups. 

Let $T$ be a torus over $\mathbb{Q}$ that is not retract rational in the sense that no nonempty open $U \subset T$ admits a factorization of $\mathrm{id}_U$ as $U \rightarrow \mathbb{A}^N_{\mathbb{Q}} \rightarrow U$, such a $T$ exists by \cite{Sca20}*{Theorems 1.1 and~1.3}. Concretely, by a result of Voskresenskii--Saltman, we may let $T$ be the kernel of the norm map $\mathrm{Res}_{K/\mathbb{Q}}(\mathbb{G}_m) \rightarrow \mathbb{G}_m$ for a $(\mathbb{Z}/2\mathbb{Z})^2$-extension $K/\mathbb{Q}$ (compare with \cite{Sal84}*{Theorem 3.14}). By \cite{Gil09}*{Proposition 5.1, proof of 2)$\Rightarrow$1)}, the failure of retract rationality of $T$ supplies an essentially smooth, local $\mathbb{Q}$-algebra $R$ and an $f \in R$ for which the map $T(R) \rightarrow T(R/fR)$ is not surjective. Let $G$ be the $R$-group that is the dilatation (also called a N\'eron blowup) of $T_R$ in the zero section of $T_{R/fR}$ (see \cite{MRR23}*{Definition 3.1}): concretely, by \cite{MRR23}*{Lemma 3.1, Theorem 3.2}, the $R$-group $G$ is smooth, affine, with connected $R$-fibers and on the category of those $R$-algebras $S$ in which $f$ is a nonzerodivisor, we have
$$
G(S) \cong \mathrm{Ker}(T(S)\rightarrow T(S/fS)).
$$
For an $\alpha \in T(R/fR)$, the dilatation $E_\alpha$ of $T$ in the image of the $(R/fR)$-point $\alpha$ is a $G$-torsor for the \'etale topology that is trivial if and only if $\alpha$ lifts to $T(R)$, in particular, $E_\alpha$ is not trivial for some $\alpha$. However, $E_\alpha$ is generically trivial, in fact, already $E_\alpha|_{R[\frac{1}{f}]}$ is trivial (note that $f \neq 0$ by our choice of $R$ and $f$). Thus, $G$ has nontrivial but generically trivial torsors over $R$. Moreover, all of these torsors are trivial over $R/fR$ because \cite{MRR23}*{Theorem 3.5 (1)} ensures that $G_{R/fR} \simeq \mathbb G_{a,\, R/fR}^{\mathrm{dim}(T)}$, in particular, this example is a Zariski rather than a Henselian phenomenon.
\eeg

\beg\label{eg:semilocal}
The method of dilatations suggests the following example of a nontrivial but Zariski locally trivial $G$-torsor over a semilocal, essentially smooth ring $R$. This rules out reductions from semilocal to local regular rings in the Grothendieck--Serre problem. Many further similar examples may be constructed by using \cite{Sca23}*{th\'eor\`eme 1.3}.

Let $t$ be the standard coordinate of $\mathbb P^1_{\mathbb F_2}$, let $R$ be the semilocal ring of $\mathbb G_{m,\, \mathbb F_2}$ at two distinct maximal ideals $\mathfrak m, \mathfrak m^\prime \subset \mathbb F_2[t, t^{-1}]$, and let $H \le \mathbb G^2_{a,\, R}$ be a smooth, affine $R$-subgroup cut out by the equation $x + x^4 = ty^4$. By \cite{Ros21a}*{proof of Lemma 4.2}, this $H$ has precisely two $\mathbb F_2(t)$-points: $(0, 0)$ and $(1, 0)$, which both extend to $R$-points. In particular, its $k_{\mathfrak m}$-point $(0, 0)$ lifts to an $R_{\mathfrak m}$-point, its $k_{\mathfrak m^\prime}$-point $(1, 0)$ lifts to an $R_{\mathfrak m^\prime}$-point, and together they assemble to a $(k_{\mathfrak m} \times k_{\mathfrak m^\prime})$-point $\alpha$ that does not lift to any $R$-point of $H$. Let $G$ be the smooth, affine $R$-group with connected $R$-fibers that is the dilatation of $H$ in the identity of $H_{k_{\mathfrak m} \times k_{\mathfrak m^\prime}}$. The dilatation of $H$ in the $(k_{\mathfrak m} \times k_{\mathfrak m^\prime})$-point $\alpha$ is a nontrivial $G$-torsor over $R$ that trivializes over both $R_{\mathfrak m}$ and $R_{\mathfrak m^\prime}$. 
\eeg

\beg \label{eg:condition-star}
In \cite{FG21}*{Section 7.2}, Florence and Gille gave examples of affine $k$-groups $G$ of finite type with nontrivial but generically trivial $G$-torsors over $k\llb  t\rrb  $, which shows that the condition on $\overline{k}$-tori in \Cref{thm:main-GS} cannot be dropped. We now present a further such example related to forms of pseudo-reductive groups.

By \cite{CP16}*{Proposition 6.2.2 and Example 6.2.3}, for any pseudo-semisimple $k$-group $\mathscr{H}$, the $k$-group $G\colonequals \underline{\mathrm{Aut}}_{\mathrm{gp}}(\mathscr{H})$ is affine, of finite type, and if $\mathscr H \cong \mathrm{Res}_{k^\prime/k}(H_{k^\prime})$ for a purely inseparable field extension $k^\prime/k$ of degree $p \colonequals \mathrm{char}(k)$ and a nontrivial, semisimple $k$-group $H$, then $G$ is also not smooth. We now show that, for any pseudo-semisimple $\mathscr{H}$ of this form, $G$ is not even a normal subgroup of a smooth $k$-group, in fact, some $\overline{k}$-torus of $G_{\overline{k}}$ does not lie in $(G^{\mathrm{gred}})_{\overline{k}}$ (see \Cref{rem:condition-star}). For this, by \Cref{thm:main-GS}, it suffices to exhibit a $k\llb  t\rrb  $-group form $\mathscr{H}'$of $\mathscr{H}_{k\llb  t\rrb  }$ for which $\mathscr{H}'_{k\llp t\rrp } \simeq \mathscr{H}_{k\llp t\rrp }$ but $\mathscr{H}' \not\simeq \mathscr{H}_{k\llb  t\rrb  }$. Our candidate is
$$
\mathscr H' := \mathrm{Res}_{(k\llb  t\rrb  [x]/(x^p - at^p\rrp /k\llb  t\rrb  }(H_{k\llb  t\rrb  [x]/(x^p - at^p)}), \quad \text{where} \quad k^\prime \cong k[y]/(y^p -a)
$$
for some $a \in k \setminus k^p$. Since $k'\llb  t\rrb  [x]/(x^p - at^p) \simeq k' \llb  t\rrb   \otimes_k k^\prime$ , we have $\mathscr{H}'_{k^\prime\llb  t\rrb  } \simeq \mathscr{H}_{k^\prime\llb  t\rrb  }$, so $\mathscr{H}'$ is a $k\llb  t\rrb  $-form of $\mathscr{H}_{k\llb  t\rrb  }$. Similarly, $k\llp t\rrp [x]/(x^p - at^p) \simeq k\llp t\rrp  \otimes_k k^\prime$, so $\mathscr{H}'_{k\llp t\rrp } \simeq \mathscr{H}_{k\llp t\rrp }$. Finally, $\mathscr{H}' \not\simeq \mathscr{H}_{k\llb  t\rrb  }$ since the $k$-group $\mathscr{H}'|_{t = 0} \cong \mathrm{Res}_{(k[x]/(x^p))/k}(H_{k[x]/(x^p)})$ is not pseudo-reductive (see \Cref{lem:res-unipotent}).
\eeg

\csub[The Iwasawa decomposition]\label{sec:Iwasawa}

We conclude by using the Grothendieck--Serre conclusion of \Cref{thm:GS-main} to prove the following Iwasawa decomposition theorem for arbitrary $k$-group schemes locally of finite type.

\begin{theorem}[Iwasawa Decomposition]\label{thm:Iwasawa}
For a discrete valuation ring $\mathcal{O}$ that is a geometrically regular $k$-algebra, $K \colonequals \mathrm{Frac}(\mathcal{O})$, a $k$-group scheme $G$ locally of finite type, and a pseudo-parabolic $k$-subgroup $P \le G^{\mathrm{sm},\, \mathrm{lin}}$ of the smooth linear part of $G$, we have
\[
G(K) = P(K) G(\mathcal{O}), \quad \text{in particular,} \quad G(K) = G^{\mathrm{sm},\, \mathrm{lin}}(K) G(\mathcal{O}).
\]
\end{theorem}

\bpf
The pseudo-properness of the connected components of $G/P$ supplied by \Cref{cor:GP-pseudo-proper} gives $(G/P)(\mathcal O) \cong (G/P)(K)$. Since, by \Cref{thm:GS-main} (and \S\ref{pp:pseudo-parabolic}), generically trivial $P$-torsors over $\mathcal{O}$ are trivial, we get the sought $G(K) = P(K) G(\mathcal{O})$.
\epf

\beg
 In the setting of \Cref{thm:Iwasawa}, suppose that $G$ is pseudo-reductive with a split maximal $k$-torus $S \le G$ and that the residue field of $\mathcal O$ is separable over $k$. Let $B \le G$ be a pseudo-Borel $k$-subgroup containing $S$ and let $U \le B$ be its (split) unipotent $k$-radical (see \S\ref{pp:pseudo-parabolic}). By \cite{CGP15}*{beginning of the proof of Theorem C.1.9} and \S\ref{pp:Gm-actions}, we have $B \cong U \rtimes Z_G(S)$ with $Z_G(S)/S$ wound unipotent. Thus, \Cref{prop:PP-groups}~\ref{m:PPG-i} ensures that $Z_G(S)(K) = X_*(S)Z_G(S)(\mathcal O)$ where a cocharacter $\lambda \colon \mathbb G_{m,\, k} \rightarrow S$ in $X_*(S)$ is interpreted to give the $K$-point $\lambda(\pi)$ for a fixed uniformizer $\pi \in \mathcal O$. Thus, overall, in this case the Iwasawa decomposition takes the form
$$
G(K) = U(K)X_*(S)G(\mathcal O).
$$
\eeg



\end{document}